\documentclass{svproc}
\usepackage{url}

\usepackage[linesnumbered,ruled]{algorithm2e}
\usepackage{amsmath} % assumes amsmath package installed
\usepackage{amssymb}  % assumes amsmath package installed
\usepackage{bm}
\usepackage{cite}
\usepackage{color}
\usepackage{comment}

\usepackage[font=scriptsize]{caption}
\usepackage{graphicx} % for pdf, bitmapped graphics files
\usepackage[font=scriptsize]{subfig}
\usepackage{floatrow} %to position captions
\usepackage{siunitx}
\DeclareMathOperator*{\argmin}{arg\,min}
\usepackage[resetlabels]{multibib}
\newcites{S}{References}
\newcommand{\M}[1]{\ensuremath{\textbf{\text{#1}}}} % Matrices
\newcommand{\V}[1]{\ensuremath{\textbf{\text{#1}}}} % Vectors
\newcommand{\Exp}[2]{\ensuremath{\mathop{\mathbb{E}}_{\substack{{#1}}} \left[{#2}\right]}} % Expectation
\newcommand{\Var}{\ensuremath{\operatorname{Var}}} % Variance
\newcommand{\tr}[1]{{#1}^{\ensuremath{\mathsf{T}}}} % transpose 
\newcommand{\inv}[1]{{#1}^{\ensuremath{\!\!\mathsf{-1}\,}}} % inverse 
\newcommand{\Th}{\ensuremath{{}^{\textrm{th}}}}

\newcommand{\mxx}[1]{{\color{blue}#1\ }}  % Edits

\begin{document}
\mainmatter              % start of a contribution
%
% Previous title:
% Decoupling stochastic optimal control problems for efficient solution: insights from experiments across a wide range of noise regimes
% Suman's title:
%\title{On the Tractable Search for Feedback in Robotic Planning under Uncertainty: Insights from Experiments over a wide noise regime}
%
%a That sounds like it is a search for feedback, where that search is tractable.
% Dylan's title:
\title{Experiments with Tractable Feedback in Robotic Planning under Uncertainty: Insights over a wide range of noise regimes (Extended Report)} 
\titlerunning{Experiments with Tractable Feedback under Uncertainty}  % abbreviated title (for running head)
%\toctitle{Decoupling stochastic optimal control problems for efficient solution: insights from experiments across a wide range of noise regimes}
\author{Mohamed Naveed Gul Mohamed \and Suman Chakravorty \and
Dylan A. Shell}
\authorrunning{M.N. Gul Mohamed et al.} % abbreviated author list (for running head)
\tocauthor{Mohamed Naveed Gul Mohamed, Suman Chakravorty and
Dylan A. Shell}
\institute{Texas A\&M University, College Station TX 77843, USA,\\
\email{mohdnaveed96@gmail.com, schakrav@tamu.edu, dshell@tamu.edu}.}

\maketitle              % typeset the title of the contribution

\begin{abstract}
We consider the problem of robotic planning under uncertainty. 
This problem may be posed as a stochastic optimal control problem, complete solution to which is fundamentally intractable owing to the infamous curse of dimensionality. 
We report the results of an extensive simulation study in which we have compared two methods, both of which aim to salvage tractability by using alternative, albeit inexact, means for treating feedback.
The first is a recently proposed method based on a near-optimal ``decoupling principle'' for tractable feedback design, wherein a nominal open-loop problem is solved, followed by a linear feedback design around the open-loop. 
The second is Model Predictive Control (MPC), a widely-employed method that uses repeated re-computation of the nominal open-loop problem during execution to correct for noise, though when interpreted as feedback, this can only said to be an implicit form.
We examine a much wider range of noise levels than have been previously reported and empirical evidence suggests that the decoupling method allows for tractable planning over a wide range of uncertainty conditions without unduly sacrificing performance.
\end{abstract}
\section{Introduction}

Planning under uncertainty is a central problem in robotics.  The space of
current methods includes several contenders, each with different simplifying
assumptions, approximations, and domains of applicability.  This is a natural
consequence of the fact that the challenge of dealing with the continuous
state, control and observation space problems, for non-linear systems and
across long-time horizons with significant noise, and potentially multiple
agents, is fundamentally intractable.

% MPC is one means for tackling these problems. It addresses uncertainty
% through repeated computation.  Feedback arises implicitly through
% computations that repeatedly solve over a time horizon; prohibitive
% computational cost, and also limited theoretical statements about
% performance, means it falls short of an ideal solution.
% 
 
Model Predictive Control is one popular means for tackling optimal control
problems~\cite{Mayne_1,Mayne_2}.  The MPC approach solves a finite horizon
``deterministic'' optimal control problem at every time step given the current
state of the process, performs only the first control action and then repeats
the planning process at the next time step. In terms of computation, this can
be a costly endeavor and, when a stochastic control problem is well
approximated by the deterministic problem (when the noise is meager), much of
this computation is superfluous.

In this paper we consider the generalization of a recently proposed method~\cite{D2C1.0} that uses a
local feedback to control noise induced deviations from the deterministic (that
we term the \emph{nominal}) trajectory.  When the deviation is too large for the
feedback to manage, replanning is triggered and it computes a fresh nominal.
Otherwise, the feedback tames the perturbations during execution and no
computation is expended in replanning.  Put another way, the method decouples
feedback and planning/nominal control but will fall back to replanning when
perturbations are excessive. Thus, by considering every deviation to
necessitate replanning, this approach will essentially reduce to MPC itself.

\begin{figure}[h]
    \centering
    \subfloat[A single agent.]{        \includegraphics[width=.25\textwidth]{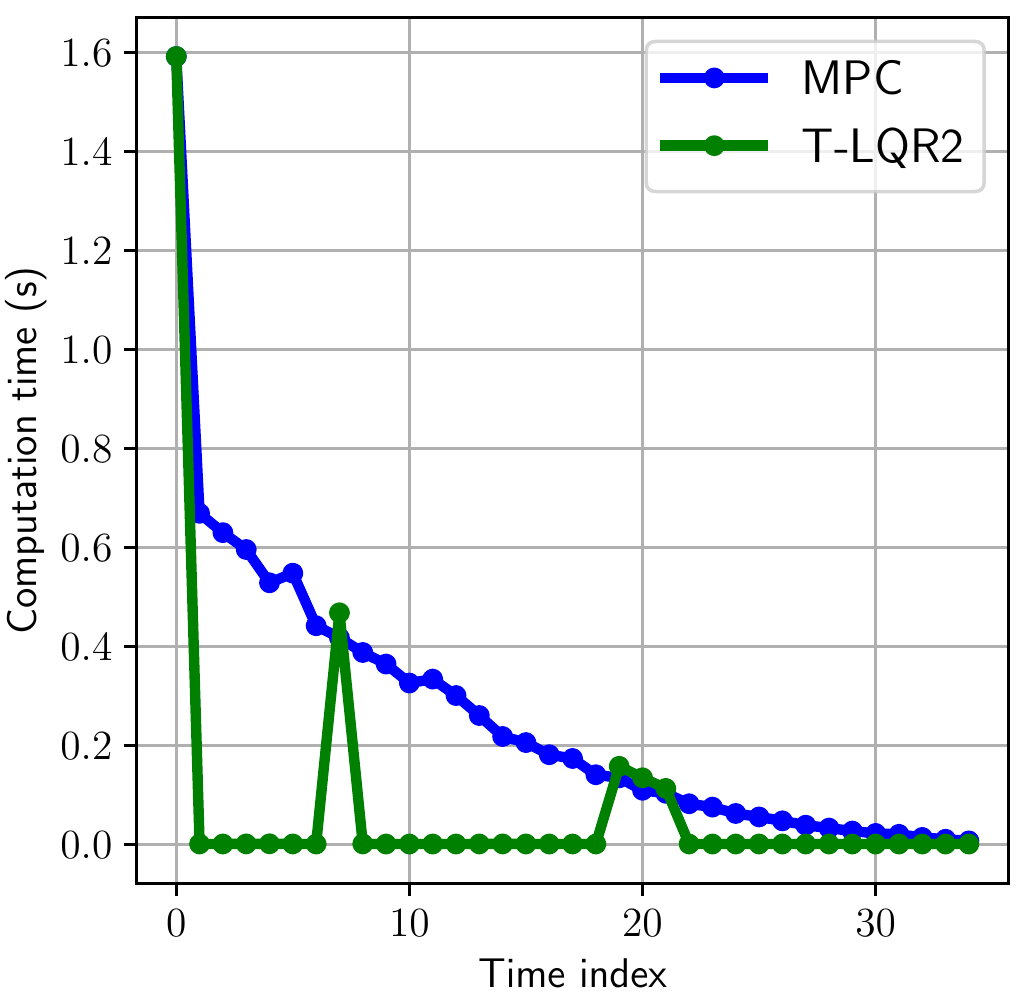}
    \label{1_agent_replan_time}}
    \subfloat[Three agents.]{
    \includegraphics[width=.25\textwidth]{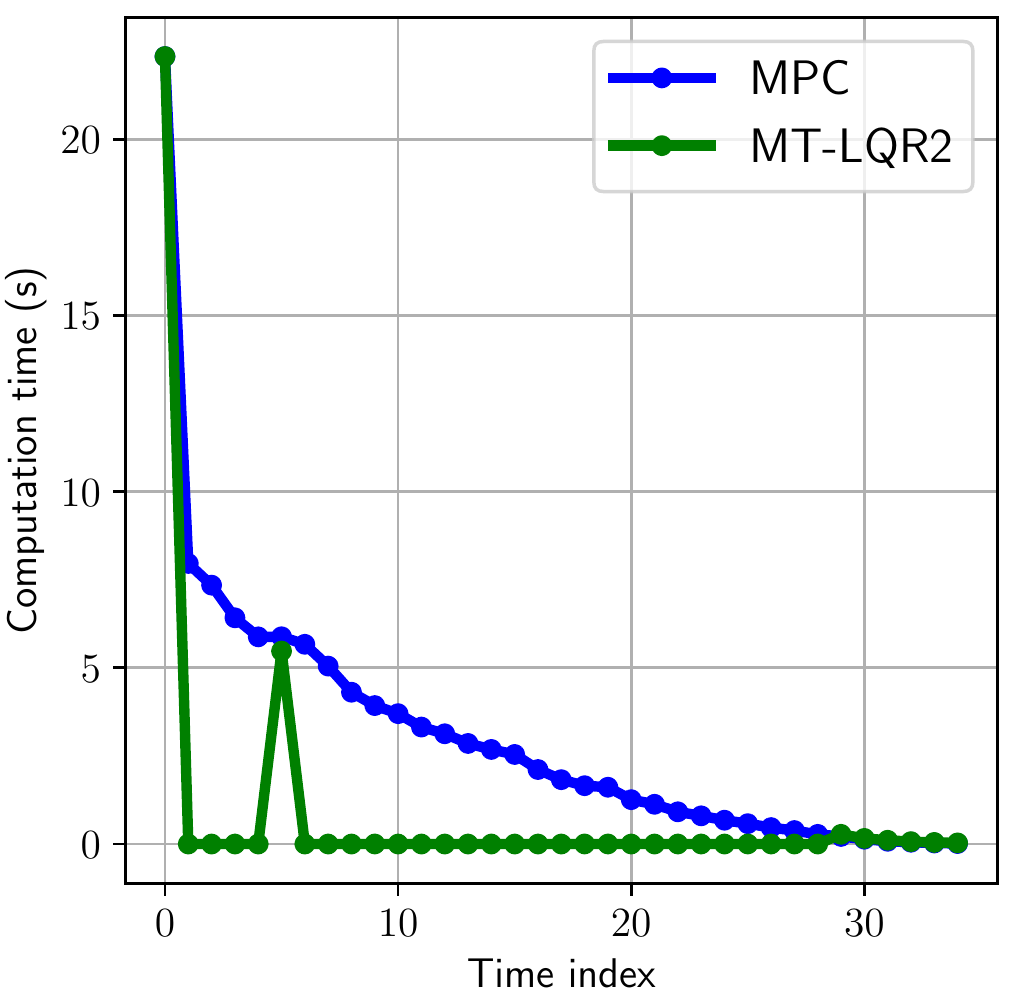}
    \label{3_agent_replan_time}}
\caption{Computation time expended by MPC (in blue) and the decoupling algorithms described herein (in green), at each time step for a sample experiment involving navigation. Both cases result in nearly identical motions by the robot.The peaks in T-LQR2 and MT-LQR2 happen only when replanning takes place. Computational effort decreases for both methods because the horizon diminishes as the agent(s) reach their goals.(To relate to subsequent figures: noise parameter $\epsilon = 0.4$ and the replan threshold = 2\% of cost deviation.)}
\label{fig:timeplot}
\end{figure}

We present an empirical investigation of this decoupling approach, exploring
dimensions that are important in characterizing its performance---key among
these being the triggering of replanning. Hence, the primary focus of the study
is on understanding the performance across a wide range of noise conditions with comparison to the ``gold standard" of MPC.
Figure~\ref{fig:timeplot} gives an overall summary of the paper's findings: the
areas under the respective curves give the total computational resources
consumed---the savings by the decoupling method over MPC are seen to be
considerable.%
\vspace{-15pt}
\vspace{-20pt}
\vspace{-20pt}
\subsection{Related Work}
Robotic planning problems under uncertainty can be posed as a stochastic optimal control problem that requires the solution of an associated Dynamic Programming (DP) problem, however, as the state dimension  increases, the computational complexity
goes up exponentially \cite{bertsekas1}, Bellman's infamous ``curse of dimensionality".  There has been recent success using sophisticated (Deep) Reinforcement Learning (RL) paradigm to solve DP problems, where deep neural networks are used as the function approximators \cite{RLHD1, RLHD2,RLHD3, RLHD4, RLHD5}, however, the training time required for these approaches is still prohibitive to permit real-time robotic planning that is considered here. 

In the case of continuous state, control and observation space problems, the
Model Predictive Control \cite{Mayne_1, Mayne_2} approach has been used with a
lot of success in the control system and robotics community.   For deterministic
systems, the process results in solving the original DP problem in a recursive
online fashion. However, stochastic control problems, and the control of
uncertain systems in general, is still an unresolved problem in MPC. As
succinctly noted in \cite{Mayne_1}, the problem arises due to
the fact that in stochastic control problems, the MPC optimization at every
time step cannot be over deterministic control sequences, but rather has to be
over feedback policies, which is, in general, difficult to
accomplish since a tractable parametrization of such policies to perform the optimization over, is, in general, unavailable. Thus, the tube-based MPC approach, and its stochastic counterparts,
typically consider linear systems \cite{T-MPC1, T-MPC2,T-MPC3} for which a
linear parametrization of the feedback policy suffices but the methods become intractable when dealing with nonlinear systems. In recent work, we have introduced a ``decoupling principle'' that allows us to tractably solve such stochastic optimal control problems in a near optimal fashion, with applications to highly efficient RL and MPC implementations \cite{D2C1.0,T-PFC}. However, this prior work required a small noise assumption. In this work, we relax this small noise assumption to show, via extensive empirical evaluation, that even when the noise is not small, a \emph{replan-when-necessary} modification of the decoupled planning approach, akin to event-triggered MPC \cite{ETMPC1, ETMPC2}, suffices to keep the planning computationally efficient while retaining performance comparable to MPC. We note that event-triggered MPC inherits the same issues mentioned above with respect to the stochastic control problem, and consequently, the techniques are only tractable for linear systems. \textit{Lest it seem that we are being unduly critical of MPC, that is definitely not our intention: we believe that MPC type replanning is unavoidable in uncertain systems, instead we additionally believe that such replanning can be substantially reduced utilizing decoupling while tractably and rigorously extending MPC to stochastic systems, i.e., the decoupled approach is not competition, but rather complimentary, to MPC. Please also see the deeper historical context to this discussion at the end of Section 3.1, after we have presented the basic near-optimality result.}

The problem of multiple agents further and severely compounds the planning
problem since now we are also faced with the issue of a control space that
grows exponentially with the number of agents in the system. Moreover, since
the individual agents never have full information regarding the system state,
the observations are partial. Furthermore, the decision making has to be done
in a distributed fashion which places additional constraints on the networking
and communication resources.  In a multi-agent setting, the stochastic optimal
problem can be formulated in the space of joint policies.  Some variations of
this problem have been successfully characterized and tackled based on the
level of observability, in/dependence of the dynamics, cost functions and
communications
\cite{seuken2008formal,oliehoek2016concise,pynadath2002communicative}. This has
resulted in a variety of solutions from fully-centralized
\cite{boutilier1996planning} to fully-decentralized approaches with many
different subclasses \cite{amato2013decentralizedB,oliehoek2012decentralized}.

The major concerns of the
multi-agent problem are tractability of the solution and the level of
communication required during the execution of the policies. In this paper, we also consider a generalization of the decoupling principle to a multi-agent, fully observed setting. We show that this leads to a spatial decoupling between agents in that they do not need to communicate for long periods of time during execution. Albeit, we do not consider the problem of when and how to replan in this paper, assuming that there exists a (yet to be determined) distributed mechanism that can achieve this, we nonetheless show that there is a highly significant increase in planning efficiency over a wide range of noise levels.  
%
\begin{comment}
\subsection{Outline of Paper}
The rest of the document is organised as follows: Section~\ref{section:prob} states the problem, Section~\ref{section:decoupling} gives background on the decoupling principle, Section~\ref{section:algorithms} explains the planning algorithms used, Section~\ref{section:results} discusses the results and observations and Section~\ref{section:conclusion} concludes.  
\end{comment}
%
\section{Problem Formulation}
\label{section:prob}
The problem of robot planning and control under noise can be formulated as a stochastic optimal control problem in the space of feedback policies. We assume here that the map of the environment is known and state of the robot is fully observed. Uncertainty in the problem lies in the system's actions. %
\subsection{System Model:}
For a dynamic system, we denote the state and control vectors by $\V{x}_t \in \ \mathbb{X} \subset \ \mathbb{R}^{n_x}$ and $\V{u}_t \in \ \mathbb{U} \subset \ \mathbb{R}^{n_u}$ respectively at time $t$. The motion model $f : \mathbb{X} \times \mathbb{U} \times \mathbb{R}^{n_u}   \rightarrow \mathbb{X} $ is given by the equation 
\begin{equation}
    \V{x}_{t+1}= f(\V{x}_t, \V{u}_t, \epsilon\V{w}_t); \  \V{w}_t \sim \mathcal{N}(\V{0}, {\mathbf \Sigma}_{\V{w}_t}) 
    \label{eq:model},
\end{equation}
where \{$\V{w}_t$\} are zero mean independent, identically distributed (i.i.d) random sequences with variance ${\mathbf\Sigma}_{\V{w}_t}$, and $\epsilon$ is a small parameter modulating the noise input to the system. 
\subsection{Stochastic optimal control problem:} %
The stochastic optimal control problem for a dynamic system with initial state $\V{x}_0$ is defined as:
\begin{equation}
    J_{\pi^{*}}(\V{x}_0) = \min_{\pi} \ \Exp{}{\sum^{T-1}_{t=0} c(\V{x}_t, \pi_t (\V{x}_t)) + c_T(\V{x}_T)},
\end{equation}
\begin{equation}
   s.t.\ \V{x}_{t+1} = f(\V{x}_t, \pi_t (\V{x}_t), \epsilon\V{w}_t),
\end{equation}
where:
\begin{description}
\item~
\begin{itemize}
    \item the optimization is over feedback policies $\pi := \{ \pi_0, \pi_1, \ldots, \pi_{T-1} \} $ and $\pi_t(\cdot)$: $\mathbb{X} \rightarrow \mathbb{U}$ specifies an action given the state, $\V{u}_t = \pi_t(\V{x}_t)$;
    \item $J_{\pi^{*}}(\cdot): \mathbb{X} \rightarrow \mathbb{R}$  is the cost function on executing the optimal policy $\pi^{*}$; 
    \item $c_t(\cdot,\cdot): \mathbb{X} \times \mathbb{U} \rightarrow \mathbb{R} $  is the one-step cost function;
    \item $c_T(\cdot): \mathbb{X} \rightarrow \mathbb{R}$ is the terminal cost function;
    \item $T$ is the horizon of the problem;
    \item the expectation is taken over the random variable $\V{w}_t$.
\end{itemize}
\end{description}

\section{A Decoupling Principle}
\label{section:decoupling}
Now, we give a brief overview of a ``decoupling principle'' that allows us to substantially reduce the complexity of the stochastic planning problem given that the parameter $\epsilon$ is small enough. We only provide an outline here and the relevant details can be found in our recent work \cite{D2C1.0}. We shall also present a generalization to a class of multi-robot problems. Finally, we preview the results in the rest of the paper.
\subsection{Near-Optimal Decoupling in Stochastic Optimal Control}
Let $\pi_t(\V{x}_t)$ denote a control policy for the stochastic planning problem above, not necessarily the optimal policy. Consider now the control actions of the policy when the noise to the system is uniformly zero, and let us denote the resulting ``nominal'' trajectory and controls as $\overline{\V{x}}_t$ and $\overline{\V{u}}_t$ respectively, i.e., $\overline{\V{x}}_{t+1} = f(\overline{\V{x}}_t, \overline{\V{u}}_t, 0)$, where $\overline{\V{u}}_t = \pi_t(\overline{\V{x}}_t)$. Note that this nominal system is well defined. \\
Further, let us assume that the closed-loop (i.e., with $\V{u}_t = \pi_t(\V{x}_t)$), system equations, and the feedback law are smooth enough that we can expand the feedback law about the nominal as $\pi_t(\V{x}_t) = \overline{\V{u}}_t + \M{K}_t\delta \V{x}_t + \M{R}_t^{\pi}(\delta \V{x}_t)$, where $\delta \V{x}_t = \V{x}_t - \overline{\V{x}}_t$, i.e., the perturbation from the nominal, $\M{K}_t$ is the linear gain obtained by the Taylor expansion about the nominal in terms of the perturbation $\delta \V{x}_t$, and $\M{R}_t^{\pi}(\cdot)$ represents the second and higher order terms in the expansion of the feedback law about the nominal trajectory. Further we assume that the closed-loop perturbation state can be expanded about the nominal as: $\delta \V{x}_t = \M{A}_t \delta \V{x}_t + \M{B}_t \M{K}_t \delta \V{x}_t + \M{R}_t^f (\delta \V{x}_t) + \epsilon \M{B}_t \V{w}_t$, where the $\M{A}_t$, $\M{B}_t$ are the system matrices obtained by linearizing the system state equations about the nominal state and control, while $\M{R}_t^f(\cdot)$ represents the second and higher order terms in the closed-loop dynamics in terms of the state perturbation $\delta \V{x}_t$. Moreover, let the nominal cost be given by $\overline{J}^{\pi} = \sum_{t=0}^T \overline{c}_t$, where $\overline{c}_t = c(\overline{\V{x}}_t,\overline{\V{u}}_t)$, for $t\leq T-1$, and $\overline{c}_T = c_T(\overline{\V{x}}_T,\overline{\V{u}}_T)$. Further, assume that the cost function is smooth enough that it permits the expansion $J^{\pi} = \overline{J} + \sum_t \M{C}_t \delta \V{x}_t + \sum_t \M{R}_t^c(\delta \V{x}_t)$ about the nominal trajectory, where $\M{C}_t$ denotes the linear term in the perturbation expansion and $\M{R}_t^c(\cdot)$ denote the second and higher order terms in the same. Finally, define the exactly linear perturbation system $\delta \V{x}_{t+1}^\ell = \M{A}_t \delta \V{x}_t^\ell + \M{B}_t\M{K}_t \delta \V{x}_t^\ell + \epsilon \M{B}_t \V{w}_t$. Further, let $\delta J_1^{\pi,\ell}$ denote the cost perturbation due to solely the linear system, i.e., $\delta J_1^{\pi,\ell} = \sum_t \M{C}_t \delta \V{x}_t^\ell$.  Then, the decoupling result states the following \cite{D2C1.0}:
\begin{theorem}
The closed-loop cost function $J^{\pi}$can be expanded as $J^{\pi} = \overline{J}^{\pi} + \delta J_1^{\pi,\ell} + \delta J_2^{\pi}$. Furthermore, $\Exp{}{J^{\pi}} = \overline{J}^{\pi} + O(\epsilon^2)$, and $\Var[J^{\pi}] = \Var[\delta J_1^{\pi,\ell}] + O(\epsilon^4)$, where $\Var[\delta J_1^{\pi,\ell}]$ is $O(\epsilon^2)$.
\end{theorem}

Thus, the above result suggest that the mean value of the cost is determined almost solely by the nominal control actions while the variance of the cost is almost solely determined by the linear closed-loop system. Thus, the decoupling result says that the feedback law design can be decoupled into an open-loop and a closed-loop problem.\\
\textit{Open-Loop Problem:} This problem solves the deterministic/ nominal optimal control problem:
\begin{equation}
    \overline{J}= \min_{\overline{\V{u}}_t} \sum_{t=0} ^{T-1} c(\overline{\V{x}}_t,\overline{\V{u}}_t) + c_T(\overline{\V{x}}_T),
\label{DOCP}
\end{equation}
subject to the nominal dynamics: $\overline{\V{x}}_{t+1} = f(\overline{\V{x}}_t, \overline{\V{u}}_t)$. \\
\textit{Closed-Loop Problem:} One may try to optimize the variance of the linear closed-loop system
\begin{equation}
    \min_{\M{K}_t} \Var[\delta J_1^{\pi,\ell}]
\end{equation}
subject to the linear dynamics $\delta \V{x}_{t+1}^\ell = \M{A}_t \delta \V{x}_t^\ell + \M{B}_t \M{K}_t \delta \V{x}_t^\ell + \epsilon \M{B}_t \V{w}_t$.
However, the above problem does not have a standard solution but note that we are only interested in a good variance for the cost function and not the optimal one. Thus, this may be accomplished by a surrogate LQR problem that provides a good linear variance as follows.\\
\textit{Surrogate LQR Problem:} Here, we optimize the standard LQR cost:
% \begin{equation}
% \delta J_{\textsc{lqr}} = \min_{\V{u}_t} \sum_{t=0}^{T-1} \tr{\V{x}}_t\M{Q}\V{x}_t + \tr{\V{u}}_t\M{R}\V{u}_t + \tr{\V{x}}_T\M{Q}_f \V{x}_T,
% \end{equation}
%\textcolor{orange}{The version above was the original; I believe that it is supposed to over an expectation of the noise (it is a linear-quadratic stochastic control problem), and using $\delta$ variables (because it is around the nominal), as written next. If this is wrong, I think I need it explained to me.}
\begin{equation}
    \delta J_{\textsc{lqr}} =\min_{\V{u}_t} \Exp{\V{w}_t}{\sum_{t=0}^{T-1} \delta {\tr{\V{x}}_t} \M{Q} \delta \V{x}_t + \delta \tr{\V{u}}_t\M{R}\delta \V{u}_t + \delta \tr{\V{x}}_T \M{Q}_f \delta \V{x}_T},
    \label{LQRcost}
\end{equation}
subject to the linear dynamics $\delta \V{x}_{t+1}^\ell = \M{A}_t \delta \V{x}_t^\ell + \M{B}_t \delta \V{u}_t + \epsilon \M{B}_t \V{w}_t$. In this paper, this decoupled design shall henceforth be called the trajectory-optimized LQR (T-LQR) design.  \\

%\begin{remark}
\textbf{A Historical Context.} The above decoupled design might seem like a perturbation feedback design outlined in classical optimal control texts such as (Ch.~6, \cite{bryson}) and we are certainly not claiming that we are the first to discover it. However, the perturbation design was always thought to be heuristic and its ``goodness'' for the stochastic optimal control problem was essentially unexplored.  Notable as an exception is the reference \cite{fleming1971stochastic} that considers the problem of how good the deterministic feedback law is for the stochastic system, which is shown to be $O(\epsilon^4)$. However, that paper assumes the availability of the optimal deterministic feedback law which is the solution of the deterministic Hamilton-Jacobi-Bellman (HJB) equation (the DP equation in continuous time problems), which, in itself, is intractable as noted by Fleming as the ``practical difficulty'' in this work (pgs.~475--476 of \cite{fleming1971stochastic}). However, MPC, by repeatedly solving the deterministic optimal control problem at every time step, implicitly furnishes the \textit{deterministic feedback law}, and thus, offers the solution to the \textit{practical dilemma} above. The field of MPC, of course, was developed almost two decades after Fleming's work, while stochastic MPC/ MPC-under-uncertainty was explored only starting at the turn of millennium \cite{Mayne_1}. Thus, this connection was lost and never really explored in the MPC literature. This connection is critical if we want to \textit{tractably} extend MPC to stochastic systems in a \textit{theoretically justifiable fashion}, in the sense that in much of the stochastic MPC literature, these two aspects are at cross purposes to each other thereby preventing a satisfactory resolution. Thus, the MPC replanning logic is well justified theoretically, even when applied to a stochastic system. 

In fact, with a few further developments, and adaptation of Fleming's work to discrete time finite horizon problems, and if the linear feedback gain is modified suitably, the perturbation design also becomes $O(\epsilon^4)$ near-optimal. Due to paucity of space, we postpone this result to a future paper, however, for the sake of completeness and the reader's benefit, the result is included in the supplementary document. \textit{The ultimate takeaway is that the implicit MPC feedback law is an excellent approximation to the optimal stochastic policy, however, a T-LQR type perturbation feedback design is much cheaper computationally, while retaining identical near-optimality guarantees as MPC.}
%\end{remark}
\subsection{Multi-agent setting}
Now, we generalize the above result to a class of multi-agent problems. 
We consider a set of agents that are transition independent, i.e, their dynamics are independent of each other. For simplicity, we also assume that the agents have perfect state measurements. Let the system equations for the agents be given by:
$
\V{x}_{t+1}^j = f(\V{x}_t^j) + \M{B}^j_t(\V{u}_t^j + \epsilon \V{w}_t^j),
$
where $j = 1,2,\dots,M$ denotes the $j\Th$ agent. (We have assumed the control affine dynamics for simplicity). Further, let us assume that we are interested in the minimization of the joint cost of the agents given by $\mathcal{J} = \sum_{t=0}^{T-1} c(\M{X}_t,\M{U}_t) + \Phi(\M{X}_T)$, where $\M{X}_t = [\V{x}_t^1,\dots,\V{x}_t^M]$, and $\M{U}_t = [\V{u}_t^1,\dots, \V{u}_t^M]$  are the joint state and control action of the system. The objective of the multi-agent problem is minimize the expected value of the cost $\Exp{}{\mathcal{J}}$ over the joint feedback policy $\M{U}_t(\cdot)$. The decoupling result holds here too and thus the multi-agent planning problem can be separated into  an open and closed-loop problem. The open-loop problem consists of optimizing the joint nominal cost of the agents subject to the individual dynamics.\\
\textit{Multi-Agent Open-Loop Problem:}
\begin{align}\label{OL-MA}
\overline{\mathcal{J}} = \min_{\overline{\M{U}}_t} \sum_{t=0}^{T-1} c(\overline{\M{X}}_t,\overline{\M{U}}_t) + \Phi(\overline{\M{X}}_T), %\nonumber\\
\end{align}
subject to the nominal agent dynamics
$
\overline{\V{x}}_{t+1}^j = f(\overline{\V{x}}_t^j) + \M{B}^j_t\overline{\V{u}}_t^j.
$
The closed-loop, in general, consists of optimizing the variance of the cost $\mathcal{J}$, given by $\Var[\delta \mathcal{J}^\ell_1]$, where $\delta \mathcal{J}_1^\ell = \sum_t \M{C}_t \delta \M{X}_t^l$ for suitably defined $\M{C}_t$, and $\delta \V{X}_t^\ell = [\delta \V{x}_t^1,\dots, \delta \V{x}_t^M]$, where the perturbations $\delta \V{x}_t^j$ of the $j\Th$ agent's state  is governed by the decoupled linear multi-agent system $\delta \V{x}_t^j = \M{A}_t\delta \V{x}_t^j + \M{B}_t^j \delta \V{u}^j_t + \epsilon \M{B}_t^j \V{w}_t^j.$ This design problem does not have a standard solution but recall that we are not really interested in obtaining the optimal closed-loop variance, but rather a good variance. Thus, we can instead solve a surrogate LQR problem given the cost function $\delta \mathcal{J}_{\textsc{mtlqr}} = \sum_{t=0}^{T-1} \sum_j \delta {\tr{\V{x}_t^j}} \M{Q}^j \delta \V{x}_t^j + \delta \tr{\V{u}_t^j}\M{R}\delta \V{u}_t^j + \sum_j\delta \tr{\V{x}_T^j} \M{Q}^j_f \delta \V{x}_T^j$. Since the cost function itself is decoupled, the surrogate LQR design degenerates into a decoupled LQR design for each agent.\\
\textit{Surrogate Decoupled LQR Problem:}
\begin{equation*}
    \delta \mathcal{J}^j =\min_{\V{u}_t^j} \Exp{\V{w}_t^j}{\sum_{t=0}^{T-1} \delta {\tr{\V{x}_t^j}} \M{Q}^j \delta \V{x}_t^j + \delta \tr{\V{u}_t^j}\M{R}\delta \V{u}_t^j + \delta \tr{\V{x}_T^j} \M{Q}^j_f \delta \V{x}_T^j}, \text{subject to the}
\end{equation*}
linear decoupled agent dynamics \mbox{$\delta \V{x}_t^j = \M{A}_t\delta \V{x}_t^j + \M{B}_t^j \delta \V{u}^j_t + \epsilon \M{B}_t^j \V{w}_t^j.$}\\
\begin{remark}
Note that the above decoupled feedback design results in a spatial decoupling between the agents in the sense that, at least in the small noise regime, after their ``initial joint plan" is made, the agents never need to communicate with each other in order to complete their missions. However, note that the joint plan requires communication.
\end{remark}

\subsection{Planning Complexity versus Uncertainty}
The decoupling principle outlined above shows that the complexity of planning can be drastically reduced while still retaining near optimal performance for sufficiently small noise (i.e., parameter  $\epsilon \ll 1$).  Nonetheless, the skeptical reader might argue that this result holds only for low values of $\epsilon$ and thus, its applicability for higher noise levels is suspect.  Still, because the result is second order, it hints that near optimality might be over a reasonably large $\epsilon$.  Naturally, the question is \textsl{`will it hold for medium to higher levels of noise?'} We purposely leave the terminology of medium to high noise nebulous but what we mean shall become clear from our experiments.\\

\textit{Preview of the Results.} In this paper, we illustrate the
degree to which the above result still holds when we allow periodic replanning
of the nominal trajectory in T-LQR in an event triggered fashion, dubbed
T-LQR2.  Here, we shall use MPC as a ``gold standard'' for comparison since the
true stochastic control problem is intractable, and owing to Fleming's result \cite{fleming1971stochastic}, the MPC policy is
$O(\epsilon^4)$ near-optimal when compared to the true stochastic policy. In fact, we can make an identical strong $O(\epsilon^4)$ claim for T-LQR as well if the \textit{linear feedback gain is designed carefully}, but owing to the paucity of space, testing with this careful feedback design is left to a future paper. 
Here, we show
that though the number of replanning operations in T-LQR2 increases the
planning burden over T-LQR, it is still much reduced when compared to MPC,
which replans continually.  The ability to trigger replanning means that T-LQR2
can always produce solutions with the same quality as MPC, albeit by demanding
the same computational cost as MPC in instances when replanning is triggered. But for moderate levels
of noise, T-LQR2 can produce comparable quality output to MPC with substantial
computational savings.

In the high noise regime, replanning is more frequent but we
shall see that there is another consideration at play. Namely, that the
effective planning horizon decreases and there seems no benefit in planning all
the way to the end rather than considering only a few steps ahead, and in fact, in some cases, it can be harmful to consider the distant future. Noting that
as the planning horizon decreases, planning complexity decreases, this helps recover tractability even in this regime.

Thus, while lower levels of noise render the planning problem tractable due to the decoupling result requiring no replanning, planning under medium noise remains tractable due to only occasional replanning, while for high levels of noise, tractability ensues because the planning horizon should shrink as the uncertainty increases.  When noise inundates the system, long-term predictions become so
uncertain that the best-laid plans will very likely run awry, and thus, it would be wasteful to
invest significant time thinking very far ahead. 
To examine this somewhat intuitive truth more quantitatively,
the parameter $\epsilon$ will be a knob we adjust, exploring these
aspects in the subsequent empirical analysis. \textit{We reiterate that the notion of low, medium and high noise regimes may seem somewhat vague, however, we provide precise definitions of these regimes using our empirical results later in this paper.}

%\subsubsection*{Small noise} The decoupling principle in \cite{decoupling} states that under small noise, as $\epsilon \downarrow 0$, the design of the feedback law can be conducted decoupled from the design of the open-loop optimized trajectory. The expected stochastic cost is equal to the nominal cost with very high probability as $\epsilon \rightarrow 0^{+}$. Since the deviation in the cost is minimal, the decoupled feedback law provides a near-optimal solution without needing for any replans.  
%\subsubsection*{Medium noise}
%At medium noise levels, the probability of the cost deviating from the nominal grows higher. So, replanning becomes necessary to avoid the cost from deviating significantly from the nominal. As noise increases, the probability of the cost deviation increases and replanning frequency increases. 
%\subsubsection*{High noise (hypothesis)} At high noise levels, constant replanning becomes indispensable where the cost is dominated by the noise. At such high levels ($\epsilon > 0.8$), our model in \eqref{eq:model} becomes significantly affected by noise so that planning ahead will have no significance when compared to planning only for a short horizon. Moreover, since planning for a short horizon produces greedy solutions, which is preferred here, as what lies in the future is highly uncertain.        

%
\section{The Planning Algorithms}
\label{section:algorithms}
The preliminaries and the algorithms are explained below:% 
\vspace{-10pt}
\subsection{Deterministic Optimal Control Problem:}
Given the initial state $\V{x}_0$ of the system, the solution to the deterministic OCP is given by \eqref{DOCP},
%
\begin{comment}
\begin{equation}
    J^{*}(\V{x}_0) = \min_{\V{u}_{0:T-1}} \left[\sum^{T-1}_{t=0} c_t(\V{x}_t, \V{u}_t) + c_T(\V{x}_T)\right], 
    \label{DOCP}
\end{equation}%
\end{comment}
\begin{comment}
\begin{align*}
subject to \ \V{x}_{t+1} = f(\V{x}_t) + \M{B}_t \V{u}_t,\\
   \V{u}_{\text{min}} \leq \V{u}_t \leq \V{u}_{\text{max}},\\
    | \V{u}_{t} - \V{u}_{t-1}| \leq \Delta \V{u}_{\text{max}}.
\end{align*}
\end{comment}
s.t. $\V{x}_{t+1} = f(\V{x}_t) + \M{B}_t \V{u}_t,\ \V{u}_{\text{min}} \leq \V{u}_t \leq \V{u}_{\text{max}},\ | \V{u}_{t} - \V{u}_{t-1}| \leq \Delta \V{u}_{\text{max}}.$ 
The last two constraint model physical limits that impose upper bounds and lower bounds on control inputs and rate of change of control inputs. 
The solution to the above problem gives the open-loop control inputs $\overline{\V{u}}_{0:T-1}$ for the system.
For our problem, we take a quadratic cost function for state and control as 
$
    c_t(\V{x}_t,\V{u}_t) = \tr{\V{x}}_t  \M{W}^x\V{x}_t + \tr{\V{u}}_t\M{W}^u \V{u}_t,
$
$
    c_T(\V{x}_T) = \tr{\V{x}}_T\M{W}^x_f\V{x}_T,
$
where $\M{W}^x,\ \M{W}^x_f  \succeq \M{0}$ and $\M{W}^u \ \succ \ \M{0}$.
\vspace{-10pt}
\subsection{Model Predictive Control (MPC):}
We employ the non-linear MPC algorithm due to the non-linearities associated
with the motion model. The MPC algorithm implemented here solves the
deterministic OCP~\eqref{DOCP} at every time step, applies the control inputs
computed for the first instant and uses the rest of the solution as an initial
guess for the subsequent computation. In the next step, the current state of
the system is measured and used as the initial state and the process is
repeated.
\vspace{-10pt}
\subsection{Short Horizon MPC (MPC-SH):}
We also implement a variant of MPC which is typically used in practical
applications where it solves the OCP only for a short horizon rather than the
entire horizon at every step. At the next step, a new optimization is solved
over the shifted horizon. This implementation gives a greedy solution but is
computationally easier to solve. It also has certain advantageous properties in
high noise cases which will be discussed in the results section. We denote the short planning horizon as $H_c$ also called as the control horizon, upto which the controls are computed. 
%A generic algorithm for MPC is shown in Algorithm~\ref{MPC_algo}. 
%As mentioned previously, MPC is the ``gold standard" for comparison in this owing to Fleming's paper showing the $O(\epsilon^4)$ near optimality of the MPC law compared to the true stochastic law \cite{fleming}.%
%
\begin{comment}
\begin{algorithm}[h]
\SetAlgoLined
    \KwIn{$\V{x}_0$ -- initial state, $\V{x}_g$ -- final state, $T$ -- time horizon, $H_c$ -- control horizon, $\Delta t$ -- time step, $\mathcal{P}$ -- system and environment parameters.}
    \For{$t \leftarrow 0$ \KwTo $T-1$}{
        $\V{u}_{t:t+H_c-1} \leftarrow$ OCP($\V{x}_{t},\V{x}_g, \min(H_c, T\!-\!t),\V{u}_{t-1}, \V{u}_{\textrm{guess}},\mathcal{P}$)\\
        $\phantom{xxxx}\V{x}_{t+1} \leftarrow$ $f(\V{x}_{t}) + \M{B}_t(\V{u}_t + \epsilon\V{w}_t)$ 
    }
    \caption{MPC algorithm\label{MPC_algo}.}
\end{algorithm}

\vspace{-20pt}
\end{comment}
\vspace{-10pt}
\subsection{Trajectory Optimised Linear Quadratic Regulator~\mbox{(T-LQR)}:}
\label{sec:lqr_gains}

As discussed in Section~\ref{section:decoupling}, stochastic optimal control problem can be decoupled and solved by designing an optimal open-loop (nominal) trajectory and a decentralized LQR policy to track the nominal. 
\\
\textit{Design of nominal trajectory}: The nominal trajectory is generated by first finding the optimal open-loop control sequence by solving the deterministic OCP~\eqref{DOCP} for the system. Then, using the computed control inputs and the noise-free dynamics, the sequence of states traversed $\overline{\V{x}}_{0:T}$ can be calculated.\\
\textit{Design of feedback policy:} In order to design the LQR controller, the system is first linearised about the nominal trajectory ($\overline{\V{x}}_{0:T}$, $\overline{\V{u}}_{0:T-1}$). Using the linear time-varying system, the feedback policy is determined by minimizing a quadratic cost as shown in~\eqref{LQRcost}. % 
% %
% \begin{equation}
%     \min_{\pi} \Exp{\V{w}_{t}}{\sum^{T-1}_{t=0} \tr{(\tilde{\V{x}}^\ell_t)} \M{W}^x (\tilde{\V{x}}^\ell_t) + \tr{(\tilde{\V{u}}^\ell_t)} \M{W}^u (\tilde{\V{u}}^\ell_t) + \tr{(\tilde{\V{x}}^\ell_T)} \M{W}^x_f (\tilde{\V{x}}^\ell_T)}.
% \end{equation}
% The linearised system is 
% \begin{equation}
%     \tilde{\V{x}}^\ell_{t+1} = \M{A}_t \tilde{\V{x}}^\ell_{t} +   \M{B}_t \tilde{\V{u}}^\ell_{t} + \epsilon\M{B}_t\V{w}_t ,
% \end{equation}
% %
% where $ \M{A}_t \ =\ \nabla_{\V{x}} f(\V{x},\V{u})|_{\overline{\V{x}}_t,\overline{\V{u}}_t}$, $ \M{B}_t \ =\ \nabla_{\V{u}} f(\V{x},\V{u})|_{\overline{\V{x}}_t,\overline{\V{u}}_t}$ and $\M{W}^x_f, \M{W}^x \succeq \V{0}, \M{W}^u \ \succ \V{0}$. Note that $\tilde{\V{x}}^\ell_{t}$ and $\tilde{\V{u}}^\ell_{t}$ denote the deviations from the nominal (the $\delta$s are elided for clarity).\\ 
The linear quadratic stochastic control problem~\eqref{LQRcost} can be easily solved using the Riccati equation and the resulting policy is $\delta{\V{u}}_{t} = -\M{L}_t\delta{\V{x}}^\ell_t$. The feedback gain and the Riccati equations are given by
$
\M{L}_t = \inv{(\M{R} + \M{B}^T_t\M{P}_{t+1} \M{B}_t)} \tr{\M{B}}_t \M{P}_{t+1}  \M{A}_t
$
and
$
\M{P}_{t} = \tr{\M{A}}_t \M{P}_{t+1} \M{A}_t - \tr{\M{A}}_t \M{P}_{t+1} \M{B}_t  \M{L}_t +  \M{Q},
$
respectively where $\M{Q}_f, \M{Q} \succeq \V{0}, \M{R} \ \succ \V{0}$ are the weight matrices for states and control and the terminal condition is
$\M{P}_{T} = \M{Q}_f $.   
\begin{comment}
\begin{equation}
    \M{L}_t = \inv{(\M{R} + \M{B}^T_t\M{P}_{t+1} \M{B}_t)} \tr{\M{B}}_t \M{P}_{t+1}  \M{A}_t,
    \label{LQR_gain}
\end{equation}

\begin{equation}
    \M{P}_{t} = \tr{\M{A}}_t \M{P}_{t+1} \M{A}_t - \tr{\M{A}}_t \M{P}_{t+1} \M{B}_t  \M{L}_t +  \M{Q}, 
    \label{Riccati}
\end{equation}
\end{comment}
\vspace{-10pt}
\subsection{T-LQR with Replanning (\mbox{T-LQR2}):}
T-LQR performs well at low noise levels, but at medium and high noise levels the system tends to deviate from the nominal. So, we define a threshold $J_{\textrm{thresh}}= \frac{J_{0:t} - \overline{J}_{0:t}}{\overline{J}_{0:t}}$, where $J_{0:t}$ denotes the actual cost during execution till time $t$. The factor $J_{\textrm{thresh}}$ measures the percentage deviation of the online trajectory from the nominal, and replanning is triggered for the system from the current state for the remainder of the horizon if the deviation exceeds it. Other replanning criteria such as state deviation can also be considered but we stick to the cost deviation in the following \footnote{In the absence of a running cost, a criterion such as state deviation could be used. Since we aim to optimize the cost, a criterion based on cost seems more reasonable.}. Note that if we set $J_{\textrm{thresh}}=0$, T-LQR2 reduces to MPC. The calculation of the new nominal trajectory and LQR gains are carried out similarly to the explanation in Section~\ref{sec:lqr_gains}. A generic algorithm for T-LQR and T-LQR2 is shown in Algorithm~\ref{TLQR_algo}. The implementations of all the algorithms are available at \url{https://github.com/MohamedNaveed/Stochastic_Optimal_Control_algos}.

\begin{comment}
\subsection{Short Horizon T-LQR with Replanning (T-LQR2-SH):}
A T-LQR equivalent of MPC-SH is also implemented where the nominal is planned only for a short horizon and it is tracked with a feedback policy as described in T-LQR. It also inherits the replanning property of T-LQR2.\\
\end{comment}

%\begin{comment}
%\vspace{-12pt}
\begin{algorithm}[!htpb]
\SetAlgoLined
\SetKwProg{Fn}{Function}{ is}{end}
\SetKwComment{Comment}{/*}{*/}   
    \KwIn{initial state $\V{x}_0$, final state $\V{x}_g$, time horizon $T$, replan threshold $J_{\textrm{thresh}}$, time step $\Delta t$, system and environment parameters $\mathcal{P}$.}
    
    \Fn{Plan($\V{x}_{0},\V{x}_g, T, \V{{u}}_{\textrm{init}}$, $\V{u}_{\textrm{guess}}$,$\mathcal{P}$)}{
    $\overline{\V{{u}}}_{0:T-1}$ $\leftarrow$ OCP($\V{x}_{0},\V{x}_g, T,\V{{u}}_{\textrm{init}}$, $\V{u}_{\textrm{guess}}$,$\mathcal{P}$)
    
    % \For{$t \leftarrow 0$ \KwTo $T-1$}{
    %     $\overline{\V{x}}_{t+1} \leftarrow f(\overline{\V{x}}_{t}) + \M{B}_t\overline{\V{u}}_t$
    % }
    $\overline{\V{x}}_{t+1} \leftarrow f(\overline{\V{x}}_{t}) + \M{B}_t\overline{\V{u}}_t$;\quad  $t=0,1, \cdots,T-1.$
    
    $\M{L}_{0:T-1} \leftarrow $
    $Compute\_LQR\_Gain(\overline{\V{x}}_{0:T-1},\overline{\V{u}}_{0:T-1}$)
    
    return $\overline{\V{x}}_{0:T},\overline{\V{u}}_{0:T-1},\M{L}_{0:T-1}$ 
    }
    \Fn{Main()}{
    $\overline{\V{x}}_{0:T}$,$\overline{\V{u}}_{0:T-1}$,$\M{L}_{0:T-1}$ $\leftarrow$ $\textrm{Plan}(\V{x}_{0},\V{x}_g, T,\mathbf{0}, \V{u}_{\textrm{guess}},\mathcal{P})$
    
    \For{$t \leftarrow 0$ \KwTo $T-1$}{
        
        $\V{u}_{t} \leftarrow \textrm{Constrain}(\overline{\V{u}}_{t} - \M{L}_{t}(\V{x}_{t} - \overline{\V{x}}_{t}))$ \tcp*[f]{ Enforce limits}
        
        %$\V{u}_{t} \leftarrow \textrm{Constrain}(\V{u}_{t})$ 
        
        $\V{x}_{t+1} \leftarrow f(\V{x}_{t}) + \M{B}_t(\V{u}_t + \epsilon\V{w}_t)$ 
        
        \If(\tcp*[f]{Replan?}){$ (J_{0:t} - \overline{J}_{0:t})/\overline{J}_{0:t}  > J_{\textrm{thresh}}$}{
            $\overline{\V{x}}_{t+1:T}, \overline{\V{u}}_{t+1:T-1}, \M{L}_{t+1:T-1} \leftarrow \textrm{Plan}(\V{x}_{t+1},\V{x}_g, T\!-\!t\!-\!1,\V{u}_t, \V{u}_{\textrm{guess}},\mathcal{P})$
            
        }
    }
    }
    \caption{T-LQR2 algorithm}\label{TLQR_algo}

\end{algorithm}

%\end{comment}
\vspace{-10pt}
\subsection{Multi-Agent versions}
The MPC version of the multi-agent planning problem is reasonably straightforward except that the complexity of the planning increased (exponentially) in the number of agents. Also, we note that the agents have to always communicate with each other in order to do the planning. The Multi-agent Trajectory-optimised LQR (MT-LQR) version is also relatively straightforward in that the agents plan the nominal path jointly once, and then the agents each track their individual paths using their decoupled feedback controllers. There is no communication whatsoever between the agents during this operation.

The MT-LQR2 version is a little more subtle. The agents have to periodically
replan when the total cost deviates more than $J_{\textrm{thresh}}$ away from the
nominal, i.e., the agents do not communicate until the need to replan arises.
In general, the system would need to detect this in a distributed fashion, and
trigger replanning. 
\begin{comment}
\textit{We postpone consideration of this important aspect of the problem to a
subsequent paper more directly focused on networking considerations and instead, we will
assume that there exists a (yet to be determined) distributed strategy that
would perform the detection and replanning.}
\end{comment}
%
\section{Simulation Results:}
\label{section:results}

We test the performance of the algorithms extensively in three different non-linear models namely, the car-like robot model, car with two trailers and a quadrotor. Due to space constraints, only the results for the car-like robot are shown below, however, that the trends are generalizable can be seen from the results on the other models that are shown in the supplementary material. Numerical optimization is carried out using \texttt{CasADi} framework \cite{Andersson2018} with \texttt{Ipopt} \cite{Ipopt} NLP solver in \texttt{Python}. To provide a good estimate of the performance, the results presented were averaged from 100 simulations for every value of noise considered. Simulations were carried out in parallel across 100 cores in a cluster equipped with Intel Xeon 2.5GHz E5-2670 v2 10-core processors.
\subsection*{Car-like robot model:}
The car-like robot considered in our work has the motion model described by 
$x_{t+1} = x_t + v_{t}\cos(\theta_t)\Delta t,\ y_{t+1} = y_t + v_{t}\sin(\theta_t)\Delta t,\ \theta_{t+1} = \theta_{t} + \frac{v_t}{L}\tan(\phi_t)\Delta t,\ \phi_{t+1} = \phi_{t} + \omega_t \Delta t,$
\begin{comment}
\begin{align*}
    x_{t+1} &= x_t + v_{t}\cos(\theta_t)\Delta t,  &  \theta_{t+1} &= \theta_{t} + \frac{v_t}{L}\tan(\phi_t)\Delta t, \\
    y_{t+1} &= y_t + v_{t}\sin(\theta_t)\Delta t,  &  \phi_{t+1} &= \phi_{t} + \omega_t \Delta t,
\end{align*}
\end{comment}
%
where $\tr{(x_t, y_t, \theta_t, \phi_t)}$ denote the robot's state vector
namely, robot's $x$ and $y$ position, orientation and steering angle at time $t$.
Also, $\tr{(v_t, \omega_t)}$ is the control vector and denotes the robot's linear
velocity and angular velocity (i.e., steering). Here $\Delta t$ is the
discretization of the time step.  
\vspace{-10pt}
\subsection*{Noise characterization:} 
We add zero mean independent identically distributed (i.i.d), random sequences ($\V{w}_t$) as actuator noise to test the performance of the control scheme. The standard deviation of the noise is $\epsilon$ times the maximum value of the corresponding control input, where $\epsilon$ is a scaling factor which is varied during testing, that is:
$
\V{w}_t = \V{u}_{\textrm{max}} \bm{\nu};  \quad \bm{\nu} \sim \mathcal{N}(\V{0}, \M{I}) 
$
and the noise is added as $\epsilon \V{w}_t$. Note that, we enforce the constraints in the control inputs before the addition of noise, so the controls can even take a value higher after noise is added. The analyses can be done with process noise as well, but $\epsilon$ loses meaning in such a scenario and the plots would just shift depending on the variance of $\V{w}_t$. Since all the algorithms use the same noise model and having been tested in an extensive range of values, the requirement for a process noise model is not really necessary.  
\subsection{Single agent setting:}
\begin{figure}[h]
    \centering
    \subfloat[Full noise spectrum]{
        \includegraphics[width=.30\textwidth]{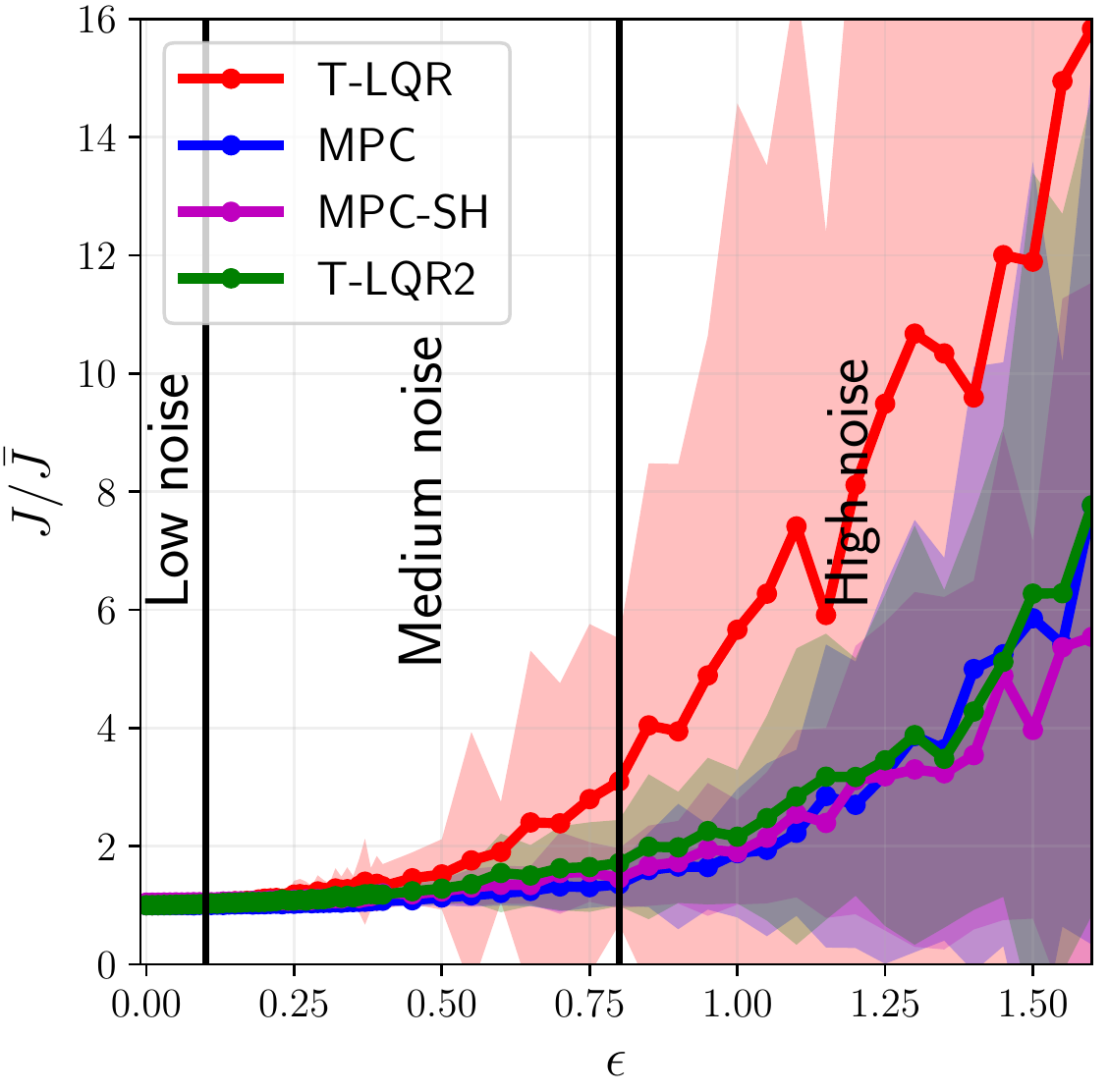}
        \label{1 agent cost full}}
    \subfloat[Enhanced detail $0\leq\epsilon\leq 0.4$]{
            \includegraphics[width=.30\textwidth]{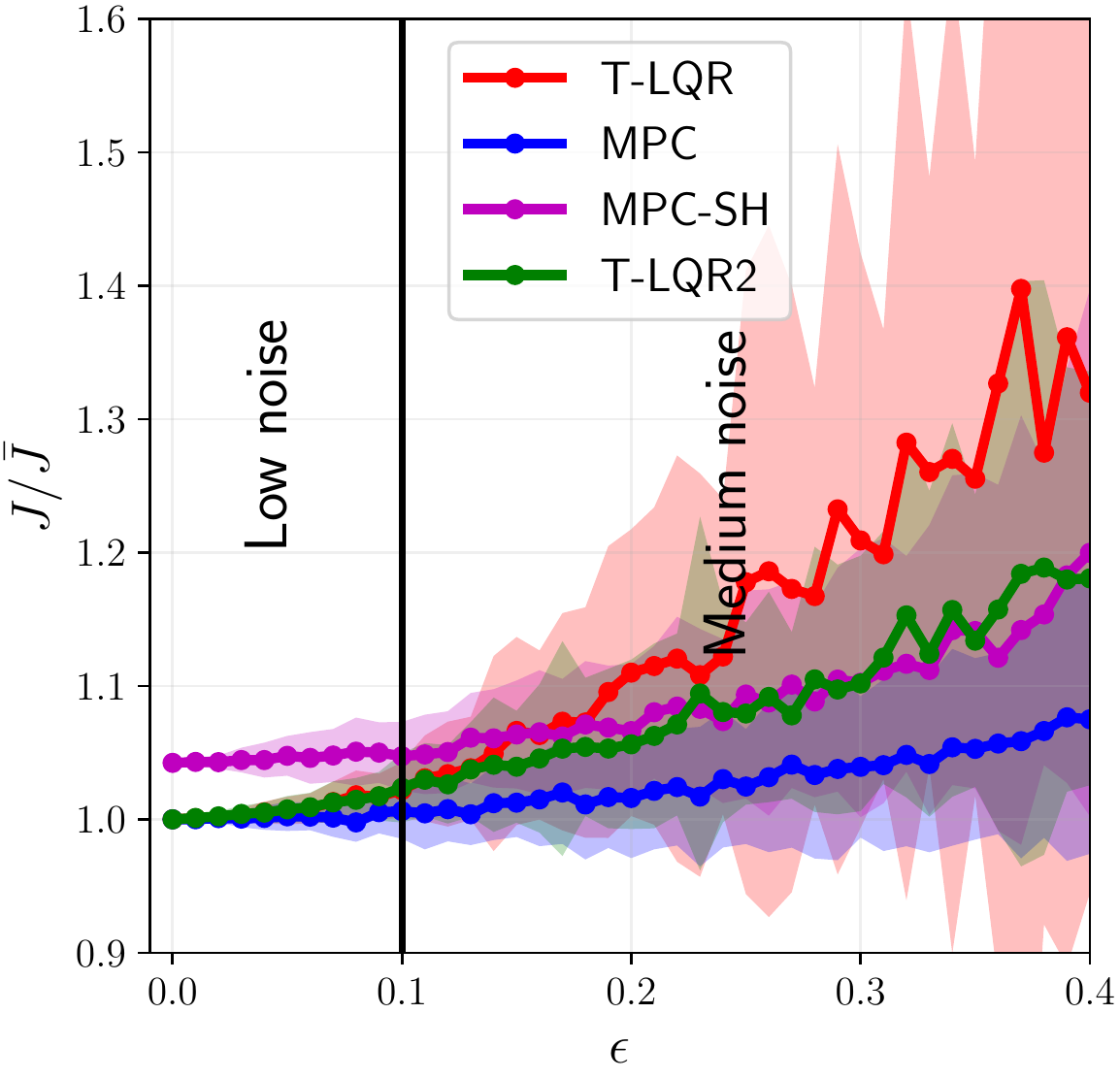}
        \label{1 agent cost low}}
    \subfloat[Replanning operations]{
        \includegraphics[width=.30\textwidth]{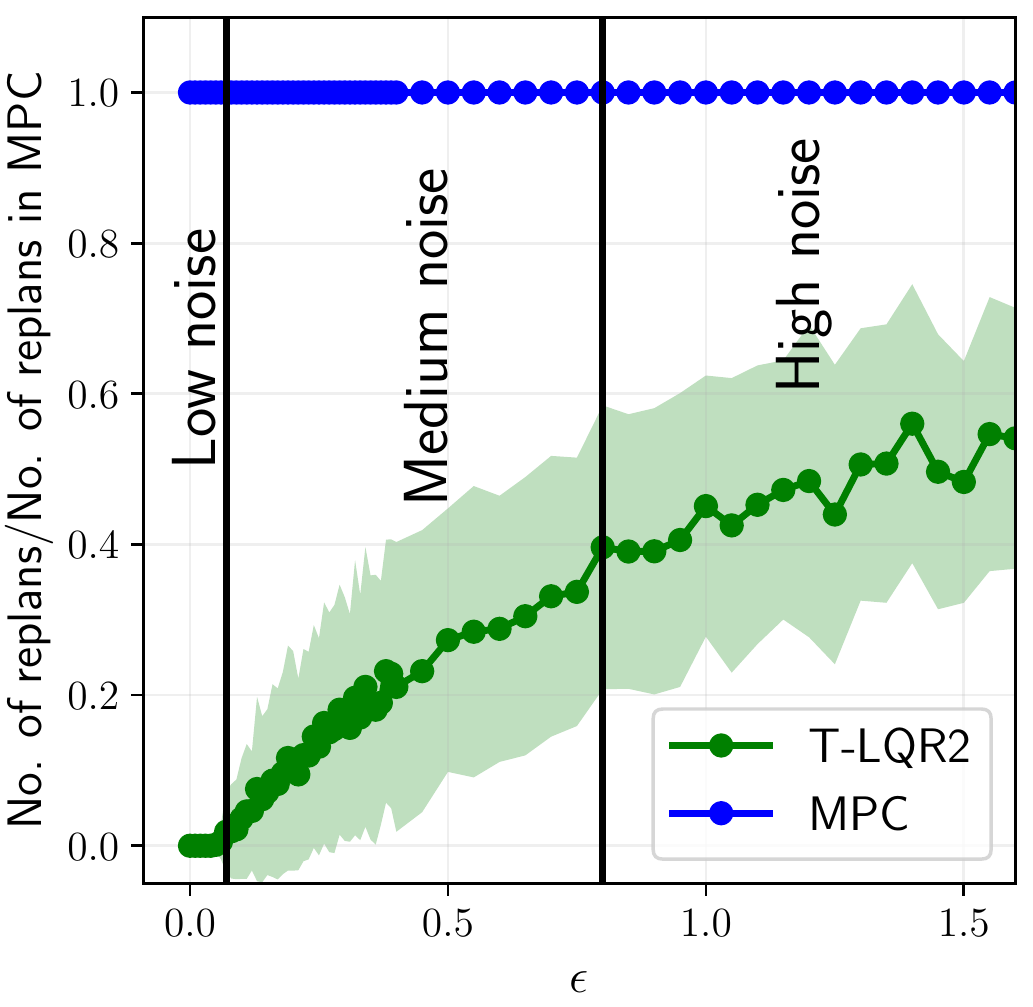}
        \label{1_agent_replan}}
        
    \caption{Cost evolution of the different algorithms for varying noise for a single agent. Control Horizon ($H_c$) used for MPC-SH was 7. $J_{\textrm{thresh}}$ = 2\% was the replanning threshold used. $J/\overline{J}$ is the ratio of the cost incurred during execution to the nominal cost and is used as the performance measure throughout the paper. The nominal cost $\overline{J}$ which is calculated by solving the deterministic OCP for the total time horizon, just acts as a normalizing factor here. The solid line in the plots indicates the mean and the shade indicates the standard deviation of the corresponding metric.}
    %\vspace*{-16pt}
    \label{1_agent_cost}
\end{figure}

\begin{figure}[!htbp]
    \subfloat[MPC]{
        \includegraphics[width=.225\textwidth]{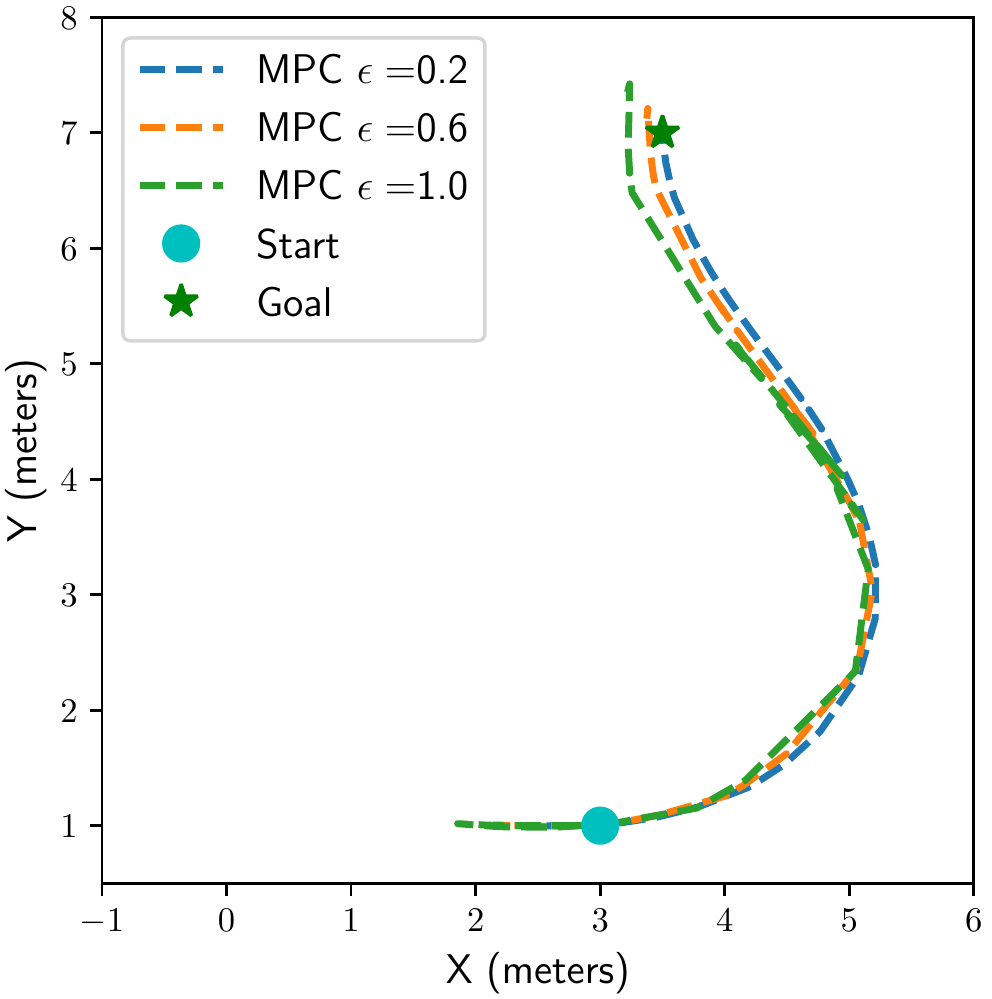}
        \label{test_cases_high_mpc}}
    \subfloat[T-LQR2]{
        \includegraphics[width=.225\textwidth]{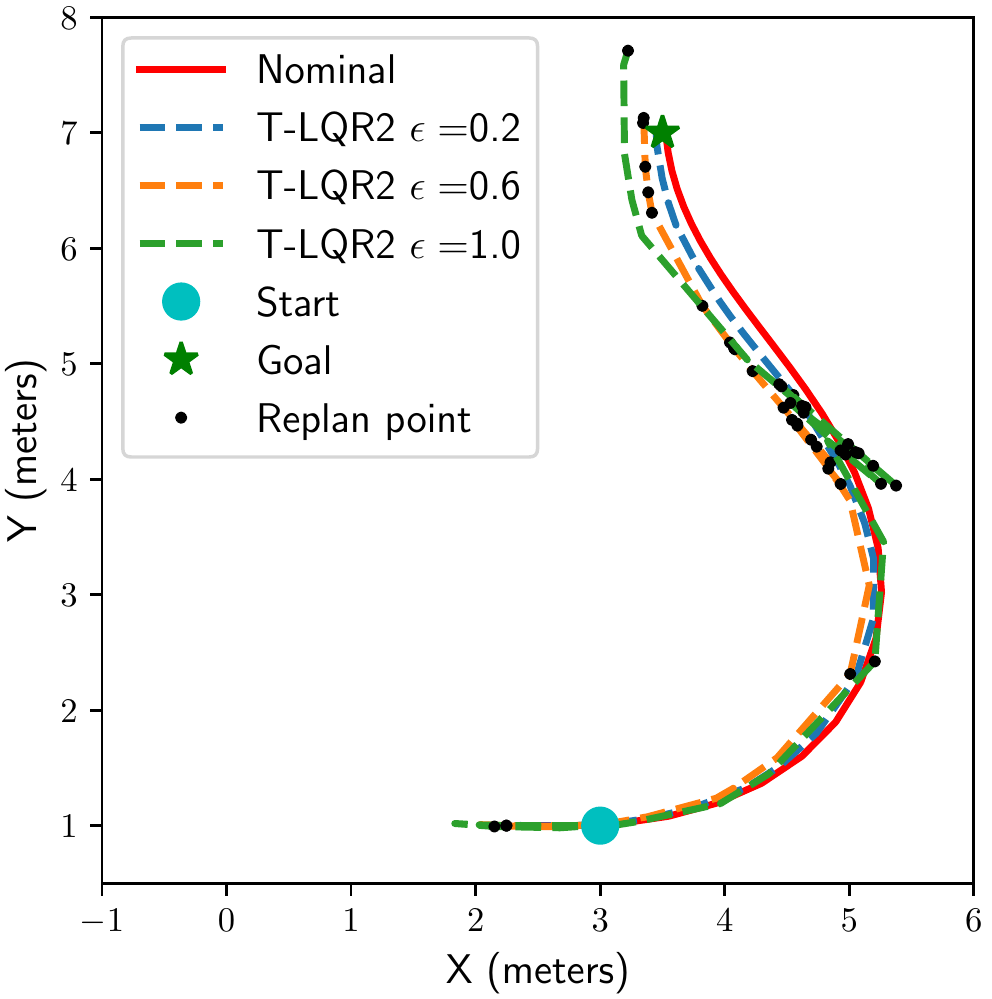}
        \label{test_cases_high_tlqr_replan}}
    \subfloat[MPC-SH]{        \includegraphics[width=.225\textwidth]{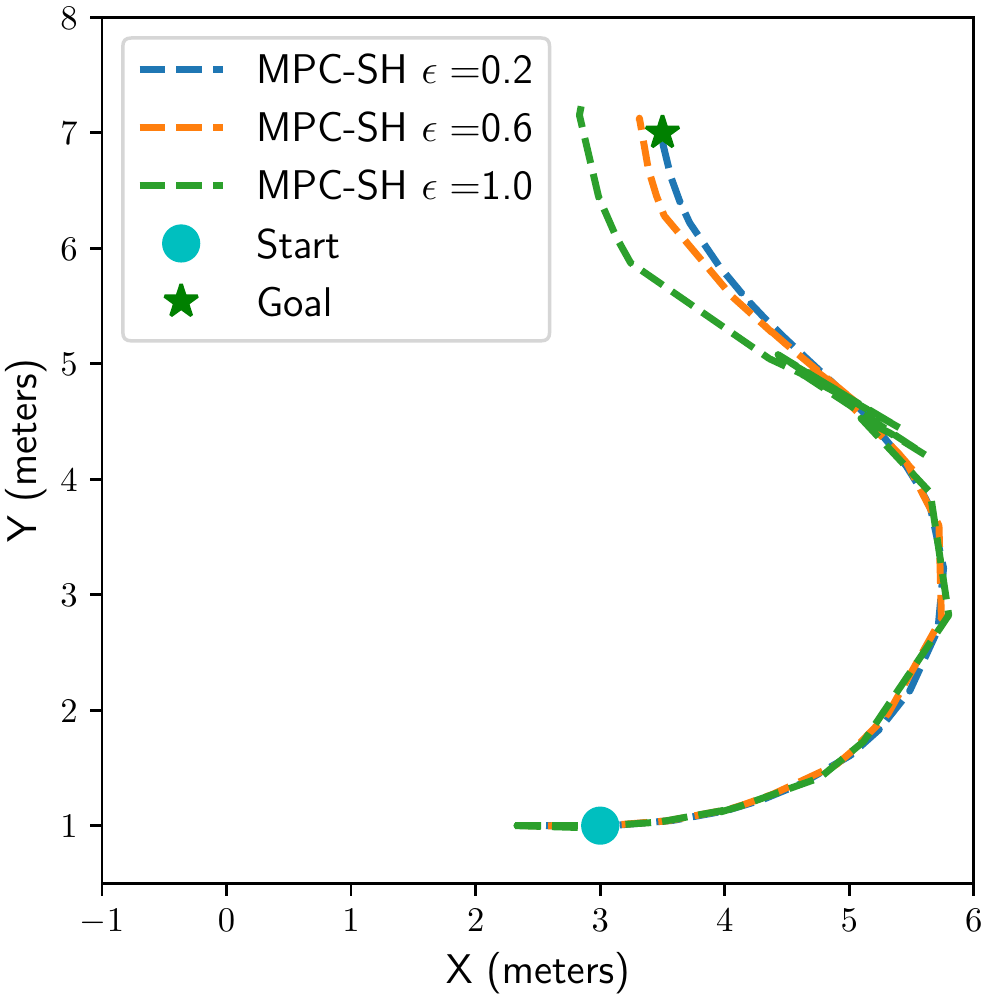}
        \label{test_cases_high_shmpc}}
    \subfloat[T-LQR]{
        \includegraphics[width=.225\textwidth]{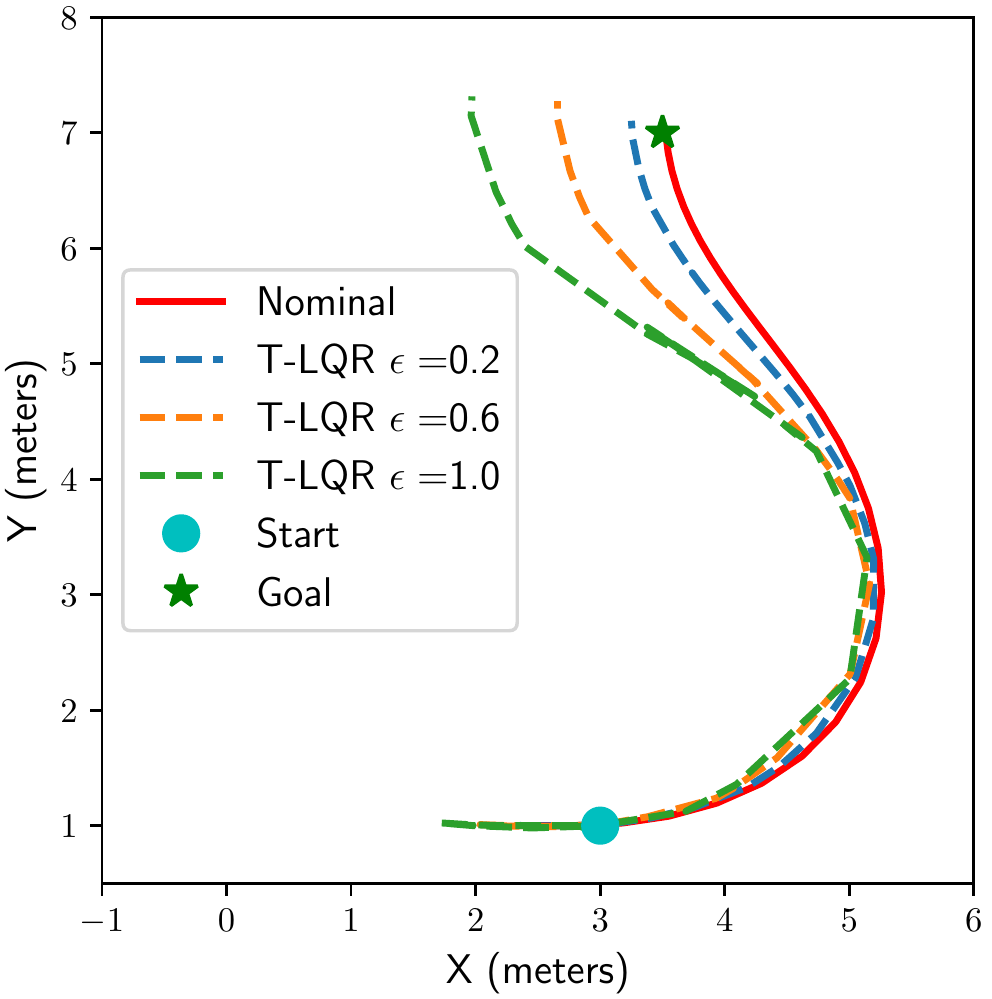}
        \label{test_cases_high_tlqr}}
    \caption{The figure shows the paths taken by the robot using a particular algorithm and how they change as $\epsilon$ varies. The key takeaway is that the paths taken under T-LQR are close to MPC for small values of $\epsilon$ and deviate as $\epsilon$ increases, while those of T-LQR2 are very close to the paths taken under MPC. It validates our claim about near-optimal performance of the decoupled feedback law under small noises and how performance can be preserved by replanning whenever necessary under medium noises. A equivalent plot for a scenario with obstacles is shown in the supplementary material.}
    \label{fig:test_cases_high}
\end{figure}

A car-like robot is considered and is tasked to move from a given initial pose
to a goal pose. The environment of the robot is shown in
Figure~\ref{fig:test_cases_high}. The experiment is done for all the control
schemes discussed and their performance for different levels of noise are shown
in Figure~\ref{1_agent_cost}. 
\subsection{Multi-agent setting:}
\begin{figure}[h]%\label{e-plot}
    \centering
    \subfloat[Full noise spectrum]{
        \includegraphics[width=.30\textwidth]{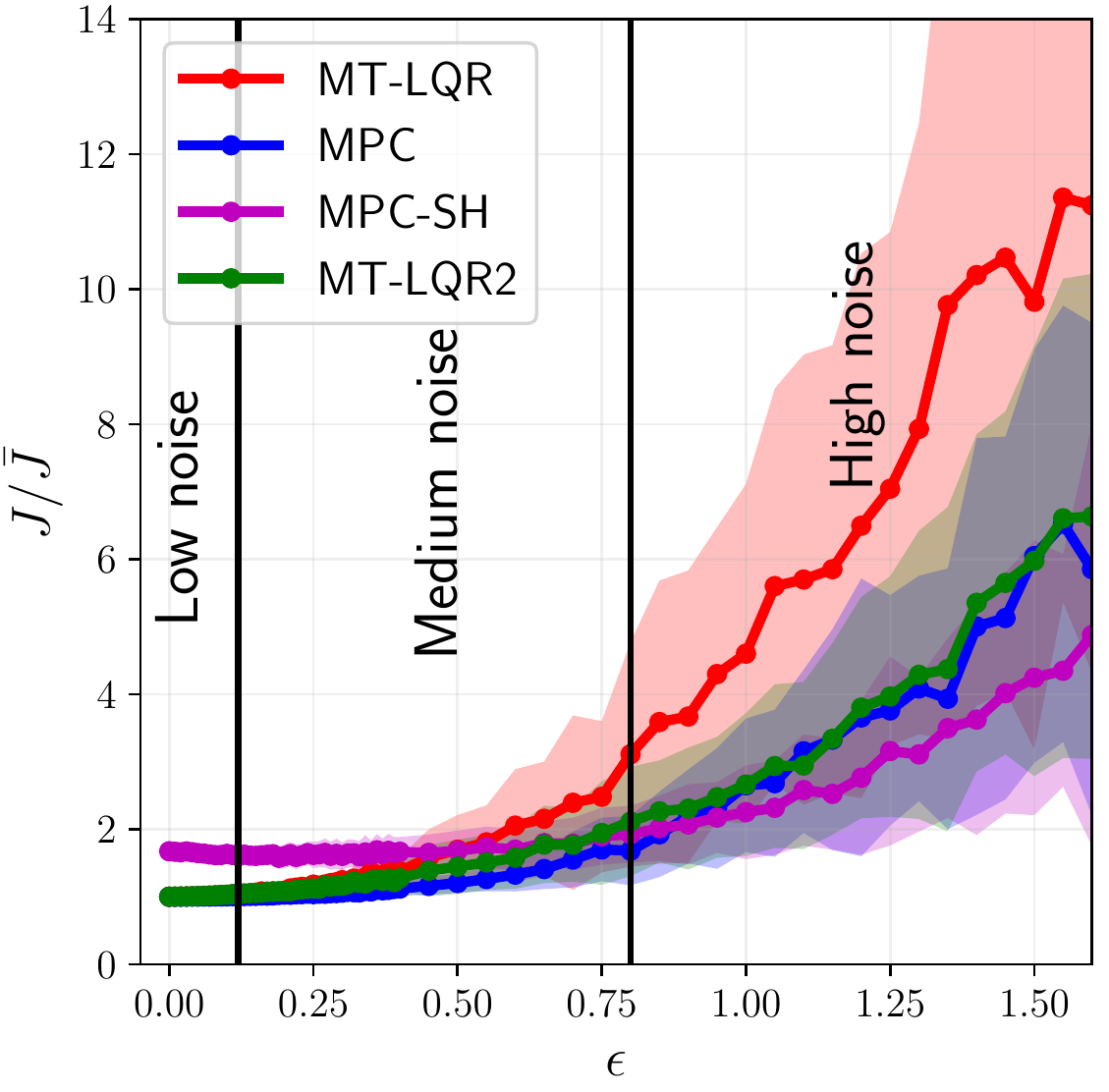}
        \label{3 agent cost full}}
    \subfloat[Enhanced detail $0\leq\epsilon\leq0.4$]{
        \includegraphics[width=.30\textwidth]{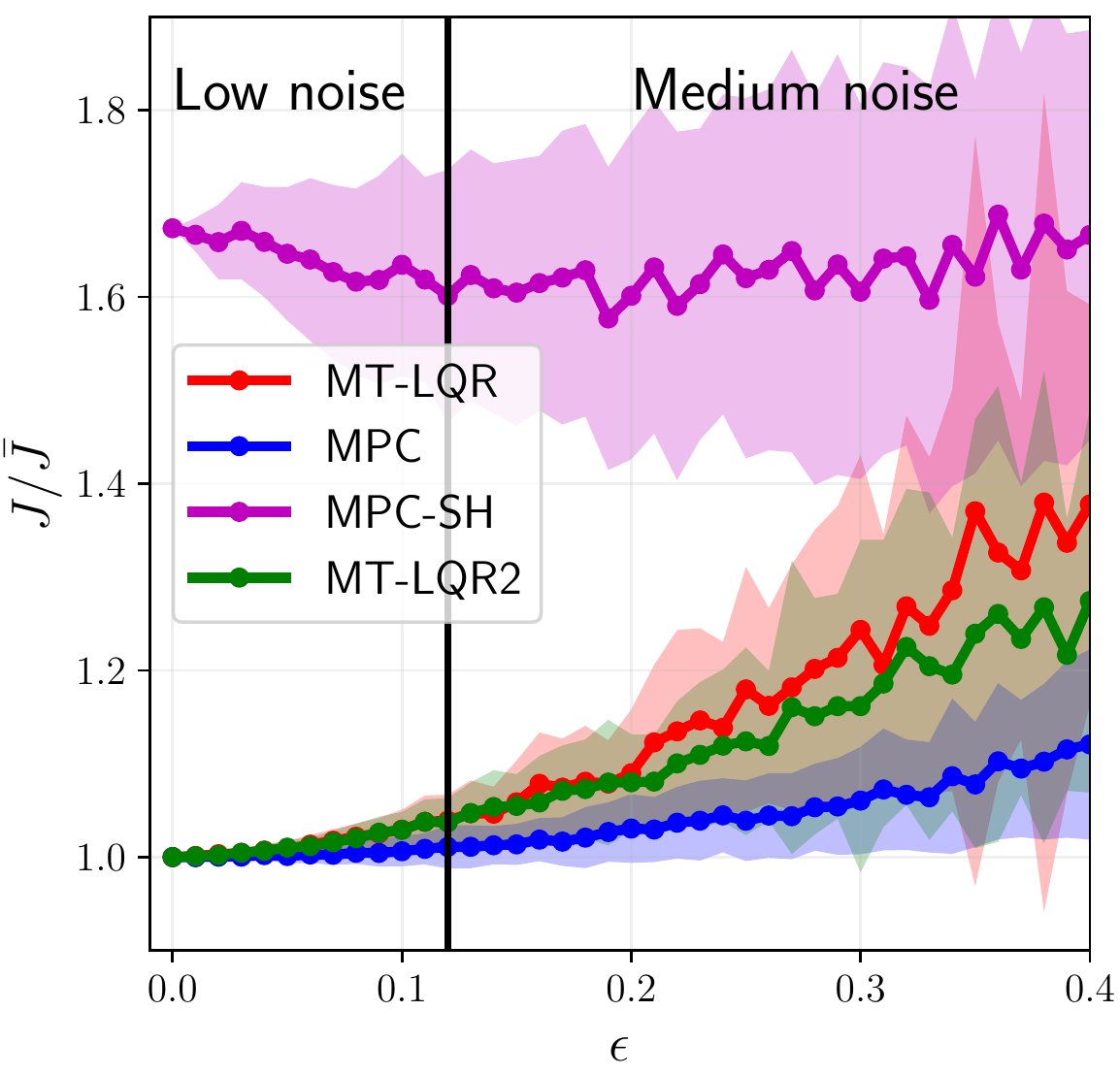}
        \label{3 agent cost low}}
    \subfloat[Replanning operations]{
        \includegraphics[width=.30\textwidth]{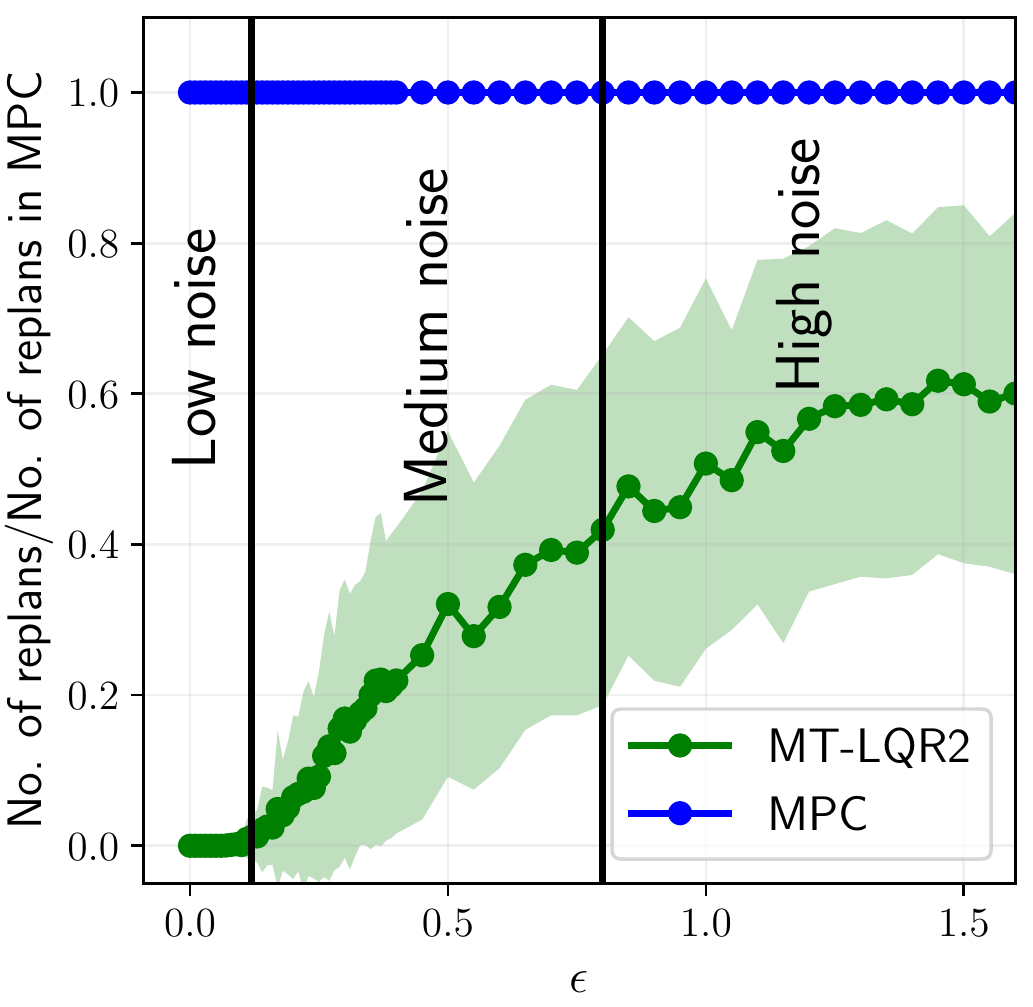}
        \label{3_agent_replan}}
     \caption{Cost evolution of the different algorithms for varying noise for 3 agents. Control Horizon ($H_c$) used for MPC-SH was 7. $J_{\textrm{thresh}}$ = 2\% was the replanning threshold used.}
    \label{3_agent_cost}
\end{figure}

A labelled point-to-point transition problem with 3 car-like robots is considered where each agent is assigned a fixed destination which cannot be exchanged with another agent. The performance of the algorithms is shown in Figure~\ref{3_agent_cost}. The cost function involves the state and control costs for the entire system similar to the single agent case. One major addition to the cost function is the penalty function to avoid inter-agent collisions which is given by 
$
    \Psi^{(i,j)} = \textrm{M}\exp\left(-(\Vert \V{p}_t^i - \V{p}_t^j\Vert_2^2 - r_{\textrm{thresh}}^2)\right) 
$
where $\textrm{M} > 0$ is a scaling factor, $\V{p}^i_t = (x^i_t, y^j_t)$ and $r_{\textrm{thresh}}$ is the desired minimum distance the agents should keep between themselves.    
%
%\vspace{-15pt}
 %test cases high noise
\begin{figure}[!htbp]
\floatbox[{\capbeside\thisfloatsetup{capbesideposition={right,top}}}]{figure}[]
 {
    \begin{subfloatrow}
        \ffigbox[\FBwidth][]
        {\caption{Solution quality for 1 agent with $\epsilon = 0.1$.}}
        {\includegraphics[width=.225\textwidth]{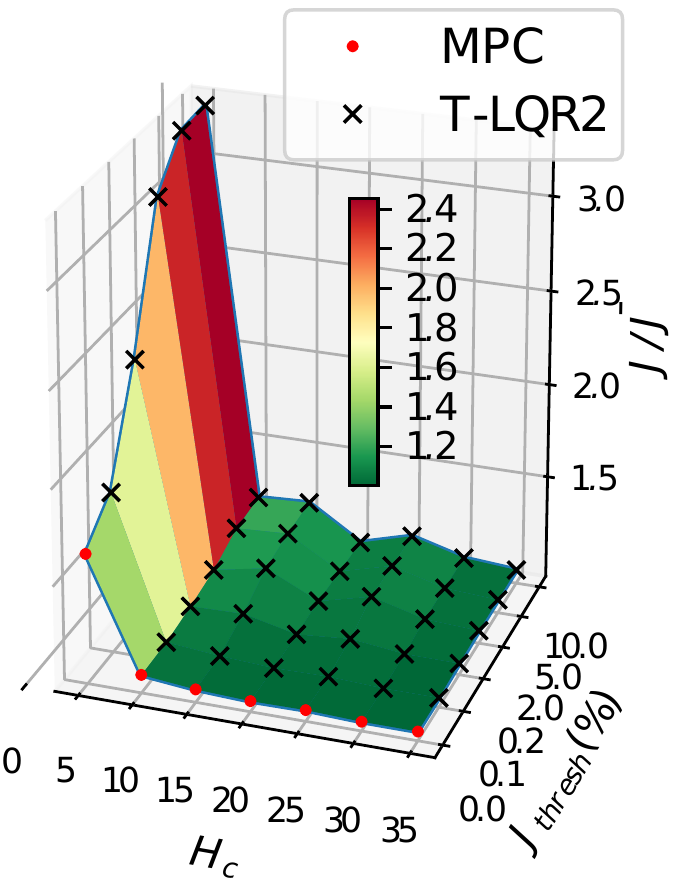}
        \label{fig:costvsHc_1_1}}
        \ffigbox[\FBwidth][]
        {\caption{Compute time for 1 agent with $\epsilon = 0.1$.}}
        {\includegraphics[width=.225\textwidth]{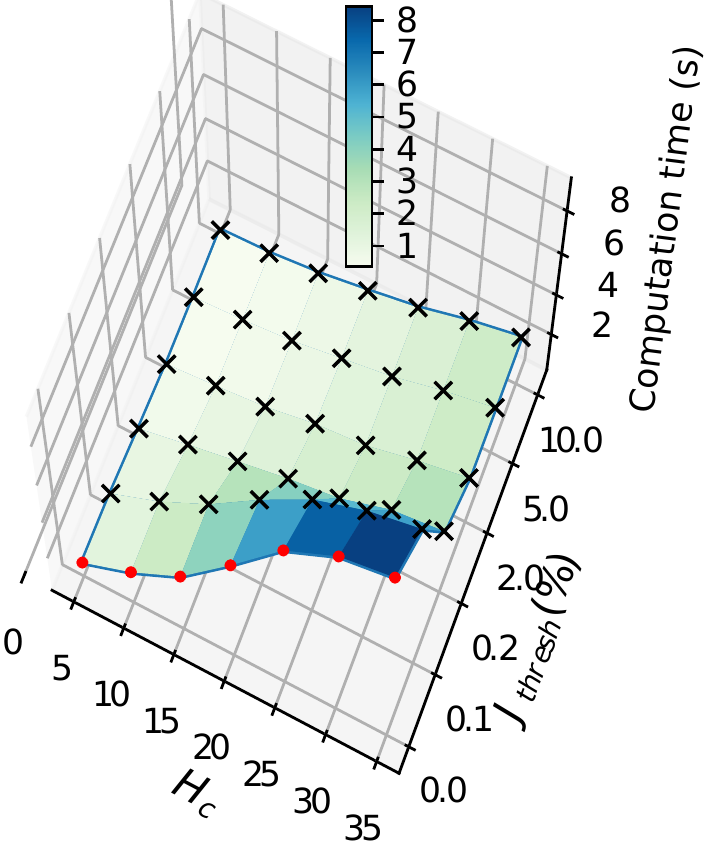}
        \label{fig:timevsHc_1_1}}
    \end{subfloatrow}
    \begin{subfloatrow}
        \ffigbox[\FBwidth][]
        {\caption{Solution quality for 1 agent with $\epsilon = 0.7$.}}
        {\includegraphics[width=.225\textwidth]{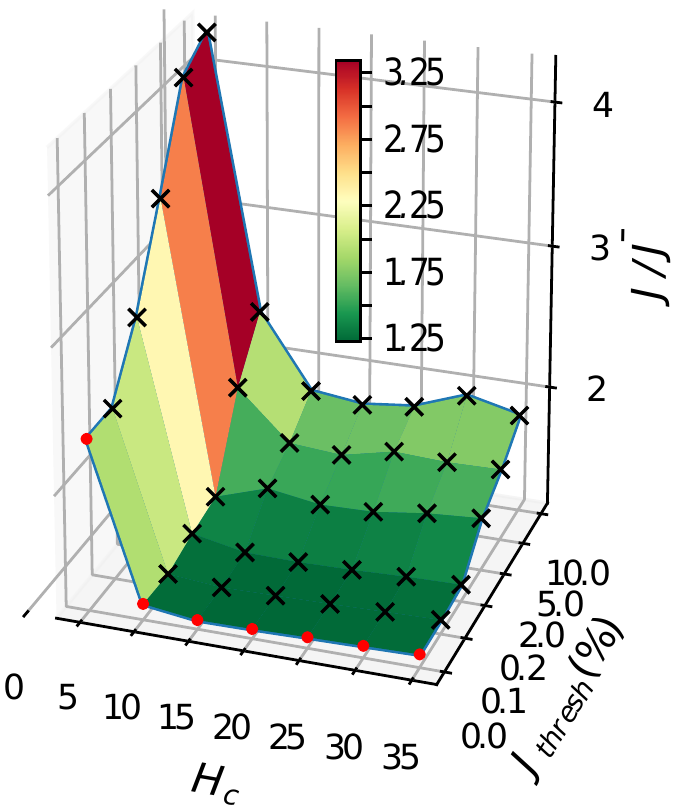}
        \label{fig:costvsHc_1_7}}
        \ffigbox[\FBwidth][]
        {\caption{Compute time for 1 agent with $\epsilon = 0.7$.}}
        {\includegraphics[width=.225\textwidth]{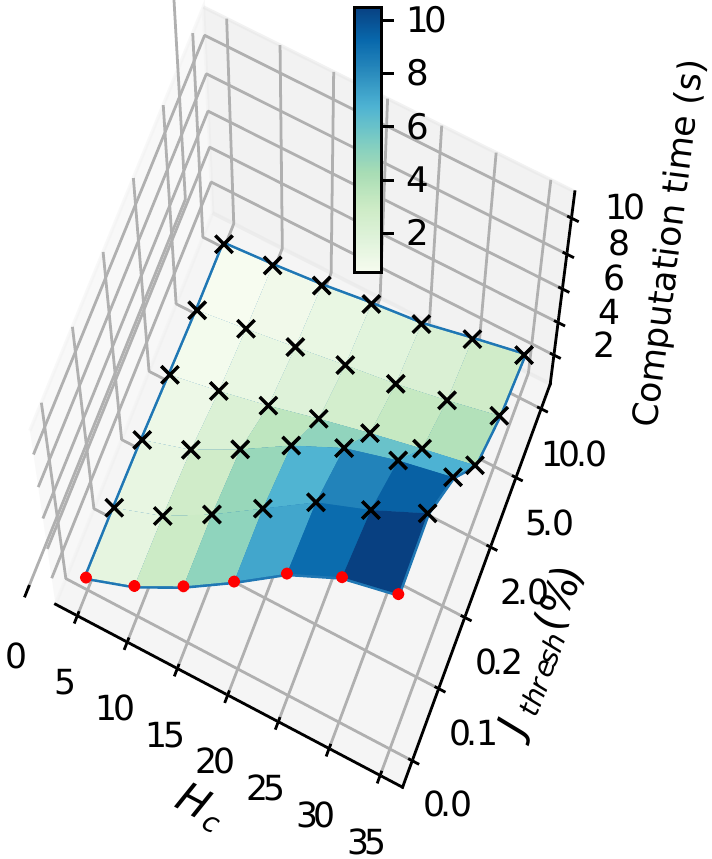}
        \label{fig:timevsHc_1_7}}
    \end{subfloatrow}
    \begin{subfloatrow}
        \ffigbox[\FBwidth][]
        {\caption{Solution quality for 3 agents with $\epsilon = 0.1$.}}
        {\includegraphics[width=.225\textwidth]{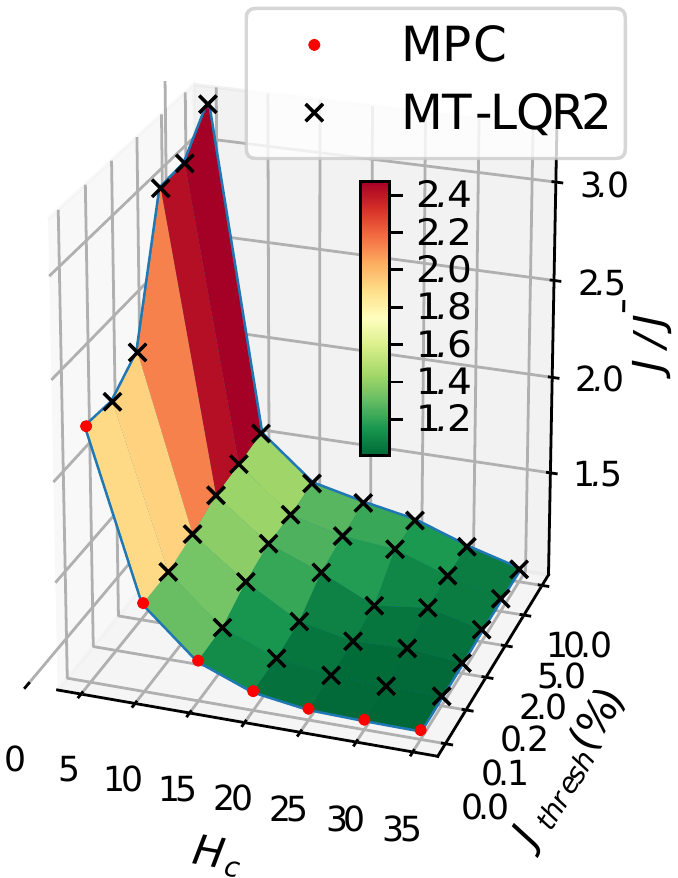}
        \label{fig:costvsHc_3_1}}
        \ffigbox[\FBwidth][]
        {\caption{Compute time for 3 agents with $\epsilon = 0.1$.}}
        {\includegraphics[width=.225\textwidth]{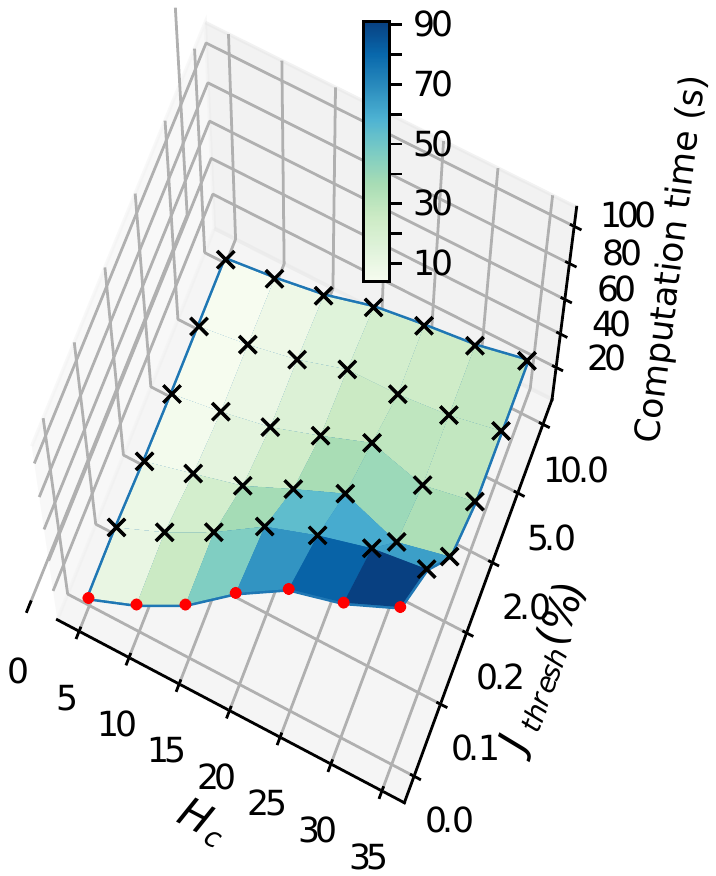}
        \label{fig:timevsHc_3_1}}
    \end{subfloatrow}
    \begin{subfloatrow}
        \ffigbox[\FBwidth][]
        {\caption{Solution quality for 3 agents with $\epsilon = 0.7$.}}
        {\includegraphics[width=.225\textwidth]{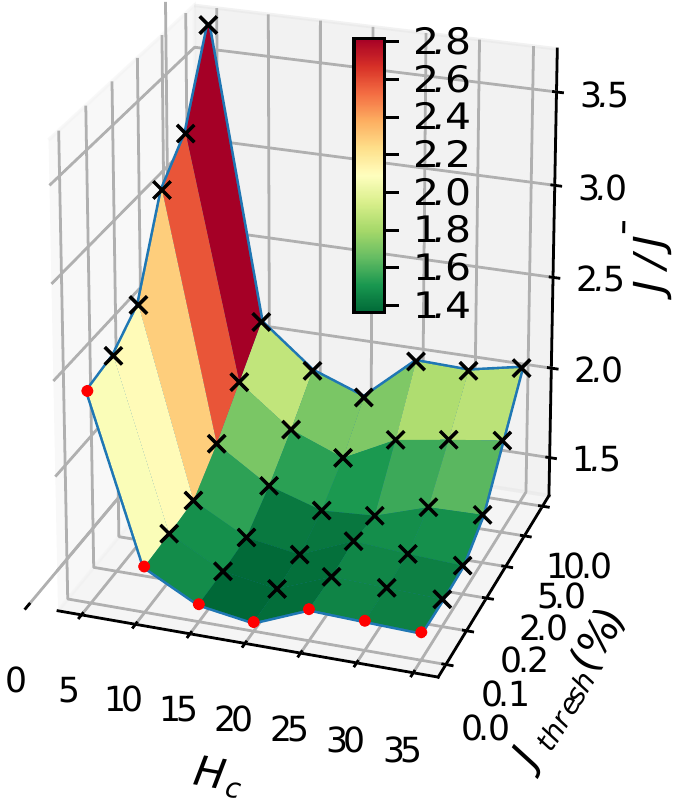}
        \label{fig:costvsHc_3_7}}
        \ffigbox[\FBwidth][]
        {\caption{Compute time for 3 agents with $\epsilon = 0.7$.}}
        {\includegraphics[width=.225\textwidth]{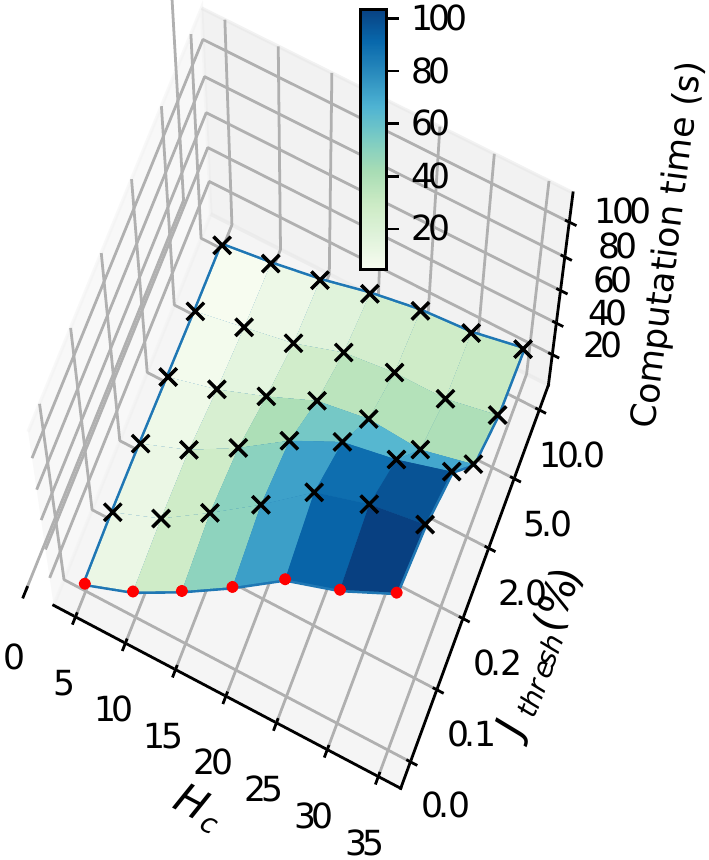}
        \label{fig:timevsHc_3_7}}
    \end{subfloatrow}}
    {\caption{We study the variation seen in cost incurred and computation time by changing the $J_{\textrm{thresh}}$ and control horizon ($H_c$) in T-LQR2/MT-LQR2 and MPC for a single agent and a case with 3 agents. Figures \ref{fig:costvsHc_1_1},~\ref{fig:timevsHc_1_1},~\ref{fig:costvsHc_1_7},~\ref{fig:timevsHc_1_7} show for a single agent and \ref{fig:costvsHc_3_1},~\ref{fig:timevsHc_3_1},~\ref{fig:costvsHc_3_7},~\ref{fig:timevsHc_3_7} for 3 agents. \ref{fig:costvsHc_1_1} and \ref{fig:timevsHc_1_1} show the performance in terms of cost and computation time, respectively, for the same experiment at $\epsilon = 0.1$ (low noise). Similarly, (c) and (d) show for $\epsilon=0.7$ (medium noise). To put into context, T-LQR2/MT-LQR2 replans less as $J_{\textrm{thresh}}$ is increased which in turn leads to decrease in computation time. The computation time also decreases with decrease in $H_c$. Though MPC does not have a threshold for replanning, it is plotted at $J_{\textrm{thresh}} = 0\%$ since it replans at every time step. We desire the performance indices, cost ($J/\bar{J})$ and computation time to be small. 
    \\
    As seen from \ref{fig:costvsHc_1_1}, for a fixed $H_c$, the performance of T-LQR2 is same as MPC and does not degrade as $J_{\textrm{thresh}}$ is increased (except when $H_c$ is too small). \ref{fig:timevsHc_1_1} shows the computational savings, where T-LQR2 is much better than MPC. From \ref{fig:costvsHc_1_1} and \ref{fig:timevsHc_1_1} it can be inferred that the decoupled approach provides good performance with substantial savings in computation time. Similarly, \ref{fig:costvsHc_1_7} shows that T-LQR2 performance is near-optimal to MPC for small $J_{\textrm{thresh}}$ and degrades as it is relaxed, which indicates the necessity of replanning to maintain optimality in the medium noise regime. The corresponding variation in computation time is shown in \ref{fig:timevsHc_1_7}, where T-LQR2 does slightly better. The computational savings can be increased further but by trading away optimality. \\
    Similar interpretation can be made for the multi-agent case as seen from Figures \ref{fig:costvsHc_3_1},~\ref{fig:timevsHc_3_1},~\ref{fig:costvsHc_3_7},~\ref{fig:timevsHc_3_7}. It can be seen that the positives of using the decoupled approach are amplified in the multi-agent case. \\
    Now, we analyse how $H_c$ affects the performance. Decreasing $H_c$ leads to greedy sub-optimal solutions. Though not clearly seen in single agent case, the decrease in performance is seen in the multi-agent case (Fig.~\ref{fig:costvsHc_3_1}). %{It can also be noted from \ref{fig:costvsHc_3_1} and \ref{fig:timevsHc_3_1} that for a fixed compute time, MT-LQR2 produces better solutions than MPC with a shorter horizon. This is because MT-LQR2 can achieve smaller compute times even while planning for the entire horizon while MPC has to shorten its horizon to meet the computation constraint.}% 
    But in high noises, decreasing the horizon does not lead to decrease in solution quality or sometimes even produces better solutions as seen in \ref{fig:costvsHc_3_7}.} 
    \label{fig:cost_3d}}
\end{figure}
\subsection{Definition of noise regimes and discussion of the results:}
Here we go on to define what exactly we mean by low, medium and high noise. The low noise regime as labelled in Figure \ref{1 agent cost low} and \ref{3 agent cost low} is the noise level at which the decoupled feedback law (T-LQR and MT-LQR) shows near-optimal performance compared to MPC and does not require any replanning operations. Beyond a limit replanning (T-LQR2 and MT-LQR2) is essential to constrain the cost from deviating away from the optimal and we call this as the medium noise regime. Figure~\ref{1_agent_replan} and~\ref{3_agent_replan} show the significant difference in the number of replanning operations, which determines the computational effort, taken by the decoupled approach compared to MPC. The significant difference in computational time between MPC and T-LQR2 can be seen from Figure~\ref{fig:timevsHc_1_1}. The trend is similar in the multi-agent case which again shows that the decoupling feedback policy is able to give computationally efficient solutions which are near-optimal in low noise cases by avoiding frequent replanning.

In the high noise regime, T-LQR2, MT-LQR2 and even MPC-SH perform on a par with MPC as seen from Figures~\ref{1 agent cost full} and~\ref{3 agent cost full} meaning, planning too far ahead is not beneficial at high noise levels. It can also be seen in Figure~\ref{fig:costvsHc_3_7} that the performance for MPC as well as MT-LQR2 is best at $H_c = 20$. Planning for a shorter horizon also eases the computation burden as seen in Figure~\ref{fig:timevsHc_3_7}. It can also be seen in Figure~\ref{3 agent cost full} where MPC-SH with $H_c=7$ outperforms MPC with $H_c=35$ at high noise levels which again show that the effective planning horizon decreases at the high noise regime.  

\section{Conclusions and Implications}
\label{section:conclusion}
In this paper, we have considered a class of stochastic motion planning problems for robotic systems over a wide range of uncertainty conditions parameterized in terms of a noise parameter $\epsilon$. We have shown extensive empirical evidence that a simple generalization of a recently developed ``decoupling principle" can lead to tractable planning without sacrificing performance for a wide range of noise levels. Future work will seek to treat the medium and high noise systems, considered here, analytically, and look to establish the near-optimality of the replanning scheme. Further, we shall consider the question of ``when and how to replan'' in a distributed fashion in the multi-agent setting, as well as relax the requirement of perfect state observation. It is also conjectured that by designing the linear feedback in a suitable fashion, the decoupling result can be made $O(\epsilon^4)$ near-optimal, thus making the algorithm theoretically as good as MPC owing to Fleming's result \cite{fleming1971stochastic}. Further, an important limitation of the method is the smoothness of the nominal trajectory such that suitable Taylor expansions are possible, this breaks down when trajectories are non-smooth such as in hybrid systems like legged robots, or maneuvers have kinks for car-like robots such as in a tight parking application. It remains to be seen as to if, and how, one may extend the decoupling to such applications. 

\bibliographystyle{abbrv}
\bibliography{bib_files/MAP_refs1,bib_files/MohammadRafi,bib_files/naveed_references}

\begin{thebibliography}{10}

\bibitem{Andersson2018}
J.~{A E Andersson}, J.~Gillis, G.~Horn, J.~B. Rawlings, and M.~Diehl.
\newblock {CasADi} -- {A} software framework for nonlinear optimization and
  optimal control.
\newblock {\em Mathematical Programming Computation}, In Press, 2018.

\bibitem{RLHD1}
R.~Akrour, A.~Abdolmaleki, H.~Abdulsamad, and G.~Neumann.
\newblock Model free trajectory optimization for reinforcement learning.
\newblock In {\em Proc. of the ICML}, 2016.

\bibitem{amato2013decentralizedB}
C.~Amato, G.~Chowdhary, A.~Geramifard, N.~K. Ure, and M.~J. Kochenderfer.
\newblock Decentralized control of partially observable markov decision
  processes.
\newblock In {\em Proc. IEEE Int. CDC}, pages 2398--2405, 2013.

\bibitem{bertsekas1}
D.~P. Bertsekas.
\newblock {\em Dynamic Programming and Optimal Control, vols I and II}.
\newblock Athena Scientific, Cambridge, MA, 2012.

\bibitem{boutilier1996planning}
C.~Boutilier.
\newblock Planning, learning and coordination in multiagent decision processes.
\newblock In {\em Proc. of the 6th conference on Theoretical aspects of
  rationality and knowledge}, pages 195--210. Morgan Kaufmann Publishers Inc.,
  1996.

\bibitem{bryson}
A.~E. Bryson and Y.~C. Ho.
\newblock {\em Applied Optimal control}.
\newblock Allied Publishers, NY, 1967.

\bibitem{T-MPC1}
L.~Chisci, J.~A. Rossiter, and G.~Zappa.
\newblock Systems with persistent disturbances: predictive control with
  restricted contraints.
\newblock {\em Automatica}, 37:1019--1028, 2001.

\bibitem{fleming1971stochastic}
W.~H. Fleming.
\newblock Stochastic control for small noise intensities.
\newblock {\em SIAM Journal on Control}, 9(3):473--517, 1971.

\bibitem{ETMPC1}
W.~Heemels, K.~Johansson, and P.~Tabuada.
\newblock An introduction to event triggered and self triggered control.
\newblock In {\em Proc. IEEE Int. CDC}, 2012.

\bibitem{RLHD4}
S.~Levine and P.~Abbeel.
\newblock Learning neural network policies with guided search under unknown
  dynamics.
\newblock In {\em Advances in NIPS}, 2014.

\bibitem{RLHD5}
S.~Levine and K.~Vladlen.
\newblock Learning complex neural network policies with trajectory
  optimization.
\newblock In {\em Proc. of the ICML}, 2014.

\bibitem{ETMPC2}
H.~Li, Y.~She, W.~Yan, and K.~Johansson.
\newblock Periodic event-triggered distributed receding horizon control of
  dynamically decoupled linear systems.
\newblock In {\em Proc. IFAC World Congress}, 2014.

\bibitem{Mayne_1}
D.~Q. Mayne.
\newblock Model predictive control: Recent developments and future promise.
\newblock {\em Automatica}, 50:2967--2986, 2014.

\bibitem{T-MPC3}
D.~Q. Mayne, E.~C. Kerrigan, E.~J. van Wyk, and P.~Falugi.
\newblock Tube based robust nonlinear model predictive control.
\newblock {\em International journal of robust and nonlinear control},
  21:1341--1353, 2011.

\bibitem{oliehoek2012decentralized}
F.~A. Oliehoek.
\newblock Decentralized pomdps.
\newblock {\em Reinforcement Learning}, 2012.

\bibitem{oliehoek2016concise}
F.~A. Oliehoek and C.~Amato.
\newblock {\em A concise introduction to decentralized POMDPs}.
\newblock Springer, 2016.

\bibitem{T-PFC}
K.~S. {Parunandi} and S.~{Chakravorty}.
\newblock T-pfc: A trajectory-optimized perturbation feedback control approach.
\newblock {\em IEEE RA-L}, 4(4):3457--3464, Oct 2019.

\bibitem{pynadath2002communicative}
D.~V. Pynadath and M.~Tambe.
\newblock The communicative multiagent team decision problem: Analyzing
  teamwork theories and models.
\newblock {\em Journal of Artificial Intelligence Research}, 16:389--423, 2002.

\bibitem{Mayne_2}
J.~B. Rawlings and D.~Q. Mayne.
\newblock {\em Model Predictive Control: Theory and Design}.
\newblock Nob Hill, Madison, WI, 2015.

\bibitem{T-MPC2}
J.~A. Rossiter, B.~Kouvaritakis, and M.~J. Rice.
\newblock A numerically stable state space approach to stable predictive
  control strategies.
\newblock {\em Automatica}, 34:65--73, 1998.

\bibitem{seuken2008formal}
S.~Seuken and S.~Zilberstein.
\newblock Formal models and algorithms for decentralized decision making under
  uncertainty.
\newblock {\em Int. Conf. on AAMAS}, 17(2):190--250, 2008.

\bibitem{RLHD3}
E.~Theodorou, Y.~Tassa, and E.~Todorov.
\newblock Stochastic differential dynamic programming.
\newblock In {\em Proc. of the ACC}, 2010.

\bibitem{RLHD2}
E.~Todorov and Y.~Tassa.
\newblock Iterative local dynamic programming.
\newblock In {\em Proc. of the IEEE Int. Symposium on ADP and RL.}, 2009.

\bibitem{Ipopt}
A.~W\"{a}chter and L.~T. Biegler.
\newblock On the implementation of a primal-dual interior point filter line
  search algorithm for large-scale nonlinear programming.
\newblock {\em Mathematical Programming}, 2006.

\bibitem{D2C1.0}
R.~Wang, K.~S. Parunandi, D.~Yu, D.~M. Kalathil, and S.~Chakravorty.
\newblock Decoupled data based approach for learning to control nonlinear
  dynamical systems.
\newblock {\em CoRR}, abs/1904.08361, 2019.

\end{thebibliography}

\newpage
\pagenumbering{arabic}
\setcounter{page}{1}
\section*{SUPPLEMENTARY MATERIAL}
\label{section:supplementary}
In this document, we provide details of the decoupling result, a  rudimentary analysis of the high noise regime, and more empirical results on different robotic models from that presented in the main paper.
\section{A NEAR OPTIMAL DECOUPLING PRINCIPLE}\label{sec3}
We discuss in detail the decoupling principle described in Section~\ref{section:decoupling}. 

We make the following assumptions for the simplicity of illustration.  We assume that the dynamics given in \eqref{eq:model} can be written in the  form 
\begin{align} 
\label{eq.0.1}
%\label{eq:dynamics1}
x_{t+1} = f(x_t) + B_{t} u_t + \epsilon w_t, 
\end{align}
where $\epsilon < 1$ is a small parameter. We also assume that the instantaneous cost $c(\cdot, \cdot)$ has the following simple form,  
\begin{align}
\label{eq:cost1}
c(x,u) = l(x) + \frac{1}{2} u'Ru.
\end{align}
We emphasis that these assumptions,  quadratic  control cost and affine in control dynamics, are purely for the simplicity of treatment. These assumptions can be  omitted at the cost of increased notational   complexity.    \\
In the following subsections, we first characterize the performance of any feedback policy. Then, we use this characterization to provide $O(\epsilon^2)$ and $O(\epsilon^4)$ near-optimality results in the subsequent subsections.

%With a slight abuse of notation, let $\pi^{*} = \pi = (\pi_{t})^{T}_{t=1}$ be the optimal control law given by \eqref{cost_sto_orig} according to the dynamics \eqref{eq.0.1} and instantaneous cost \eqref{eq:cost1}.   

\subsection{Characterizing the Performance of a Feedback Policy}

Consider a noiseless  version of the system dynamics given by \eqref{eq.0.1}. We denote the ``nominal'' state trajectory as $\bar{x}_{t}$ and the ``nominal'' control as $\bar{u}_{t}$ where $u_{t} = \pi_{t}(x_{t})$, where $\pi = (\pi_{t})^{T-1}_{t=1}$ is a given control policy. The resulting dynamics  without noise is given by $\bar{x}_{t+1} = f(\bar{x}_t) + B_{t} \bar{u}_t$.   

%Let the associated nominal incremental costs be denoted by $\bar{c}_t = l(\bar{x}_t) + \frac{1}{2} \bar{u}_t'R\bar{u}_t$, and the associated nominal terminal cost be denoted by $\bar{c}_T = c_T(\bar{x}_T)$. 

Assuming that $f(\cdot)$ and $\pi_{t}(\cdot)$ are sufficiently smooth, we can linearize the dynamics about the nominal trajectory. Denoting  $\delta x_t = x_t - \bar{x}_t, \delta u_t = u_t - \bar{u}_t$, we can express, 
\begin{align}
\delta x_{t+1} &= A_t \delta x_t + B_t \delta u_t + S_t(\delta x_t) + \epsilon w_t, \label{eq.2}\\
\delta u_{t} &=  K_t \delta x_t + \tilde{S}_t(\delta x_t), \label{eq.3}
\end{align}
where $A_t = \frac{\partial f}{\partial x}|_{\bar{x}_t}$, $K_{t} = \frac{\partial \pi_{t}}{\partial x}|_{\bar{x}_t}$, and  $S_t(\cdot), \tilde{S}_t(\cdot)$ are second and higher order terms in the respective expansions. Similarly, we can linearize the instantaneous cost $c(x_{t}, u_{t})$ about the nominal values $(\bar{x}_{t}, \bar{u}_{t})$ as,
\begin{align}
c(x_t,u_t) &= {l}(\bar{x}_{t}) + L_t \delta x_t + H_t(\delta x_t) + \nonumber\\ 
&\hspace{1cm} \frac{1}{2}\bar{u}_t'R\bar{u}_t +  \delta u_t'R\bar{u}_t + \delta u_t'R\delta u_t,\label{eq.4}\\
c_{T}(x_{T}) &= {c}_{T}(\bar{x}_{T}) + C_T \delta x_T + H_T(\delta x_T),\label{eq.5}
\end{align}
where $L_t = \frac{\partial l}{\partial x}|_{\bar{x}_t}$, $C_T  = \frac{\partial c_T}{\partial x}|_{\bar{x}_t}$, and $H_t(\cdot)$ and $H_T(\cdot)$ are second and higher order terms in the respective expansions.

Using \eqref{eq.2} and \eqref{eq.3}, we can write the closed loop dynamics  of the trajectory $(\delta x_{t})^{T}_{t=1}$ as, 
\begin{align}
\label{eq.6}
\delta x_{t+1} = \underbrace{(A_t+B_tK_t)}_{\bar{A}_t} \delta x_t + \underbrace{\{B_t\tilde{S}_t(\delta x_t) + S_t(\delta x_t)\}}_{\bar{S}_t(\delta x_t)} + \epsilon w_t, 
\end{align} 
where $\bar{A}_t$ represents the linear part of the closed loop systems and the term $\bar{S}_t(.)$ represents the second and higher order terms in the closed loop system. Similarly, the closed loop incremental cost given in  \eqref{eq.4} can be expressed as
\begin{align}
 \label{eq.7}
c(x_t,u_t) = \underbrace{\{{l}(\bar{x}_{t}) + \frac{1}{2}\bar{u}_t'R\bar{u}_t\}}_{\bar{c}_t} + \underbrace{[L_t + \bar{u}_t'RK_t]}_{\bar{C}_t} \delta x_t \nonumber\\
+ \underbrace{(K_t\delta x_t+\tilde{S}_t(\delta x_t))'R(K_t\delta x_t + \tilde{S}_t(\delta x_t))}_{\bar{H}_t(\delta x_t)}.
\end{align}

Therefore, the cumulative cost of any given closed loop trajectory $(x_{t}, u_{t})^{T}_{t=1}$  can be expressed as,
\begin{align}
\label{eq.9a}
J^{\pi} &= \sum^{T-1}_{t=1}c(x_{t}, u_{t} = \pi_{t}(x_{t})) + c_{T}(x_{T}) \nonumber \\
&=\sum_{t=1}^T \bar{c}_t + \sum_{t=1}^T \bar{C}_t \delta x_t + \sum_{t=1}^T \bar{H}_t(\delta x_t),
\end{align}
where $\bar{c}_{T} = c_{T}(\bar{x}_{T}),  \bar{C}_{T} = C_{T}$.

We first show the following critical result. 
\begin{lemma} 
\label{L1}
Given any sample path, the state perturbation equation 
\[\delta x_{t+1} = \bar{A}_{t}  \delta x_{t} + \bar{S}_t(\delta x_t) + \epsilon w_t \]
given in  \eqref{eq.6} can be equivalently characterized  as
\begin{align}
\label{eq:mod-pert-1}
\delta x_{t}  = \delta x_t^l + e_t, ~ \delta x_{t+1}^l = \bar{A}_t \delta x_t^l + \epsilon w_t
\end{align} 
where $e_t$ is an $O(\epsilon^2)$ function that depends on the entire noise history $\{w_0,w_1,\cdots w_t\}$ and $\delta x_t^l$ evolves according to the  linear closed loop system. Furthermore, $e_t = e_t^{(2)} + O(\epsilon^3)$, where $e_t^{(2)} = \bar{A}_{t-1} e_{t-1}^{(2)} + \delta x_t^{l'}\bar{S}_{t-1}^{(2)} \delta x_t^l$, $e_0^{(2)} = 0$, and $\bar{S}_t^{(2)}$ represents the Hessian corresponding to the Taylor series expansion of the function $\bar{S}_t(.)$. 
\end{lemma}
\begin{proof}
We only consider the case when the state $x_t$ is scalar, the vector case is straightforward to derive and only requires a more complex notation.\\
We proceed by induction. The first general instance of the recursion occurs at $t=3$.
It can be shown that: 
 \begin{align}
 &\delta x_3 = \underbrace{(\bar{A}_2\bar{A}_1(\epsilon w_0) + \bar{A}_2 (\epsilon w_1) + \epsilon w_2)}_{\delta x_3^l} + \nonumber\\
 &\underbrace{\{\bar{A}_2 \bar{S}_1(\epsilon w_0) + \bar{S}_2(\bar{A}_1(\epsilon w_0) + \epsilon w_1 + \bar{S}_1(\epsilon w_0))\}}_{e_3}. 
 \end{align}
 Noting that $\bar{S}_1(.)$ and $\bar{S}_2(.)$ are second and higher order terms, it follows that $e_3$ is $O(\epsilon^2)$. \\
 Suppose now that $\delta x_t = \delta x_t^l + e_t$ where $e_t$ is $O(\epsilon^2)$. Then:
 \begin{align}
 \delta x_{t+1} &= \bar{A}_{t+1}(\delta x_t^l + e_t) + \epsilon w_t + \bar{S}_{t+1}(\delta x_t), \nonumber\\
 &= \underbrace{(\bar{A}_{t+1} \delta x_t^l + \epsilon w_t)}_{\delta x_{t+1}^l} +\underbrace{\{\bar{A}_{t+1}e_t + \bar{S}_{t+1}(\delta x_t)\}}_{e_{t+1}}.
 \end{align}
 Noting that $\bar{S}_{t}$ is $O(\epsilon^2)$ and that $e_{t}$ is $O(\epsilon^2)$ by assumption, the result follows that $e_{t+1}$ is $O(\epsilon^2)$. \\
 Now, let us take a closer look at the term $e_t$ and again proceed by induction. It is clear that $e_1 = e_1^{(2)} =0$. Next, it can be seen that $e_2 = \bar{A}_1e_1^{(2)} + S_1^{(2)}(\delta x_1^l)^2 + O(\epsilon^3) = \bar{S}_1^{(2)} (\epsilon \omega_0)^2 + O(\epsilon^3)$, which shows the recursion is valid for $t=2$ given it is so for $t=1$. \\
 Suppose that it is true for $t$. Then:
 \begin{align}
 \delta x_{t+1} &= \bar{A}_t \delta x_t + S_t(\delta x_t) + \epsilon \omega_t, \nonumber\\
 &= \bar{A}_t(\delta x_t^l + e_t) + S_t(\delta x_t^l + e_t) + \epsilon \omega _t, \nonumber\\
 &= \underbrace{(\bar{A}_t \delta x_t^l + \epsilon \omega_t)}_{\delta x_{t+!}^l}+ \underbrace{\bar{A}_t e_t^{(2)} + S_t^{(2)}(\delta x_t^l)^2}_{e_{t+1}^{(2)}} + O(\epsilon^3),
 \end{align}
 where the last line follows because $e_t = e_t^{(2)} + O(\epsilon^3)$, and $\bar{S}_t(.)$ contains second and higher order terms only. This completes the induction and the proof.
\end{proof}

%Using  \eqref{eq:mod-pert-1} in  \eqref{eq.9a}, we can obtain the cumulative cost of any given closed loop trajectory as,
%\begin{align}
%J^{\pi} = \underbrace{\sum_{t=1}^T \bar{c}_t }_{\bar{J}^{\pi}} + \underbrace{\sum_{t=1}^T \bar{C}_t \delta x_t^l}_{\delta J_1^{\pi}} + %\nonumber\\
%\underbrace{\sum_{t=1}^T \bar{H}_t(\delta x_t) + \bar{C}_t \bar{\bar{S}}_t}_{\delta J_2^{\pi}}. \label{eq.9b}
%\end{align}
Next, we have the following result for the expansion of the cost to go function $J^{\pi}$.
\begin{lemma}
\label{L2}
Given any sample path, the cost-to-go under a policy can be expanded as: 
\begin{equation} 
J^{\pi}  = \underbrace{\sum_t \bar{c}_t}_{\bar{J}^{\pi}} + \underbrace{\sum_t \bar{C}_t \delta x_t^l}_{\delta J_1^{\pi}} + \underbrace{\sum_t \delta x_t^{l'}\bar{H}_t^{(2)} \delta x_t^l+ \bar{C}_t e_t^{(2)}}_{\delta J_2^{\pi}} + O(\epsilon^3),
\end{equation}
where $\bar{H}_t^{(2)}$ denotes the second order coefficient of the Taylor expansion of $\bar{H}_t(.)$.
\end{lemma}
\begin{proof}
We have that:
\begin{align}
J^{\pi} &= \sum_t \bar{c}_t + \sum_t \bar{C}_t (\delta x_t^l + e_t) + \sum_t \bar{H}_t(\delta x_t^l + e_t), \nonumber\\
&= \sum_t \bar{c}_t + \sum_t \bar{C}_t \delta x_t^l + \sum_t \delta x_t^{l'}\bar{H}_t^{(2)} \delta x_t^l + \bar{C}_t e_t^{(2)} + O(\epsilon^3), \nonumber
\end{align}
where the last line of the equation above follows from an application of Lemma \ref{L1}.
\end{proof}

Now, we show the following important result.

\begin{proposition}  
\label{prop1} 
 \begin{align*}
\tilde{J}^{\pi} &= \mathbb{E}[J^{\pi} ] =  \bar{J}^{\pi} + O(\epsilon^2), \\
 \text{Var}(J^{\pi}) &=  \underbrace{\text{Var}(\delta J_{1}^{\pi})}_{O(\epsilon^{2})} +  O(\epsilon^4). 
 \end{align*}
 \end{proposition}
 
\begin{proof}
It is useful to first write the sample path cost in a slightly different fashion. It can be seen that given sufficient smoothness of the requisite functions, the cost of any sample path can be expanded as follows:
\begin{equation}
J^{\pi} = \bar{J}^{\pi} + \epsilon J_1^{\pi}+ \epsilon^2 J_2^{\pi} + \epsilon^3 J_3^{\pi} + \epsilon^4 J_4^{\pi} + \mathcal{R}, \nonumber\\
\end{equation}
where:
\begin{align}
J_1^{\pi} &= \mathcal{J}^1 \bar{\omega}, \nonumber\\
J_2^{\pi} &= \bar{\omega}' \mathcal{J}^2 \bar{\omega}, \nonumber
\end{align}
and so on for $J_3^{\pi}, J^{\pi}_4$ respectively, where $\mathcal{J}^i$ are constant matrices (tensors) of suitable dimensions, and $\bar{\omega} = [\omega_1, \cdots \omega_N]$. Further, the remainder function $\mathcal{R}$ is an $o(\epsilon^4)$ function in the sense that $\epsilon^{-4} \mathcal{R} \rightarrow 0$ as $\epsilon \rightarrow 0$. \\
Further, due to the whiteness of the noise sequences $\bar{\omega}$, it follows that $E[J_1^{\pi}] = 0$, and $E[J_3^{\pi}] = 0$, since these terms are made of odd valued products of the noise sequences, while $E[J_2^{\pi}], E[J_4^{\pi}]$ are both finite owing to the finiteness of the moments of the noise values and the initial condition. Further $lim_{\epsilon \rightarrow 0} \epsilon^{-4} E[\mathcal{R}] = E[\lim_{\epsilon} \epsilon^{-4} \mathcal{R}] = 0$, i.e., $E[\mathcal{R}]$ is $o(\epsilon^4)$. \\
Therefore, using Lemma \ref{L2}, and taking expectations on both sides, we obtain:
\begin{equation}
E[J^{\pi}] = \bar{J}^{\pi} + E[\epsilon J_1^{\pi}] + E[\epsilon^2 J_2^{\pi}] + O(\epsilon^4) = \bar{J}^{\pi} + O(\epsilon^2), \nonumber
\end{equation}
since $E[J_1^{\pi}] = 0$, and $E[\epsilon^ 2 J_2^{\pi,2}]$ is $O(\epsilon^2)$ due to the fact that $E[J_2^{\pi}] < \infty$.  \\
Next, using Lemma \ref{L2}, and taking the variances on both sides, and doing some work, we have:
\begin{align}
Var[J^{\pi}] &= Var[\epsilon J_1^{\pi}] + E[\epsilon J_1^{\pi}\epsilon^2J_2^{\pi}] + Var[\epsilon^2 J_2^{\pi}] + o(\epsilon^4) \nonumber\\
 &= Var[\delta J_1^{\pi}] + O(\epsilon^4),
\end{align}
where the second equality follows from the fact that  $E[\epsilon J_1^{\pi} \epsilon ^2J_2^{\pi}] = 0$ (proved in the appendix), and $Var[J_2^{\pi}]< \infty$. This completes the proof of the result.
%From Lemma \ref{L2}, we get,
%\begin{align} 
%\tilde{J}^{\pi} = \mathbb{E}[J^{\pi}] =  \mathbb{E}[ \bar{J}^{\pi} + \delta J_1^{\pi} + \delta J_2^{\pi}], \nonumber\\
%= \bar{J}^{\pi} + \mathbb{E}[\delta J_2^{\pi}]  = \bar{J}^{\pi} + O(\epsilon^2), \label{eq.10}
%= \bar{J}_0^{\pi} + \underbrace{ \mathbb{E}[\delta J_2^{\pi}]}_{\delta \tilde{J}_2^{\pi}} = \bar{J}^{\pi}_0 + O(\epsilon^2). \label{eq.10}
%\end{align}
%The first equality in the last line of the equations before follows from the fact that $\mathbb{E}[\delta x_t^l] = 0$, since its the linear part of the state perturbation driven by white noise and by definition $\delta x_1^l = 0$.The second equality follows form the fact that $\delta J_2^{\pi}$ is an $O(\epsilon^2)$ function.  Now,
%\begin{align}
%\text{Var}(J^{\pi}) = \mathbb{E}[ J^{\pi} - \tilde{J}^{\pi}]^2 \nonumber\\
%= \mathbb{E}[ \bar{J}_0^{\pi} + \delta J_1^{\pi} + \delta J_2^{\pi} - \bar{J}_0^{\pi} - \delta \tilde{J}_2^{\pi}]^2 \nonumber\\
%= \text{Var}(\delta J_1^{\pi}) + \text{Var}(\delta J_2^{\pi})  + 2 \mathbb{E}[\delta J_1^{\pi} \delta J_2^{\pi}].
%\end{align}
%Since $\delta J_2^{\pi}$ is $O(\epsilon^2)$, $\text{Var}(\delta J_2^{\pi})$ is an $O(\epsilon^4)$ function. It can be shown that $\mathbb{E}[\delta J_1^{\pi} \delta J_2^{\pi}]$ is $O(\epsilon^4)$ as well (proof is given \cite{d2cTR}). Finally $\text{Var}(\delta J_1^{\pi})$ is an $O(\epsilon^2)$ function because $\delta x^l$ is an $O(\epsilon)$ function. Combining these, we will get the desired result. 
\end{proof} 
A further consequence of the result above is the following. Suppose that given a policy $\pi_t(.)$, we only consider the linear part, i.e., the linear approximation $\pi_t^l(x_t) = \bar{u}_t + K_t\delta x_t$. However, according to Lemma \ref{L2}, the $\epsilon^2$ terms in the expansion of the cost of any sample path solely result from the linear closed loop system. Therefore, it follows that the sample path cost under the full policy $\pi_t(.)$ and the linear policy $\pi_t^l(.)$ agree up to the $\epsilon^2$ term. Therefore, it follows that $E[J^{\pi}] - E[J^{\pi^l}]  = O(\epsilon^4)$! We summarize this result in the following:

%\begin{comment}
\begin{proposition}
\label{prop2}
Let $\pi_t(.)$ be any given feedback policy. Let $\pi_t^l(x_t) = \bar{u}_t + K_t\delta x_t$ be the linear approximation of the policy. Then, the error in the expected cost to go under the two policies, $E[J^{\pi}] - E[J^{\pi^l}]  = O(\epsilon^4)$.
\end{proposition}

The above two results in Propositions \ref{prop1} and \ref{prop2} will form the basis of an $O(\epsilon^2)$ and an $O(\epsilon^4)$ decoupling result in the following subsections.
%\end{comment}

\subsection{An $O(\epsilon^2)$ Near-Optimal Decoupled Approach for Closed Loop Control} 

The following observations can now be made from Proposition \ref{prop1}. 
 
\begin{remark}[Expected cost-to-go]Recall that $u_{t} = \pi_t(x_t)$ $= \bar{u}_t + K_t\delta x_t + \tilde{S}_t(\delta x_t)$. However, note that due to Proposition \ref{prop1}, the expected cost-to-go, $\tilde{J}^{\pi}$, is determined almost solely (within $O(\epsilon^2)$)  by the nominal control action sequence $\bar{u}_t$. In other words, the linear and higher order feedback terms have only $O(\epsilon^2)$ effect on the expected cost-to-go function.
\end{remark}

\begin{remark}[Variance of cost-to-go] Given the nominal control action $\bar{u}_t$, the variance of the cost-to-go, which is $O(\epsilon^2)$, is determined overwhelmingly (within $O(\epsilon^4)$) by the linear feedback term $K_t \delta x_t$, i.e., by the variance of the linear perturbation of the cost-to-go, $\delta J_1^{\pi}$, under the linear closed loop system $\delta x_{t+1}^l = (A_t+B_tK_t)\delta x_t^l + \epsilon w_t$.
\end{remark}
Proposition \ref{prop1} and the remarks above suggest that an open loop control  super imposed with a closed loop control for the perturbed linear system may be approximately optimal. We delineate  this idea below. 

\textit{Open Loop Design.} First, we design an optimal (open loop) control sequence $\bar{u}^{*}_t$ for the noiseless system. More precisely, 
\begin{align}
\label{OL}
(\bar{u}^{*}_t)^{T-1}_{t=1} &= \arg \min_{(\tilde{u}_t)^{T-1}_{t=1}} \sum_{t=1}^{T-1} c(\bar{x}_t, \tilde{u}_t) + c_T(\bar{x}_T), \\
\bar{x}_{t+1} &= f(\bar{x}_t) + B_{t} \tilde{u}_t.  \nonumber
\end{align}
%We will discuss the details of this open loop design in Section \ref{sec4}.

\textit{Closed Loop Design.} We find the optimal feedback gain $K^{*}_t$ such that the variance of the linear closed loop system around the nominal path, $(\bar{x}_t, \bar{u}^{*}_t)$, from the open loop design above, is minimized.  
\begin{align}
(K^{*}_t)^{T-1}_{t=1} &=   \arg \min_{(K_t)^{T-1}_{t=1}} ~ \text{Var}(\delta J_1^{\pi}), \nonumber\\
 \delta J_1^{\pi} &= \sum_{t=1}^T \bar{C}_t  x_t^l,\nonumber\\
\delta x_{t+1}^l &= (A_t + B_tK_t) \delta x_t^l + \epsilon w_t. \label{CL}
\end{align} 
We now characterize the approximate closed loop policy below.

\begin{proposition}\label{propact3}
Construct a closed loop policy 
\begin{align}
\pi_t^*(x_t) = \bar{u}_t^* + K_t^*\delta x_t,
\end{align} 
where $\bar{u}_t^*$ is the solution of the open loop problem \eqref{OL}, and $K_t^*$ is the solution of the closed loop problem \eqref{CL}. Let $\pi^{o}$ be the optimal closed loop policy. Then, 
\begin{align}
 |\tilde{J}^{\pi*}  - \tilde{J}^{\pi^o}| = O(\epsilon^2).\nonumber
 \end{align}
 Furthermore, among all policies with nominal control action $\bar{u}_t^*$, the variance of the cost-to-go under policy $\pi_t^*$, is within $O(\epsilon^4)$ of the variance of the policy with the minimum variance.
 \end{proposition}
\begin{proof}
We have 
\begin{align*}
\tilde{J}^{\pi^*} - \tilde{J}^{\pi^o}  &= \tilde{J}^{\pi^*} -  \bar{J}^{\pi^*} + \bar{J}^{\pi^*} -  \tilde{J}^{\pi^o}  \\
&\leq \tilde{J}^{\pi^*} -  \bar{J}^{\pi^*} + \bar{J}^{\pi^{o}} -  \tilde{J}^{\pi^o}
\end{align*}
The  inequality above is due the fact that $\bar{J}^{\pi^*} \leq \bar{J}^{\pi^{o}}$, by definition of $\pi^{*}$. Now, using Proposition \ref{prop1}, we have that $|\tilde{J}^{\pi^*} - \bar{J}^{\pi^*}| = O(\epsilon^2)$, and $|\tilde{J}^{\pi^o} -  \bar{J}^{\pi^o}| = O(\epsilon^2)$. Also, by definition, we have $\tilde{J}^{\pi^o} \leq   \tilde{J}^{\pi^*}$. Then, from the above inequality, we get 
\begin{align*}
| \tilde{J}^{\pi^*} - \tilde{J}^{\pi^o} | \leq |\tilde{J}^{\pi^*} -  \bar{J}^{\pi^*} | + | \bar{J}^{\pi^{o}} -  \tilde{J}^{\pi^o} | = O(\epsilon^{2})
\end{align*}
A similar argument holds for the   variance as well. 
\end{proof}

Unfortunately, there is no standard solution to the closed loop problem \eqref{CL} due to the non additive nature of the cost function $\text{Var}(\delta J_1^{\pi})$. Therefore,  we solve a standard LQR problem as a surrogate, and the effect is again one of reducing the variance of the cost-to-go by reducing the variance of the closed loop trajectories.

\textit{Approximate Closed Loop Problem.} We solve the following LQR problem for suitably defined cost function weighting factors $Q_t$, $R_t$:
\begin{align}
\label{eq:acl1}
\min_{(\delta u_t)^{T}_{t=1}} ~&\mathbb{E} [\sum_{t=1}^{T-1} \delta x_t'Q_t\delta x_t + \delta u_t'R_t \delta u_t + \delta x_T'Q_T\delta x_t], \nonumber\\
&\delta x_{t+1} = A_t \delta x_t + B_t\delta u_t + \epsilon w_t.
\end{align}
The solution to the above problem furnishes us a feedback gan $\hat{K}_t^*$ which we can use in the place of the true variance minimizing gain $K_t^*$.  

\begin{remark}
Proposition \ref{prop1} states that the expected cost-to-go of the problem is dominated by the nominal cost-to-go. Therefore, even an open loop policy consisting of simply the nominal control action is within $O(\epsilon^2)$ of the optimal expected cost-to-go. However, the plan with the optimal feedback gain $K_t^*$ is strictly better than the open loop plan in that it has a lower variance in terms of the cost to go. Furthermore, solving the approximate closed loop problem using the surrogate LQR problem, we can expect a lower variance of the cost-to-go function as well.
%which is borne out empirically (albeit we cannot prove it).
\end{remark}

\subsection{An $O(\epsilon^4)$ Near-Optimal Decoupled Approach for Closed Loop Control}
In order to derive the results in this section, we need some additional structure on the dynamics. \textit{In essence, the results in this section require that the time discretization of the dynamics be small enough}. Thus, let the dynamics be given by:
\begin{equation}\label{SDE0}
x_t = x_{t-1} + \bar{f}(x_{t-1}) \Delta t + \bar{g}(x_{t-1}) u_t \Delta t + \epsilon \omega_t \sqrt{\Delta t},
\end{equation}
where $\omega_t$ is a white noise sequence, and the sampling time $\Delta t$ is small enough that $O(\Delta t^{\alpha})$ is negligible for $\alpha > 1$. The noise term above is a Brownian motion, and hence the $\sqrt{\Delta t}$ factor. Further, the incremental cost function $c(x,u)$ is given as:
$c(x,u) = \bar{l}(x) \Delta t + \frac{1}{2} u'\bar{R}u \Delta t.$
The main reason to use the above assumptions is to simplify the Dynamic Programming (DP) equation governing the optimal cost-to-go function of the system. The DP equation for the above system is given by:
\begin{equation}\label{DP}
J_t(x) = \min_{u_t} \{c(x,u) + E[J_{t+1}(x')]\},
\end{equation}
where $x' = x + \bar{f}(x)\Delta t + \bar{g}(x) u_t \Delta t + \epsilon \omega_t \sqrt{\Delta t}$ and $J_t(x)$ denotes the cost-to-go of the system given that it is at state $x$ at time $t$. The above equation is marched back in time with terminal condition $J_T(x) = c_T(x)$, and $c_T(.)$ is the terminal cost function. Let $u_t(.)$ denote the corresponding optimal policy.Then, it follows that the optimal control $u_t$ satisfies (since the argument to be minimized is quadratic in $u_t$)
\begin{equation}\label{opt_control}
u_t= -R^{-1}\bar{g}' J_{t+1}^x,
\end{equation}
where $J_{t+1}^x = \frac{\partial J_{t+1}}{\partial x}$.
 Further, let $u_t^d(.)$ be the optimal control policy for the deterministic system, i.e., Eq. \ref{SDE0} with $\epsilon =0$. The optimal cost-to-go of the deterministic system, $\phi_t(.)$ satisfies the deterministic DP equation:
\begin{equation}\label{detDP}
    \phi_t(x)= \min_{u}[c(x,u) + \phi_{t+1}(x')],
\end{equation}
where $x' = x + \bar{f}(x')\Delta t + \bar{g}(x')u\Delta t$. Then, identical to the stochastic case, $u_t^d = R^{-1}\bar{g}'\phi_t^x$.
Next, let $\varphi_t(.)$ denote the cost-to-go of the deterministic policy when applied to the stochastic system, i.e., $u_t^d$ applied to Eq. \ref{SDE0} with $\epsilon > 0$. The cost-to-go $\varphi_t(.)$ satisfies the policy evaluation equation:
\begin{equation}\label{detDP_eps}
   \varphi_t(x) = c(x, u_t^d(x)) + E[\varphi_{t+1}(x')],
\end{equation}
where now $x' = x+ \bar{f}(x)\Delta t + \bar{g}(x)u_t^d(x) \Delta t + \epsilon \omega_t \sqrt{\Delta t}$.
Note the difference between the equations \ref{detDP} and \ref{detDP_eps}.
Then, we have the following important result.
\begin{proposition}
\label{prop_eps}
The difference between the cost function of the optimal stochastic policy, $J_t$, and the cost function of the ``deterministic policy  applied to the stochastic system", $\varphi_t$, is $O(\epsilon^4)$, i.e. $|J_t (x) - \varphi_t (x)| = O(\epsilon^4)$ for all $(t,x)$.
\end{proposition}
The above result was originally proved in a seminal paper \cite{fleming1971stochastic} for continuous time, first passage problems. We have provided a simple derivation of the result, in the context of a discrete time finite horizon problem below. 
\begin{proof}
Using Proposition \ref{prop1}, we know that any cost function, and hence, the optimal cost-to-go function can be expanded as:
\begin{equation} \label{f1}
J_t(x) = J_t^0 + \epsilon^2 J_t^1 + \epsilon^4 J_t^2 + \cdots% \nonumber
\end{equation}
Thus, substituting the minimizing control in Eq. \ref{opt_control} into the dynamic programming Eq. \ref{DP} implies:
\begin{align}
J_t(x) = \bar{l}(x)\Delta t + \frac{1}{2} r(\frac{-\bar{g}}{r})^2 (J_{t+1}^x)^2 \Delta t + J_{t+1}^x \bar{f}(x)\Delta t  \nonumber\\
+ \bar{g}(\frac{-\bar{g}}{r})(J_{t+1}^x)^2 \Delta t+ \frac{\epsilon^2}{2} J_{t+1}^{xx} \Delta t + J_{t+1}(x), \label{f2}
\end{align}
where $J_t^{x}$, and $J_t^{xx}$ denote the first and second derivatives of the cost-to go function. Substituting Eq. \ref{f1} into eq. \ref{f2} we obtain that:
\begin{align}
(J_t^0 + \epsilon^2 J_t^1 + \epsilon^4 J_t^2+\cdots) = \bar{l}(x)\Delta t + \nonumber \\
\frac{1}{2}\frac{\bar{g}^2}{r}(J_{t+1}^{0,x} + \epsilon^2 J_{t+1}^{1,x}+\cdots)^2 \Delta t  \nonumber \\
+(J_{t+1}^{0,x}+ \epsilon^2 J_{t+1}^{1,x}+\cdots) \bar{f}(x) \Delta t \nonumber \\
- \frac{\bar{g}^2}{r} (J_{t+1}^{0,x}+ \epsilon^2 J_{t+1}^{1,x}+\cdots)^2\Delta t \nonumber\\
+ \frac{\epsilon^2}{2} (J_{t+1}^{0,x}+ \epsilon^2 J_{t+1}^{1,x}+\cdots)\Delta t + J_{t+1}(x). \label{f3}
\end{align}
Now, we equate the $\epsilon^0$, $\epsilon^2$ terms on both sides to obtain perturbation equations for the cost functions $J_t^0, J_t^1, J_t^2 \cdots$. \\
First, let us consider the $\epsilon^0$ term. Utilizing Eq. \ref{f3} above, we obtain:
\begin{align} \label{f4}
J_t^0  = \bar{l} \Delta t + \frac{1}{2} \frac{\bar{g}^2}{r}(J_{t+1}^{0,x})^2 \Delta t + %\nonumber\\
\underbrace{(\bar{f} + \bar{g}\frac{-\bar{g}}{r} J_t^{0,x})}_{\bar{f}^0} J_t^{0,x}\Delta t + J_{t+1}^0, 
\end{align}
with the terminal condition $J_T^0 = c_T$, and where we have dropped the explicit reference to the argument of the functions $x$ for convenience. \\
Similarly, one obtains by equating the $O(\epsilon^2)$ terms in Eq. \ref{f3} that:
\begin{align}
J_t^1  = \frac{1}{2} \frac{\bar{g}^2}{r} (2J_{t+1}^{0,x}J_{t+1}^{1,x})\Delta t + J_{t+1}^{1,x} \bar{f}\Delta t -\frac{\bar{g}^2}{r} (2J_{t+1}^{0,x}J_{t+1}^{1,x})\Delta t + \nonumber
\frac{1}{2} J_{t+1}^{0,xx} \Delta t + J_{t+1}^1,
\end{align}
which after regrouping the terms yields:
\begin{align} \label{f5}
J_t^1 = \underbrace{(\bar{f}+ \bar{g} \frac{-\bar{g}}{r} J_{t+1}^{0,x})J_{t+1}^{1,x} }_{= \bar{f}^0}\Delta t + \frac{1}{2} J_{t+1}^{0,xx} \Delta t + J_{t+1}^1,
\end{align}
with terminal boundary condition $J_T^1 = 0$.
Note the perturbation structure of Eqs. \ref{f4} and \ref{f5}, $J_t^0$ can be solved without knowledge of $J_t^1, J_t^2$ etc, while $J_t^1$ requires knowledge only of $J_t^0$, and so on. In other words, the equations can be solved sequentially rather than simultaneously.\\

Now, let us consider the deterministic policy $u_t^d(.)$ that is a result of solving the deterministic DP equation:
\begin{equation}
\phi_t(x) = \min_{u} [c(x,u) + \phi_{t+1}(x')],
\end{equation}
where $x' = x + \bar{f}\Delta t + \bar{g} u \Delta t$, i.e., the deterministic system obtained by setting $\epsilon =0$ in Eq. \ref{SDE0}, and $\phi_t$ represents the optimal cost-to-go of the deterministic system. Analogous to the stochastic case, $u_t^d = \frac{-\bar{g}}{r} \phi_t^x$.
Next, let $\varphi_t$ denote the cost-to-go of the deterministic policy $u_t^d(.)$ \textit{when applied to the stochastic system, i.e., Eq. \ref{SDE0} with $\epsilon >0$}. Then, the cost-to-go of the deterministic policy, when applied to the stochastic system, satisfies:
\begin{equation}
\varphi_t = c(x, u_t^d(x)) + E[\varphi_{t+1}(x')],
\end{equation}
where $x' = \bar{f}\Delta t + \bar{g}u_t^d \Delta t + \epsilon\sqrt{\Delta t}\omega_t$. Substituting $u_t^d(.) = \frac{-\bar{g}}{r}\phi_t^x$ into the equation above implies that:
\begin{eqnarray}
\varphi_t = \varphi_t^0 + \epsilon^2\varphi_t^1 + \epsilon^4 \varphi_t^2+ \cdots \nonumber\\
= \bar{l} \Delta t + \frac{1}{2}\frac{\bar{g}^2}{r}(\phi_{t+1}^x)^2 \Delta t \nonumber
+ (\varphi_{t+1}^{0,x} + \epsilon^2\varphi_{t+1}^{1,x} + \cdots )\bar{f} \Delta t \nonumber\\
+  \bar{g}\frac{-\bar{g}}{r} \phi_{t+1}^x (\varphi_{t+1}^{0,x} + \epsilon^2\varphi_{t+1}^{1,x} + \cdots ) \Delta t \nonumber\\
+ \frac{\epsilon^2}{2} (\varphi_{t+1}^{0,xx} + \epsilon^2\varphi_{t+1}^{1,xx} + \cdots ) \Delta t \nonumber\\
+  (\varphi_{t+1}^0 + \epsilon^2\varphi_{t+1}^1 + \cdots). \label{f6}
\end{eqnarray}
As before, if we gather the terms for $\epsilon^0$, $\epsilon^2$ etc. on both sides of the above equation, we shall get the equations governing $\varphi_t^0, \varphi_t^1$ etc.  First, looking at the $\epsilon^0$ term in Eq. \ref{f5}, we obtain:
\begin{align} \label{f7}
\varphi_t^0 = \bar{l} \Delta t + \frac{1}{2} \frac{\bar{g}^2}{r}(\phi_{t+1}^x)^2 \Delta t  %\nonumber\\
+ (\bar{f} + \bar{g} \frac{-\bar{g}}{r}\phi_{t+1}^x )\varphi_{t+1}^{0,x}\Delta t + \varphi_{t+1}^0,
\end{align}
with the terminal boundary condition $\varphi_T^0 = c_T$. However, the deterministic cost-to-go function also satisfies:
\begin{align}\label{f8}
\phi_t = \bar{l} \Delta t + \frac{1}{2} \frac{\bar{g}^2}{r} (\phi_{t+1}^x)^2 \Delta t  %\nonumber\\
+( \bar{f}   + \bar{g} \frac{-\bar{g}}{r} \phi_{t+1}^x) \phi_{t+1}^x \Delta t + \phi_{t+1},
\end{align}
with terminal boundary condition $\phi_T = c_T$. Comparing Eqs. \ref{f7} and \ref{f8}, it follows that $\phi_t = \varphi_t^0$ for all $t$. Further, comparing them to Eq. \ref{f4}, it follows that $\varphi_t^0 = J_t^0$, for all $t$. Also, note that the closed loop system above, $\bar{f}   + \bar{g} \frac{-\bar{g}}{r} \phi_{t+1}^x = \bar{f}^0$ (see Eq. \ref{f4} and \ref{f5}). \\

Next let us consider the $\epsilon^2$ terms in Eq. \ref{f6}. We obtain:
\begin{align}
\varphi_t^1 = \bar{f} \varphi_{t+1}^{1,x} \Delta t + \bar{g}\frac{-\bar{g}}{r} \phi_{t+1}^x \varphi_{t+1}^{1,x} \Delta t + \nonumber
\frac{1}{2}\varphi_{t+1}^{0,xx} + \varphi_{t+1}^1.
\end{align}
Noting that $\phi_t = \varphi_t^0 $, implies that (after collecting terms):
\begin{align} \label{f9}
\varphi_t^1 = \bar{f}^0 \varphi_{t+1}^{1,x} \Delta t + \frac{1}{2}\varphi_{t+1}^{0,xx} \Delta t + \varphi_{t+1}^1,
\end{align}
with terminal boundary condition $\varphi_N^1 = 0$. Again, comparing Eq. \ref{f9} to Eq. \ref{f5}, and noting that $\varphi_t^0 = J_t^0$,  it follows that $\varphi_t^1 = J_t^1$, for all $t$. This completes the proof of the result. 
\end{proof}

Given some initial condition $x_0$, consider a linear truncation of the optimal deterministic policy, i.e., let $u_t^l(.) = \bar{u}_t + K_t\delta x_t$, where the deterministic policy is given by $u_t^d = \bar{u}_t + K_t\delta x_t + S_t(\delta x_t)$, where $S_t(.)$ denote the second and higher order terms in the optimal deterministic feedback policy.  Using Proposition \ref{prop2}, it follows that the cost of the linear policy, say $\varphi_t^l(.)$,  is within $O(\epsilon^4)$ of the cost of the deterministic policy $u_t^d(.)$, when applied to the stochastic system in Eq. \ref{SDE0}. However, the result in Proposition \ref{prop_eps} shows that the cost of the deterministic policy is within $O(\epsilon^4)$ of the optimal stochastic policy. Taken together, this implies that the cost of the linear deterministic policy is within $O(\epsilon^4)$ of the optimal stochastic policy. This may be summarized in the following result.
\begin{proposition} \label{e4}
Let the optimal cost function under the true stochastic policy be given $J_t(.)$ Let the optimal deterministic policy be given by $u_t^d(x_t) = \bar{u}_t + K_t\delta x_t + S_t(\delta x_t)$, and the linear approximation  to the policy be $u_t^l(x_t) = \bar{u}_t + K_t\delta x_t$, and let the cost of the linear policy be given by $\varphi_t^l(x)$. Then $|J_t(x) - \varphi_t^l(x)| = O(\epsilon^4)$ for all $(t,x)$.
\end{proposition}
Now, it remains to be seen how to design the $\bar{u}_t$ and the linear feedback term $K_t$. 
The open loop optimal control sequence $\bar{u}_t$ is found identically to the previous section. However, the linear feedback gain $K_t$ is calculated in a slightly different fashion and may be done as shown in the following result. In the following, $\mathcal{F}(x) = x + \bar{f}(x)\Delta t$, $\mathcal{G}(x) = \bar{g}(x)\Delta t$, $A_t = \frac{\partial \mathcal{F}}{\partial x}|_{\bar{x}_t} + \frac{\partial \mathcal{G}\bar{u}_t}{\partial x}|_{\bar{x}_t}$, $B_t = \mathcal{G}(\bar{x}_t)$, $L_t = \frac{\partial l}{\partial x}|_{\bar{x}_t}'$ and $L_{tt} = \nabla^2_{xx} l |_{\bar{x}_t}$. Let $\phi_t(x_t)$ denote the optimal cost-to-go of the detrministic problem, i.e., Eq \ref{SDE0} with $\epsilon = 0$. \\
\begin{proposition} \label{T-PFC}
\textbf{Decoupled Design.}
Given an optimal nominal trajectory $(\bar{x}_t, \bar{u}_t)$, the backward evolutions of the first and second derivatives, $G_t = \frac{\partial \phi_t}{\partial x}|_{\bar{x}_t}'$ and $P_t = \nabla^2_{xx} \phi_t|_{\bar{x}_t}$, of the optimal cost-to-go function ${\phi}_t({ x_t})$, initiated with the terminal boundary conditions $G_N = \frac{\partial {c}_N({x_N})}{\partial { x_N}}\arrowvert_{\bar{x}_N}' $ and $P_N = \nabla^2_{x} c_N\arrowvert_{\bar{x}_N}$ respectively, are as follows: \\
\begin{equation}
\hspace*{-5.72cm}G_t = L_t + G_{t+1}A_t, \label{OLprop6}
\end{equation}
\begin{equation}
\hspace*{-1.5cm}P_t = L_{tt} + A_t' P_{t+1} A_t - K'_t S_t K_t + G_{t+1} \otimes \tilde R_{t,xx},  \label{feedback}
\end{equation}
for $t = \{0,1,...,N-1\}$, where, $S_t = (R_t + B_t' P_{t+1}B_t), K_t = -S_t^{-1}(B_t' P_{t+1} A_t + (G_{t+1} \otimes \tilde{R}_{t_{x u}})'), \tilde{ R}_{t,xx} = \nabla^{2}_{xx}\mathcal{F}({ x_t})\arrowvert_{{\bar{x}_t}} + \nabla^{2}_{xx}\mathcal{G}({x_t})\arrowvert_{{\bar{x}_t}, {\bar{u}_t}}, \tilde{ R}_{t,xu} = \nabla^{2}_{xu}(\mathcal{F}({x_t}) + \mathcal{G}({x_t}){{u}_t}) \arrowvert_{\bar{x}_t, \bar{u}_t}$ where $\nabla^2_{xx}$ represents the Hessian of a vector-valued function w.r.t $x$ and $\otimes$ denotes the tensor product. 
\end{proposition}
\begin{proof}
Consider the Dynamic Programming equation for the deterministic cost-to-go function:
\begin{align*}
{\phi}_t({x_t}) = \mathop{min}_{{u_t}} Q_t({ x_t},{u_t}) = \mathop{min}_{{u_t}}\{c_t({x_t}, {u_t}) + {\phi}_{t+1}({x_{t+1}})\}
\end{align*}
By Taylor's expansion about the nominal state at time $t+1$,
\begin{align*}
{\phi}_{t+1}({x_{t+1}})=&{\phi}_{t+1}({ \bar{x}_{t+1}}) + G_{t+1} \delta {x_{t+1}} \\
&+ \frac{1}{2}\delta {x_{t+1}}' P_{t+1} \delta {x_{t+1}} + q_{t+1}(\delta {x_{t+1}}).
\end{align*}
Substituting the linearization of the dynamics, \(\delta {x_{t+1}} = A_t \delta { x_t} + B_t \delta {u_t} + r_t(\delta {x_t}, \delta {u_t})\) in the above expansion, 
\begin{align*}
%\begin{split}
& {\phi}_{t+1}({x_{t+1}}) = {\phi}_{t+1}({ \bar{x}_{t+1}}) +  G_{t+1} (A_t \delta { x_t} + B_t \delta { u_t} + r_t(\delta {x_t}\\
&, \delta {u_t}) ) + ( A_t \delta {x_t} + B_t \delta {u_t} + r_t(\delta {x_t}, \delta {u_t}))' P_{t+1} (A_t \delta {x_t} \\
&+ B_t \delta {u_t} +  r_t(\delta {x_t}, \delta {u_t}) )   + q_{t+1}(\delta {x_{t+1}}).  
%\end{split}
\end{align*}
Similarly, expand the incremental cost at time $t$ about the nominal state,
\begin{align*}
c_t({x_t}, {u_t}) = \bar{l}_t + L_t \delta {x_t} + \frac{1}{2} \delta {x_t}' L_{tt} \delta {x_t} + \frac{1}{2} \delta {u_t}' R_t {\bar{u}_t} \\
+ \frac{1}{2}  {\bar{u}_t}' R_t \delta {u_t} + \frac{1}{2} \delta {u_t}' R_t \delta {u_t}  + \frac{1}{2}{\bar{u}_t}' R_t {\bar{u}_t} + s_t(\delta {x_t}).
\end{align*}
\begin{align*}
&Q_t({x_t},{u_t}) = \overbrace{[\bar{l}_t + \frac{1}{2} {\bar{u}_t}^\intercal R_t { \bar{u}_t} + {\phi}_{t+1}({\bar{x}_{t+1}}) ]}^{\bar{\phi}_t({\bar{x}_t}, \bar{u}_t)}\\
&+ \delta {u_t}'(B_t' \frac{P_{t+1}}{2} B_t + \frac{1}{2} R_t) \delta {u_t} + 
\delta {u_t}'(B_t' \frac{P_{t+1}}{2} A_t \delta {x_t} \\
&+ \frac{1}{2} R_t {\bar{u}_t} +B_t' \frac{P_{t+1}}{2}r_t) + (\delta { x_t}' A_t' \frac{P_{t+1}}{2}B_t 
+ \frac{1}{2} {\bar{u}_t} R_t \\
&+r_t' \frac{P_{t+1}}{2}B_t + G_{t+1}B_t) \delta {u_t} + \delta { x_t}' A_t' \frac{P_{t+1}}{2}A_t \delta {x_t} \\
&+ \delta {x_t}' \frac{P_{t+1}}{2}A_t' r_t+ (r_t' \frac{P_{t+1}}{2}A_t + G_{t+1} A_t) \delta {x_t} \\
&+ r_t' \frac{P_{t+1}}{2}r_t+ G_{t+1}r_t + q_t \equiv \bar{\phi}_t(\bar{x}_t, \bar{u}_t) + H_t(\delta x_t, \delta u_t).
\end{align*}
\begin{flalign*}
\text{Now,}
\mathop{min}_{{u_t}} Q_t({x_t}, {u_t}) &=  \mathop{min}_{{\bar{u}_t}}  \bar{\phi}_t({\bar{x}_t}, { \bar{u}_t}) + \mathop{min}_{\delta {u_t}} H_t(\delta {x_t} ,\delta {u_t})
\end{flalign*}
\textbf{First order optimality:} Along the optimal nominal control sequence \(\bar{u}_t\), it follows from the minimum principle that
\begin{align*}
\frac{\partial c_t({x_t}, {u_t})}{\partial { u_t}} + \frac{\partial g({x_t}) }{\partial { u_t}}' \frac{\partial {\phi}_{t+1}({ x_{t+1}})}{\partial {x_{t+1}}}  = 0   
\end{align*}
\begin{equation}
\Rightarrow R_t {\bar{u}_t} + B_{t}' G_{t+1}' = 0
\end{equation}
By setting \(\frac{\partial H_t(\delta {x_t}, \delta {u_t})}{\partial \delta {u_t}} = 0  \), we get:
\begin{align*}
 \delta {u^{*}_t} &=- S_t^{-1} (R_t { \bar{u}_t} + B_t' G_{t+1}') - S_t^{-1}(B_t' P_{t+1}A_t + \\
 & (G_t \otimes \tilde R_{t,xu})') \delta {x_t}-S_t^{-1}(B_t' P_{t+1}r_t) \\
 &= \underbrace{- S_t^{-1}(B_t' P_{t+1}A_t + (G_{t+1} \otimes \tilde R_{t,xu})')}_{K_t} \delta {x_t}  \\
 &+\underbrace{S_t^{-1}(-B_t' P_{t+1}r_t)}_{p_t}
 \end{align*}
 where, \(S_t = R_t + B_t' P_{t+1}B_t.\)
 \begin{align*}
    \Rightarrow{\delta {u_t}= K_t \delta {x_t} + p_t}.
 \end{align*}

\noindent Substituting it in the expansion of $J_t$ and regrouping the terms based on the order of \(\delta x_t \) (till $2^{nd}$ order), we obtain:
\begin{align*}
{\phi}_t({x_t}) = \bar{\phi}_t({\bar{x}_t}) + (L_t + (R_t {\bar{u}_t} + B_t' G_{t+1}')K_t + G_{t+1}A_t)\delta {x_t} \\
+ \frac{1}{2}\delta {x_t}' (L_{tt} + A_t' P_{t+1} A_t - K_t' S_t K_t + 
G_{t+1}\otimes{\tilde R}_{t,xx}) \delta {x_t }.
\end{align*}

Expanding the LHS about the optimal nominal state result in the recursive equations in Proposition \ref{T-PFC}.
\end{proof}

\subsection{Summary of the Decoupling Results and Implications}
The previous two subsections showed that the feedback parameterization can be written as: $\pi_t(x_t) = \bar{u}_t + K_t \delta x_t$, where $\delta x_t = x_t -\bar{x}_t$ denotes the state deviation from the nominal. Further, it was shown that the optimal open loop sequence $\bar{u}_t$ is independent of the feedback gain, while the feedback gain $K_t$ can be designed based on the optimal $\bar{u}_t$. \textit{Hence, the term decoupling, in the sense that the search for the optimal parameter $(\bar{u}_t^*, K_t^*)$ need not be done jointly.}\\
Moreover, it was shown that depending on how one designed the gain $K_t$, we can obtain either $O(\epsilon^2)$ (Proposition \ref{propact3}), or $O(\epsilon^4)$ (Propositions \ref{T-PFC}), near-optimality to the true stochastic policy.\\

\section{Analysis of the High Noise Regime}
In this section, we perform a rudimentary analysis of the high noise regime. The medium noise case is more difficult to analyze and is left for future work, along with a more sophisticated treatment of the high noise regime.\\
%\textit{High Noise Regime.} 
First, recall the Dynamic Programming (DP) equation for the backward pass to determine the optimal time varying feedback policy:
\begin{equation} \label{DP}
    J_t(\V{x}_t) = \min_{\V{u}_t}\left\{c(\V{x}_t,\V{u}_t) + \Exp{}{J_{t+1}(\V{x}_{t+1})}\right\},
\end{equation}
where $J_t(\V{x}_t)$ denotes the cost-to-go at time $t$ given the state is $\V{x}_t$, with the terminal condition $J_T(\cdot) = c_T(\cdot)$ where $c_T$ is the terminal cost function, and the next state $\V{x}_{t+1} = f(\V{x}_t) + \M{B}_t(\V{u}_t + \epsilon \V{w}_t)$. Suppose now that the noise is so high that $\V{x}_{t+1} \approx \M{B}_t \epsilon \V{w}_t$, i.e., the dynamics are completely swamped by the noise.\\
Consider now the expectation $\Exp{}{c_T(\V{x}_{t+1})}$ given some control $\V{u}_t$ was taken at state $\V{x}_t$. Since $\V{x}_{t+1}$ is determined entirely by the noise, $\Exp{}{c_T(\V{x}_{t+1})} = \int c_T(\M{B}_t\epsilon \V{w}_t)\mathbf{p}(\V{w}_t) d\V{w}_t = \overline{c_T}$, where $\overline{c_T}$ is a constant regardless of the previous state and control pair $\V{x}_t, \V{u}_t$. This observation holds regardless of the function $c_T(\cdot)$ and the time $t$.\\
Next, consider the DP iteration at time $T-1$. Via the argument above, it follows that $\Exp{}{J_T(\V{x}_T)}=\Exp{}{c_T(\V{x}_T)} = \overline{c_T}$, regardless of the state control pair $\V{x}_{T-1},\V{u}_{T-1}$ at the $(T-1)^{th}$ step, 
and thus, the minimization reduces to 
$J_{T-1}(\V{x}_{T-1}) = \min_{\V{u}} \left\{c(\V{x}_{T-1},\V{u})  + \overline{c_T}\right\}$, 
and thus, the minimizer is just the greedy action $\V{u}^*_{T-1} = \argmin_{\V{u}} c(\V{x}_{T-1},\V{u})$ due to the constant bias $\overline{c_T}$. 
The same argument holds for any $t$ since, although there might be a different $J_{t}(\cdot)$ at every time $t$, the minimizer is still the greedy action that minimizes $c(\V{x}_t,\V{u})$ as the cost-to-go from the next state is averaged out to simply some $\bar{J}_{t+1}$.\\

\section{Additional Simulation Results:}
In addition to the simulations performed on the car-like robot shown in the main article, experiments are performed on a car with trailers and a quadrotor whose results are shown in Figure~\ref{fig:cost_trailer} and \ref{fig:cost_quad} respectively. We also show a scenario on a car-like robot in an environment with obstacles to illustrate that the decoupling approach can handle such cases. The parameters used in the simulations are given in Table~\ref{table:params1} and \ref{table:params3}. As seen in car-like robot, the performance of T-LQR2 is close to MPC for a wide range of noise levels. It is also evident from the replanning operations plots in Figure~\ref{fig:trailer_replan} and \ref{fig:quad_replan} that T-LQR2 is computationally efficient when compared to MPC. MPC-SH also exhibits the similar trends as shown in the main article.  

\subsection{Car-like robot with trailers:}
\begin{figure}[h]
    \centering
    \subfloat[Full noise spectrum]{
        \includegraphics[width=.30\textwidth]{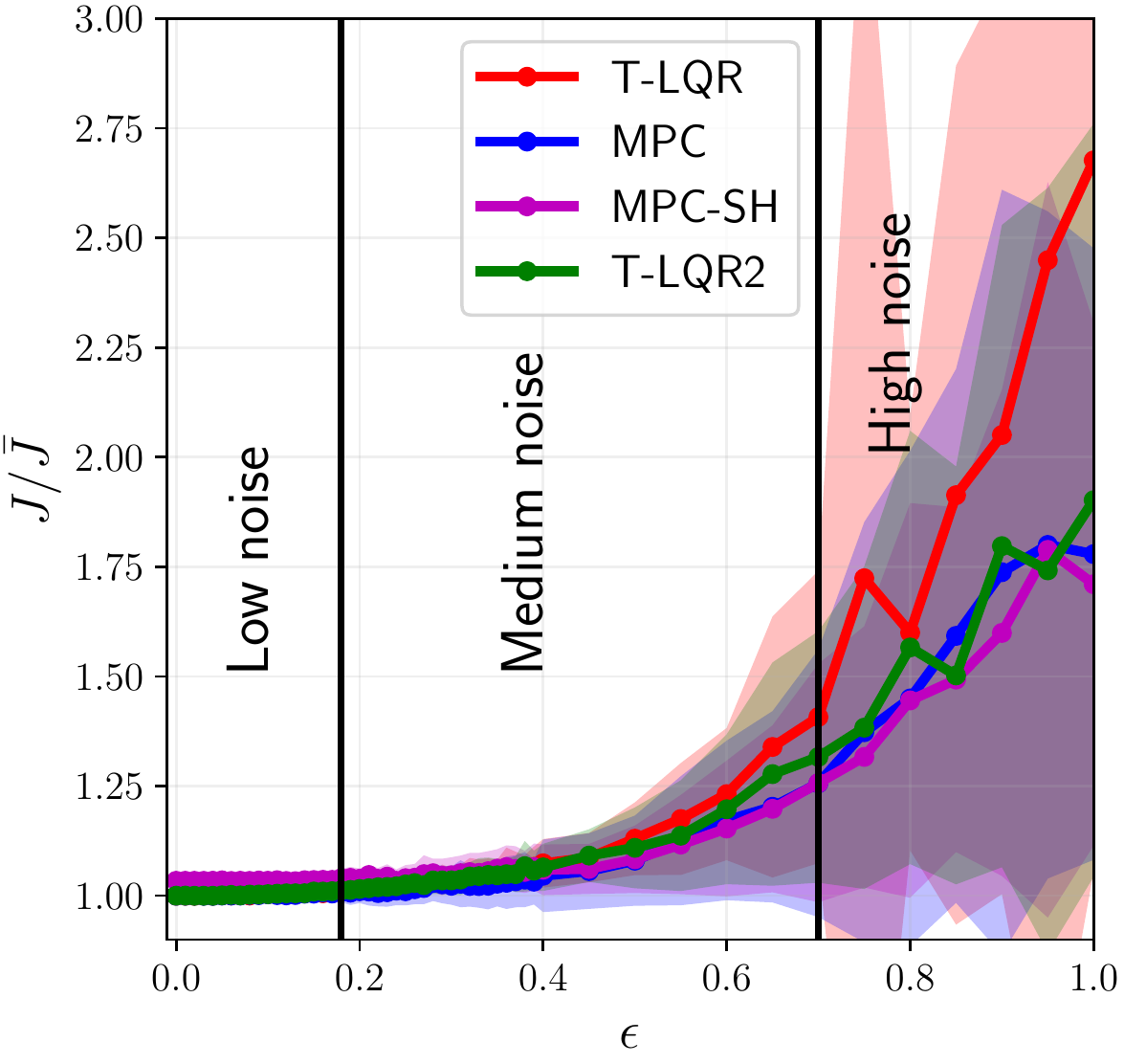}
        \label{trailer_full}}
    \subfloat[Enhanced detail $0 \leq \epsilon \leq 0.4$]{
        \includegraphics[width=.30\textwidth]{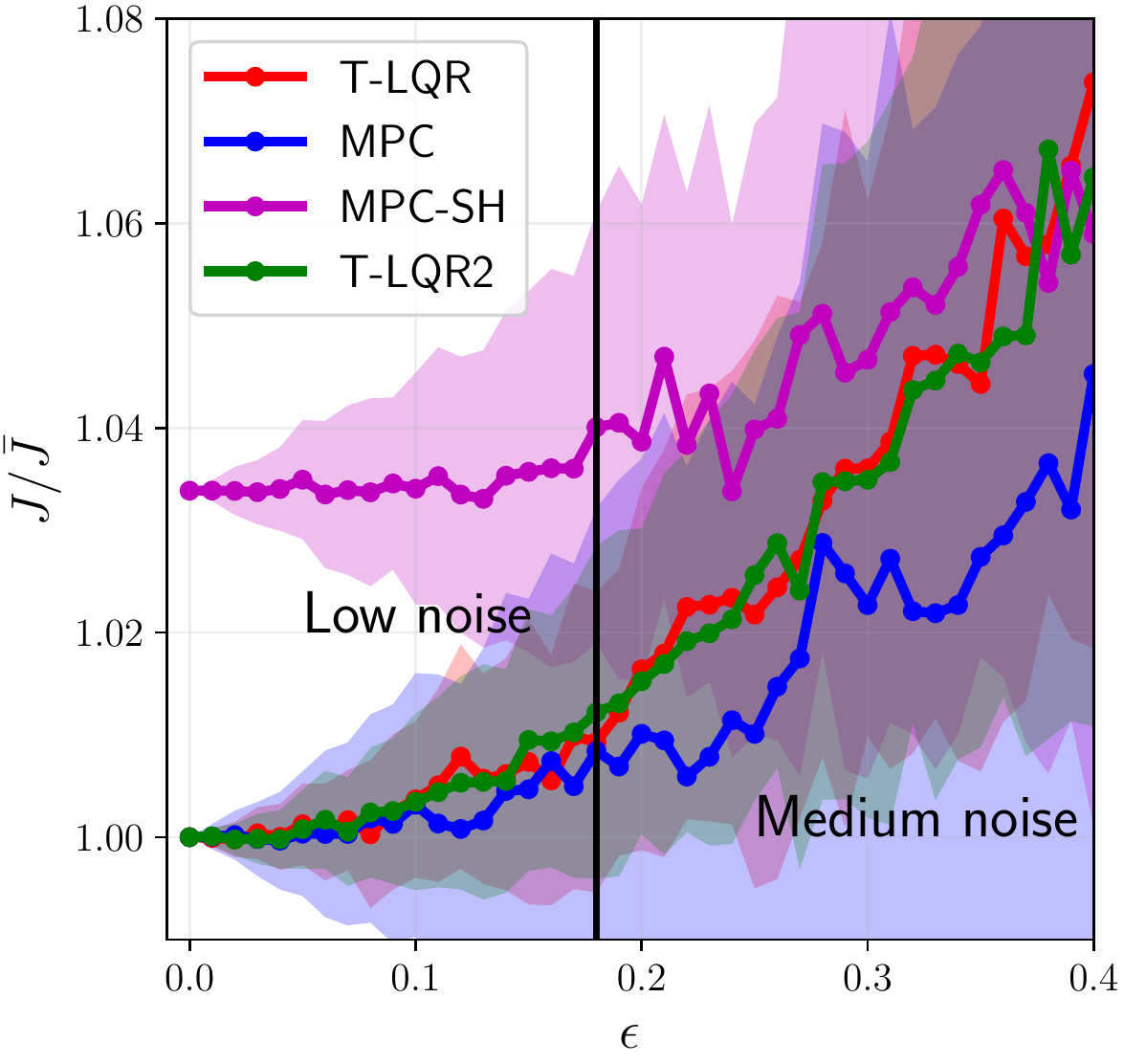}
        \label{trailer_low}}
    \subfloat[Replanning operations]{
    \includegraphics[width=.30\textwidth]{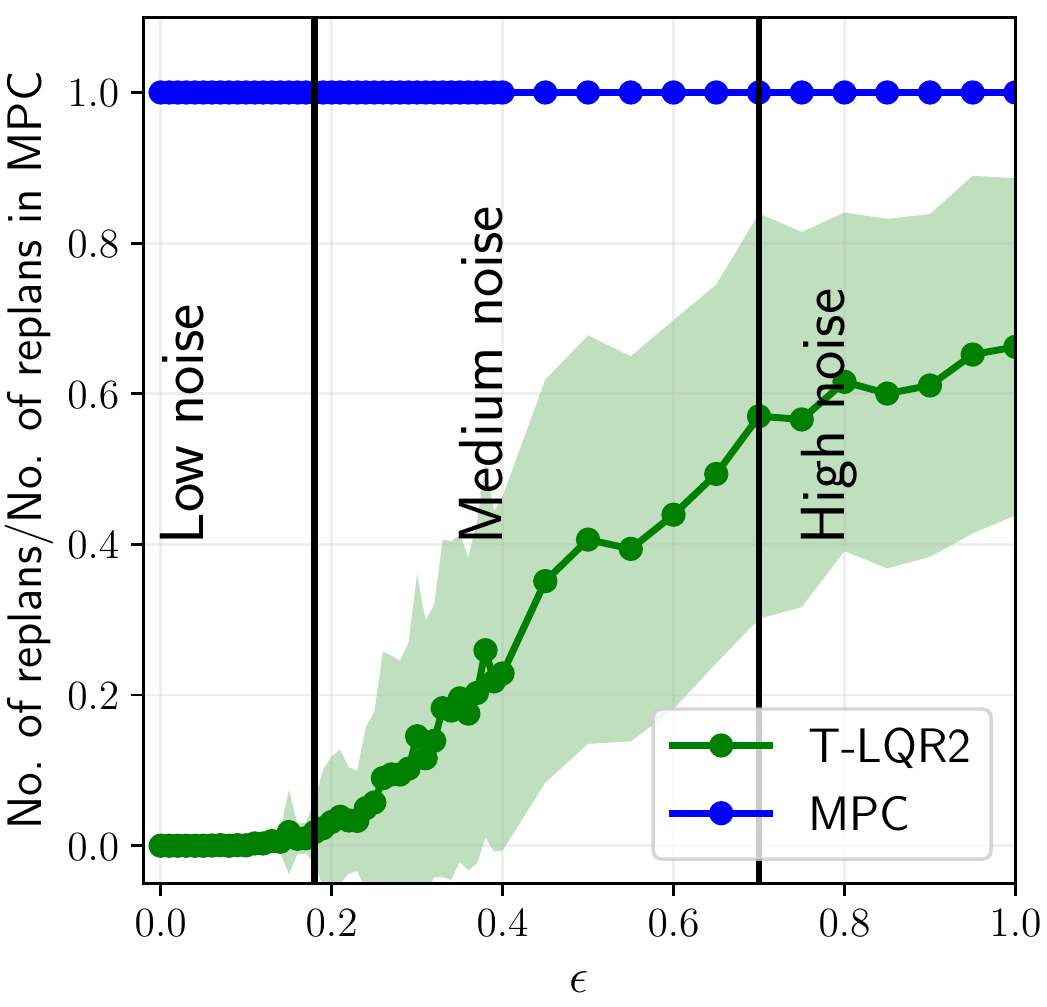}
    \label{fig:trailer_replan}}
    \caption{Cost evolution of the different algorithms for a car with 2 trailer system.}
    \label{fig:cost_trailer}
\end{figure}

Having trailers in a car-like robot makes it more complex by increasing the state dimension of it by the number of trailers attached. Here we consider 2 trailers whose heading angles are given by, 

\begin{align*}
    \theta_{1}(t+1) &= \theta_{1}(t) + \frac{v_t}{L}sin(\theta (t) - \theta_1 (t)) \Delta t, \\
    \theta_{2}(t+1) &= \theta_{2}(t) + \frac{v_t}{L}cos(\theta (t) - \theta_1 (t))sin(\theta_1 (t) - \theta_2 (t)) \Delta t.
\end{align*}
The performance is shown in Figure~\ref{fig:cost_trailer}. As seen in the car-like robot, T-LQR is near-optimal in the low noise regime, while T-LQR2 performs similar to MPC in the medium and high noise regime. In the high noise regime, as seen earlier, MPC-SH achieves similar performance to MPC and T-LQR2 despite planning only for a short horizon.
\subsection{Quadrotor:}
\begin{figure}[h]
    \centering
    \subfloat[Full noise spectrum]{
        \includegraphics[width=.30\textwidth]{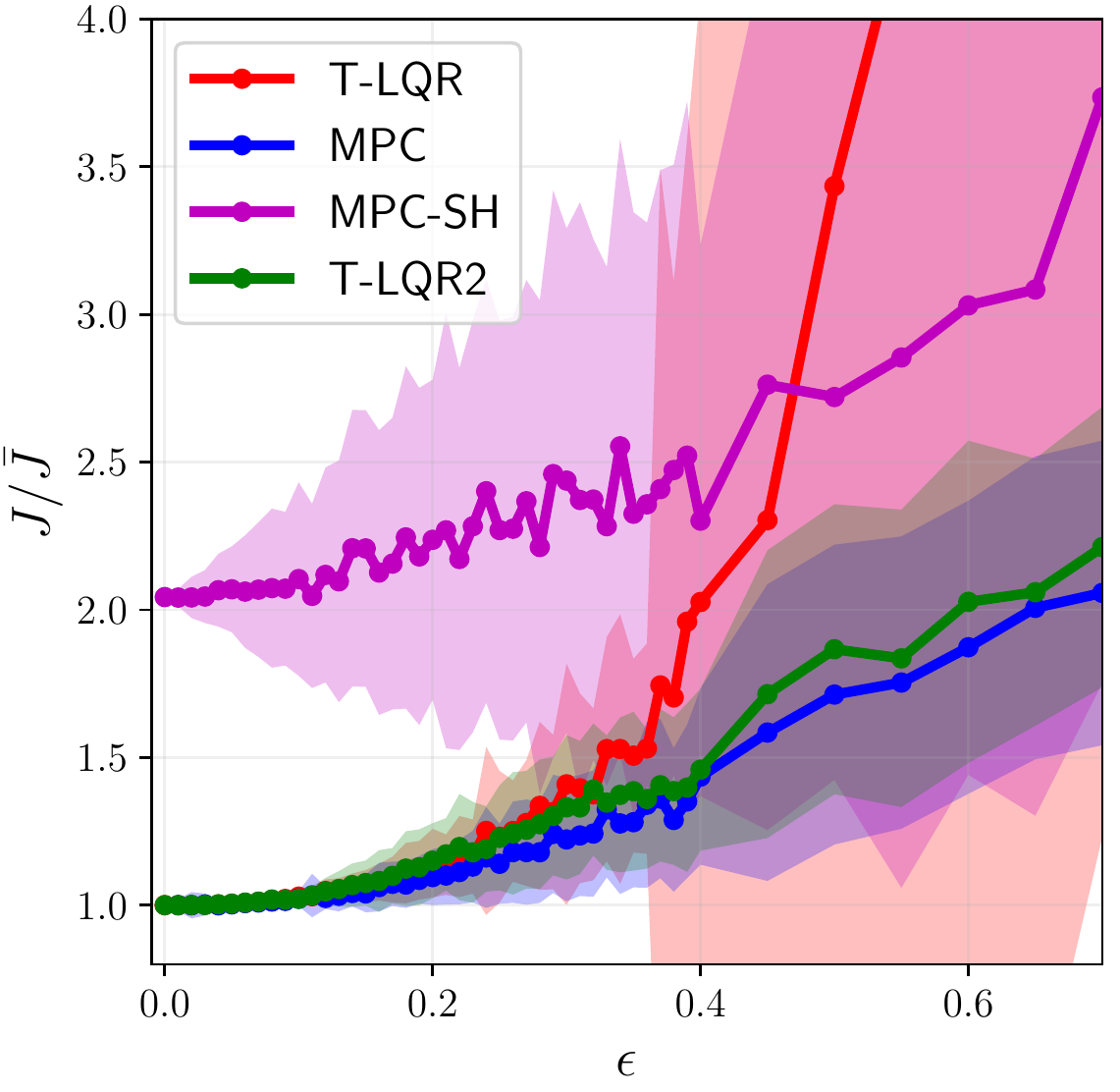}
        \label{quad_full}}
    \subfloat[Enhanced detail $0 \leq \epsilon \leq 0.4$]{
        \includegraphics[width=.30\textwidth]{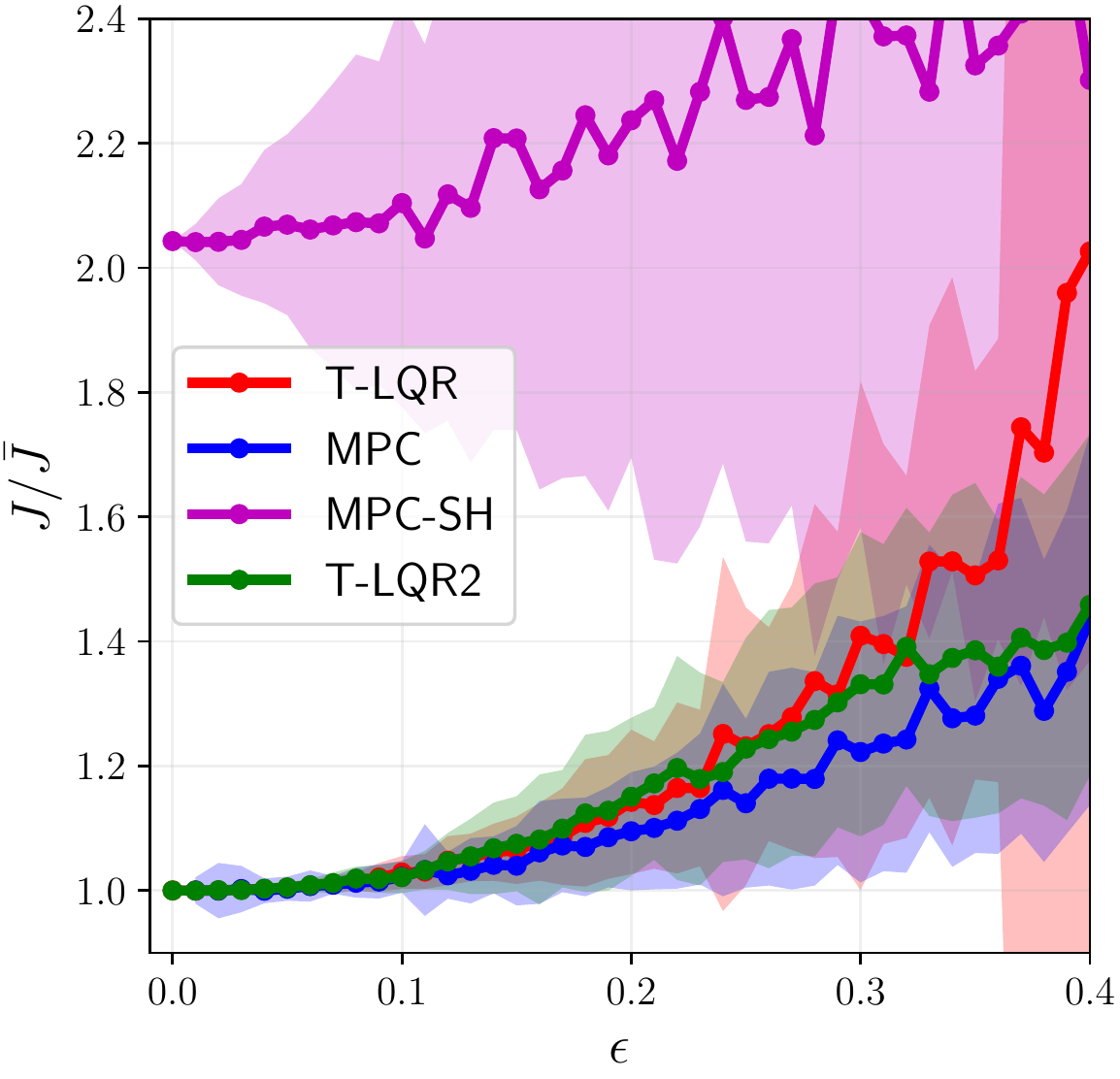}
        \label{quad_low}}
    \subfloat[Replanning operations]{
    \includegraphics[width=.30\textwidth]{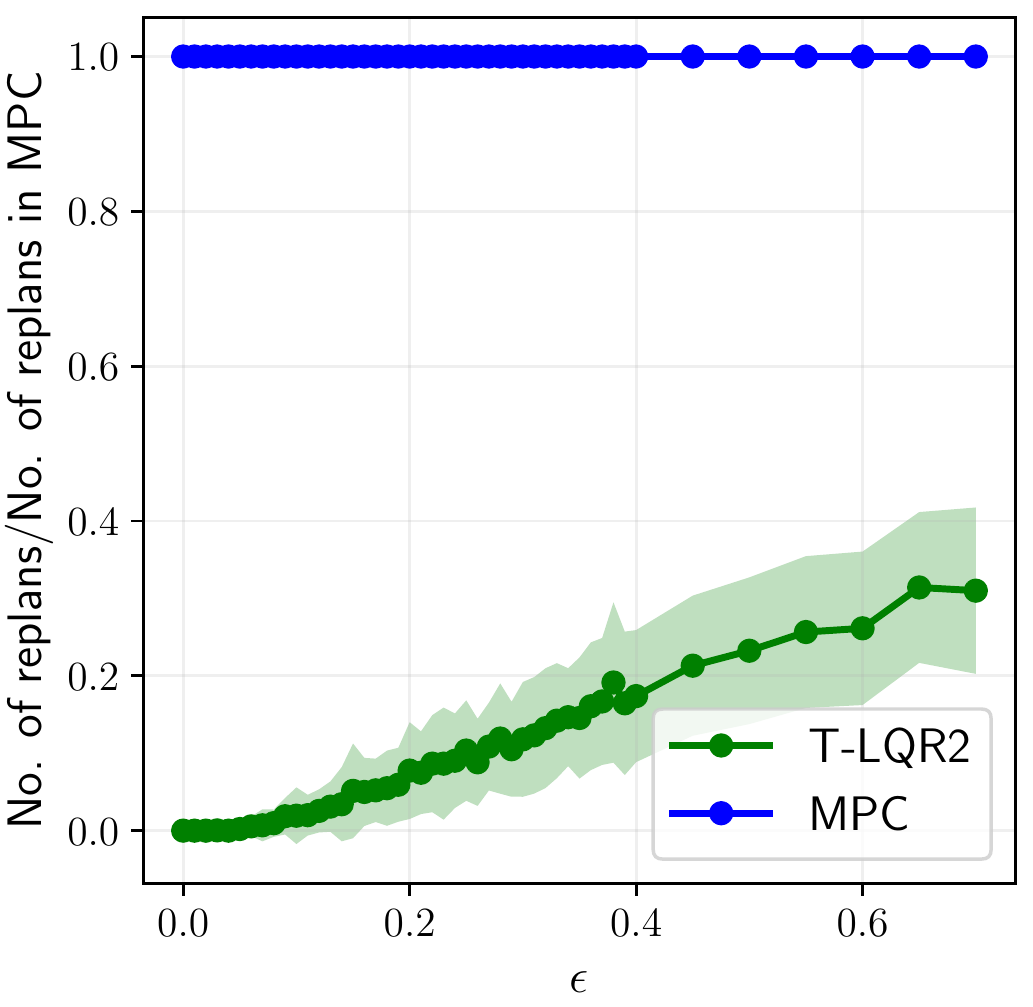}
    \label{fig:quad_replan}}
    \caption{Cost evolution of the different algorithms for a Quadrotor.}
    \label{fig:cost_quad}
\end{figure}

To evaluate in a 3D setting, we consider a quadrotor whose 12D state vector comprises of its position, orientation, linear and angular velocities - $(\V{x}_t,\mathbf{\theta}_t,\V{v}_t,\mathbf{\omega}_t)$. The model is described by
\begin{align*}
    \dot{\V{x}}_t &= \V{v}_t, & \dot{\V{v}}_t &= \V{g} + \frac{1}{m} R_{\mathbf{\theta}_t}\V{F}_t, \\
    \dot{\mathbf{\theta}}_t &= \M{W}_{\mathbf{\theta}}^{-1} \mathbf{\omega}_t, & \dot{\mathbf{\omega}}_t &= \M{I}^{-1} \mathbf{\tau}_t
\end{align*}
where, $\M{W}_{\mathbf{\theta}}$ is the transformation from the inertial to body frame and \M{I} is the inertia matrix. The model has thrust ($\V{F}_t$) and torques ($\mathbf{\tau_t}$) in its body fixed frame as the 4 control inputs. The results are shown in Figure \ref{fig:cost_quad}. Unlike a mobile robot which is stable even in high noise cases, a quadrotor is susceptible to failure or reach states from which no form of control can help it recover. So, the performance degrades earlier compared to the other two systems. But it can still be observed that T-LQR2 performs on a par with MPC in spite of the former replanning less than half the number of times compared to the latter. 

\subsection{Car-like robot in the presence of obstacles:} 
We show a case where the problem involves static obstacles in the environment. We assume the robot knows the map of the environment. The obstacles can be defined as ellipsoids. The ellipsoids can be represented with center $\V{o}^k \in \mathcal{R}^2$ and a positive definite matrix $\V{E}^k \in \mathcal{R}^{2\times 2}$. The obstacle penalty function for an agent whose position is $\V{p}_t$ in an environment with n obstacles is
\begin{equation}
    \Phi = \text{M}\ \sum_{k=1}^{n} \text{exp}(-[(\V{p}_t - \V{o}^k)^T \V{E}^k (\V{p}_t - \V{o}^k) - 1]), \nonumber
\end{equation} 
where M is a scaling factor.

\begin{figure}[!htbp]
    \subfloat[MPC]{
        \includegraphics[width=.225\textwidth]{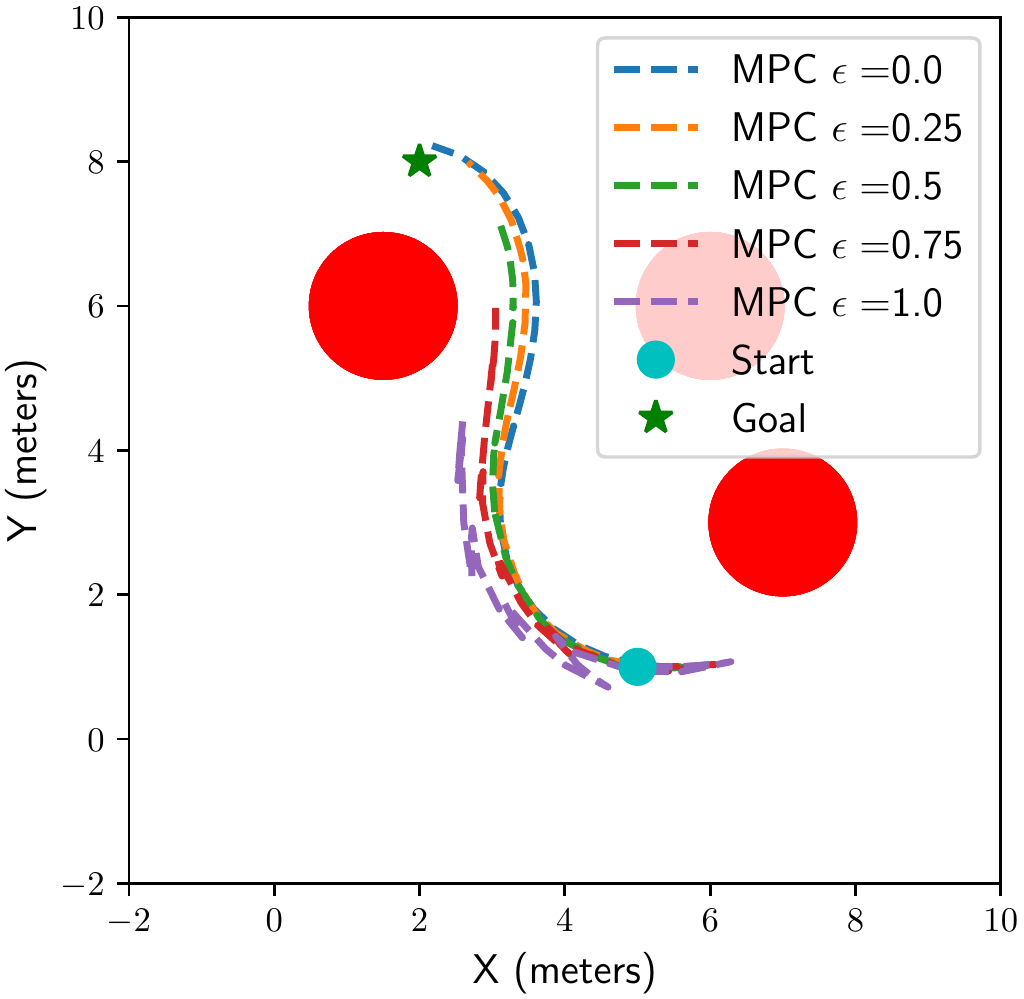}
        \label{test_cases_mpc_obs}}
    \subfloat[T-LQR2]{
        \includegraphics[width=.225\textwidth]{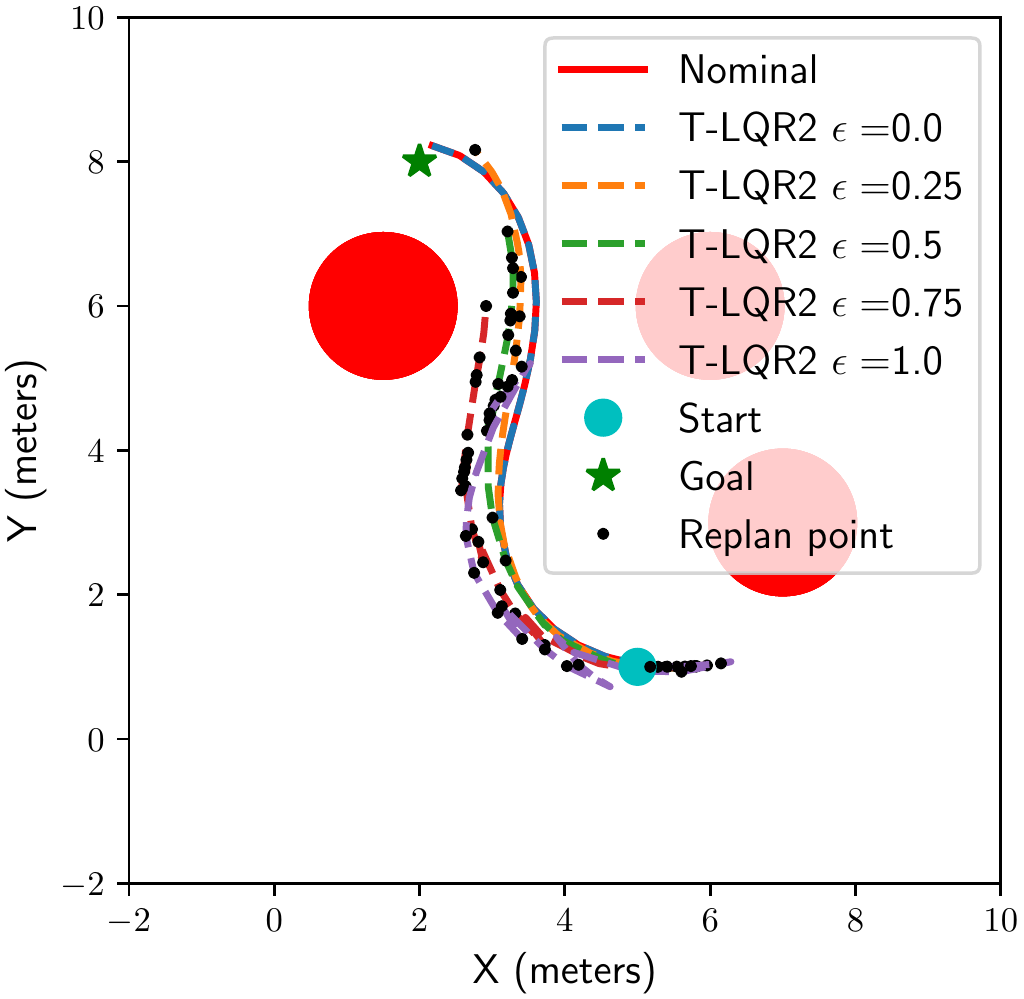}
        \label{test_cases_tlqr_replan_obs}}
    \subfloat[MPC-SH]{        \includegraphics[width=.225\textwidth]{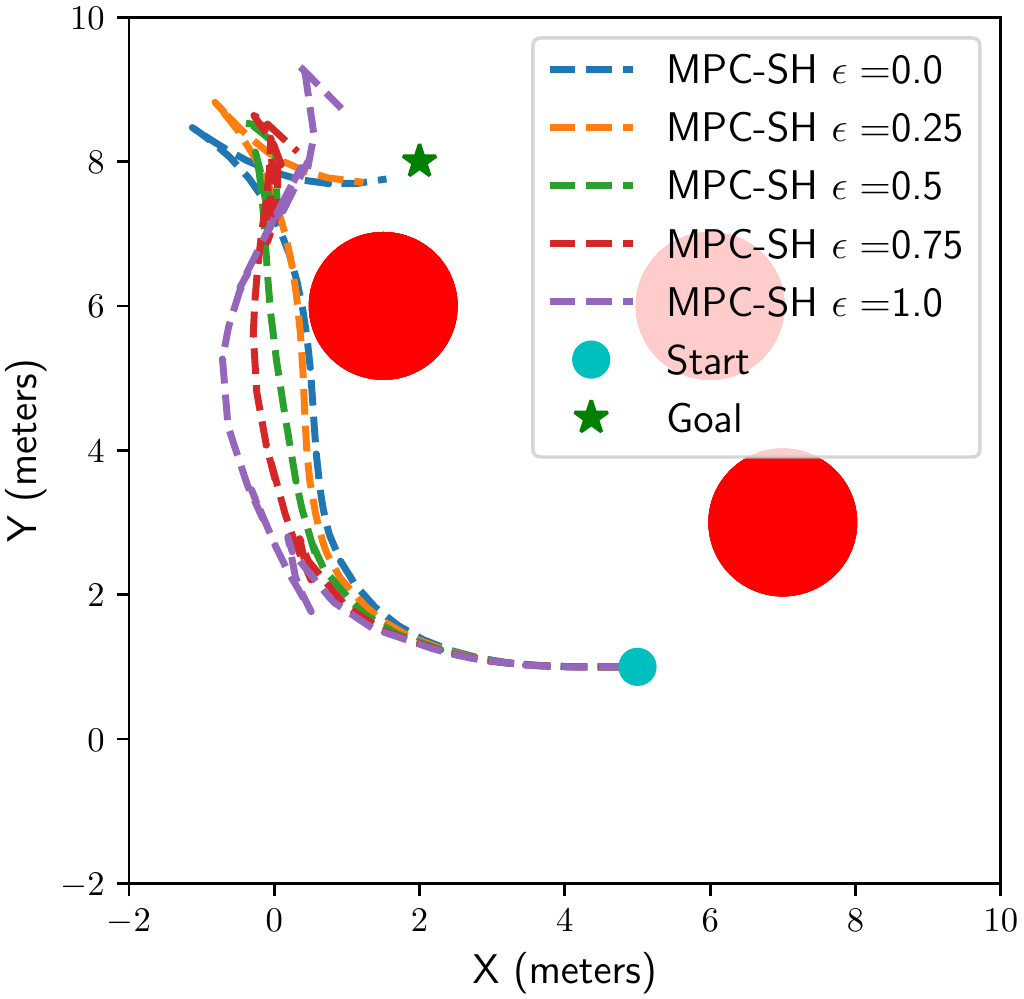}
        \label{test_cases_shmpc_obs}}
    \subfloat[T-LQR]{
        \includegraphics[width=.225\textwidth]{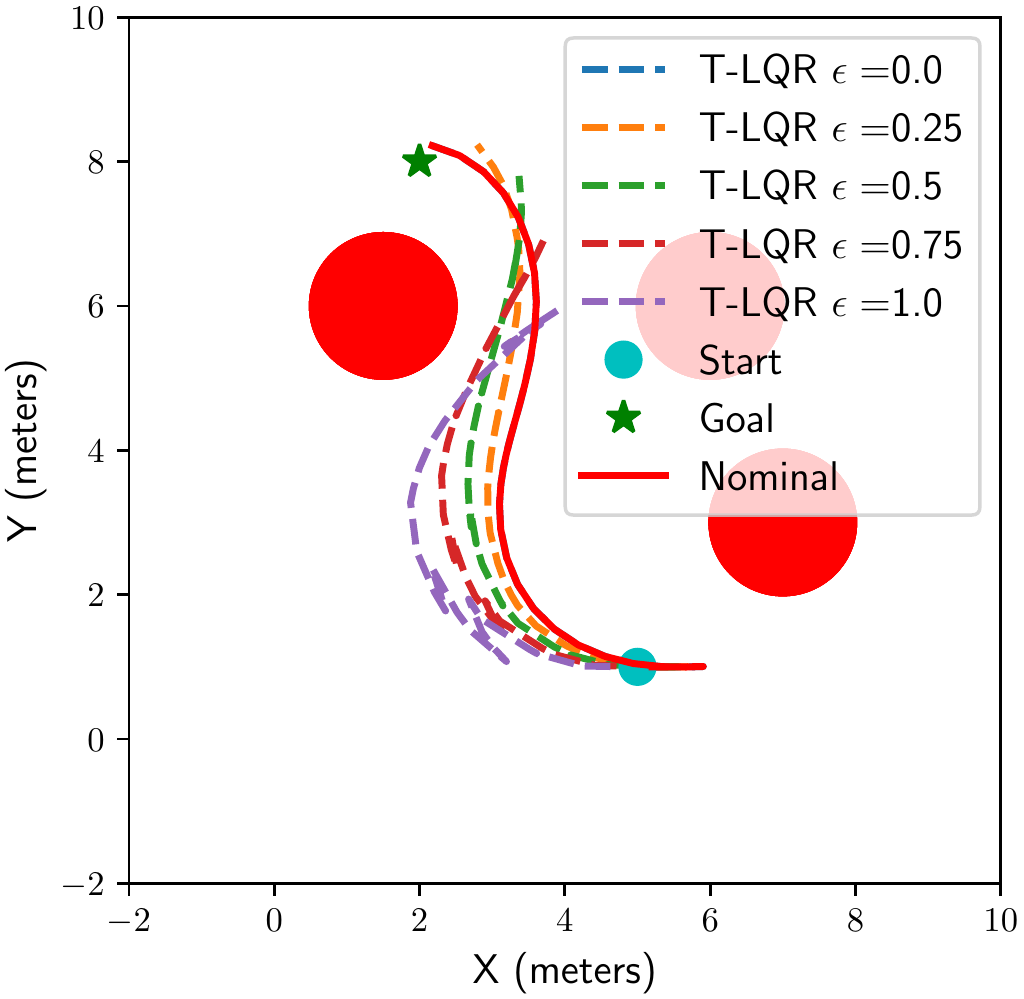}
        \label{test_cases_tlqr_obs}}
    \caption{The figure shows the paths taken by the robot using a particular algorithm and how they change as $\epsilon$ varies. A difference we see here is the paths taken by MPC-SH when compared to the others. Since it plans for a short horizon and hence greedy, it takes a different path unlike the others.}
    \label{fig:test_cases_obs}
\end{figure} 

\begin{figure}[h]
    \centering
    \subfloat[Full noise spectrum]{
        \includegraphics[width=.30\textwidth]{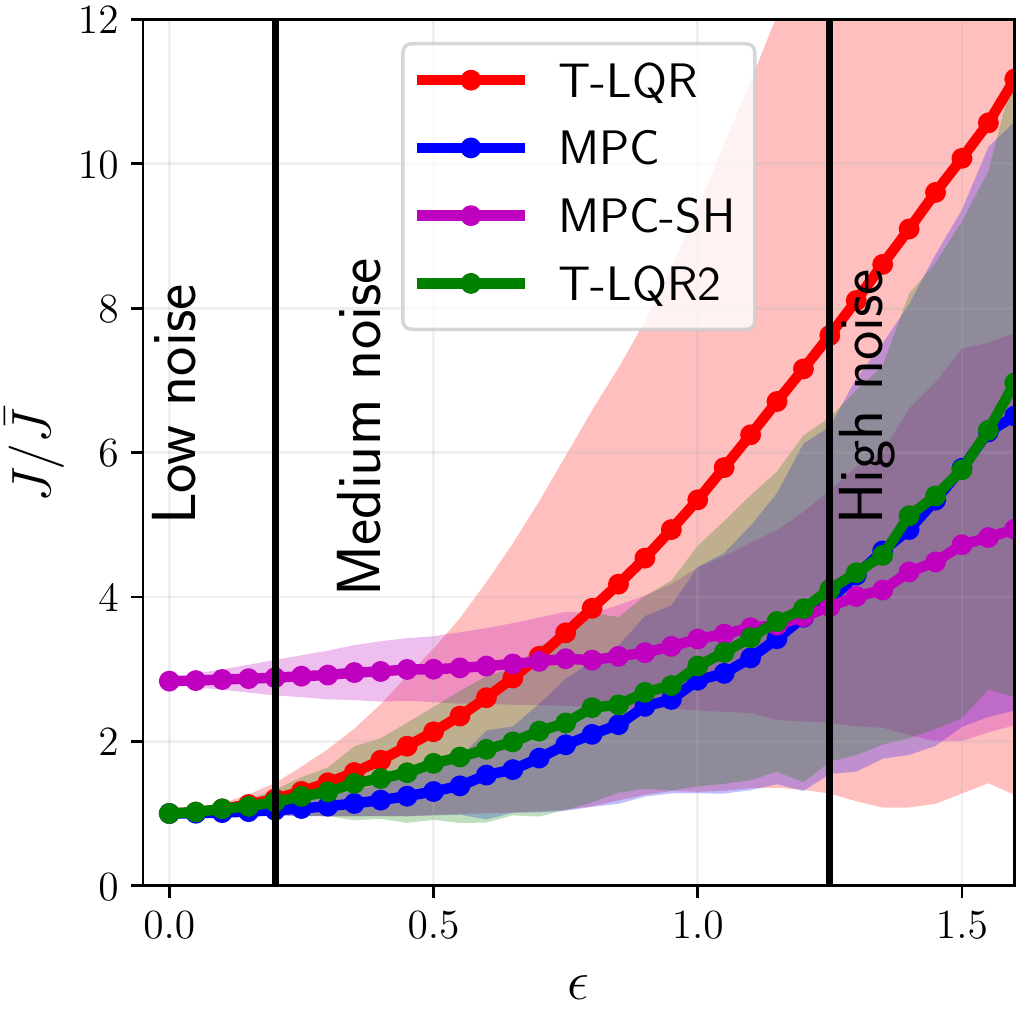}
        \label{car_obs_full}}
    \subfloat[Enhanced detail $0 \leq \epsilon \leq 0.4$]{
        \includegraphics[width=.30\textwidth]{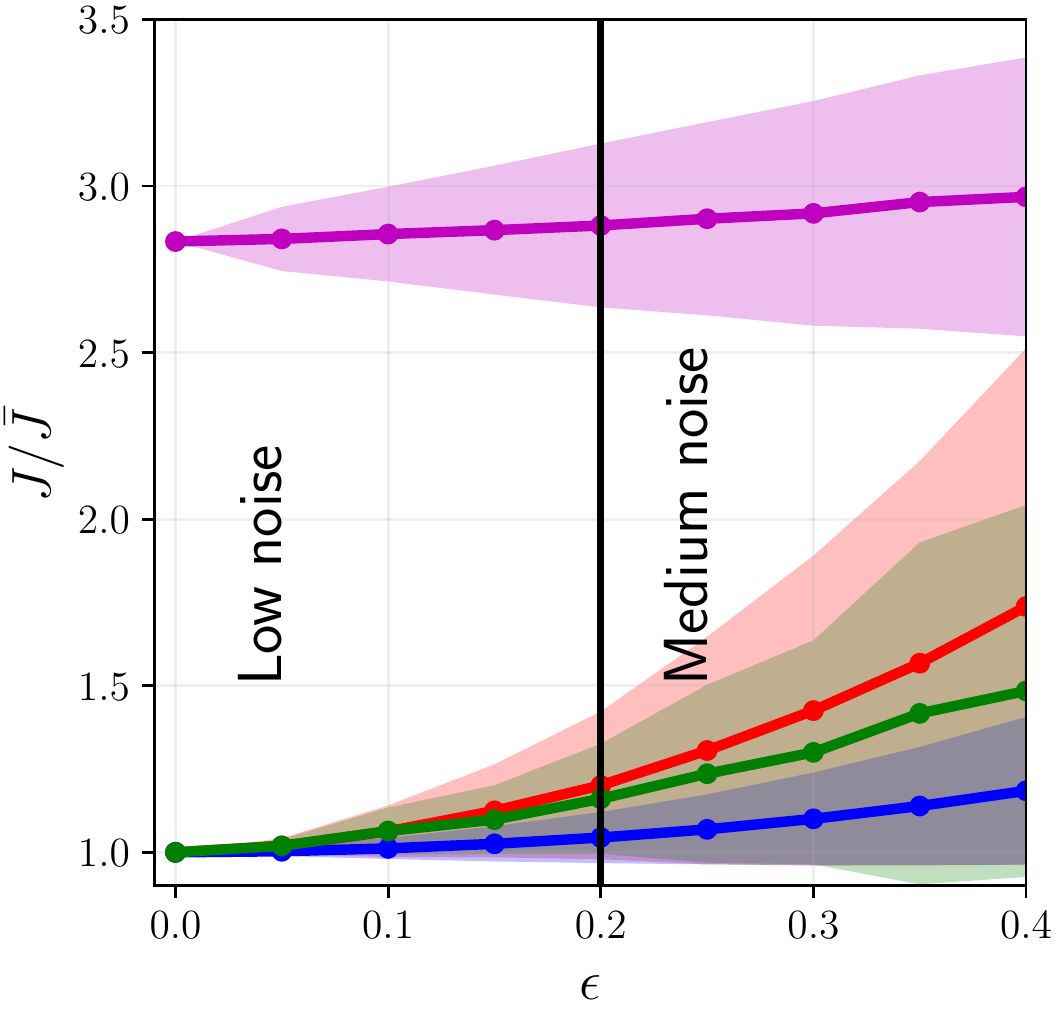}
        \label{car__obs_low}}
    \subfloat[Replanning operations]{
    \includegraphics[width=.30\textwidth]{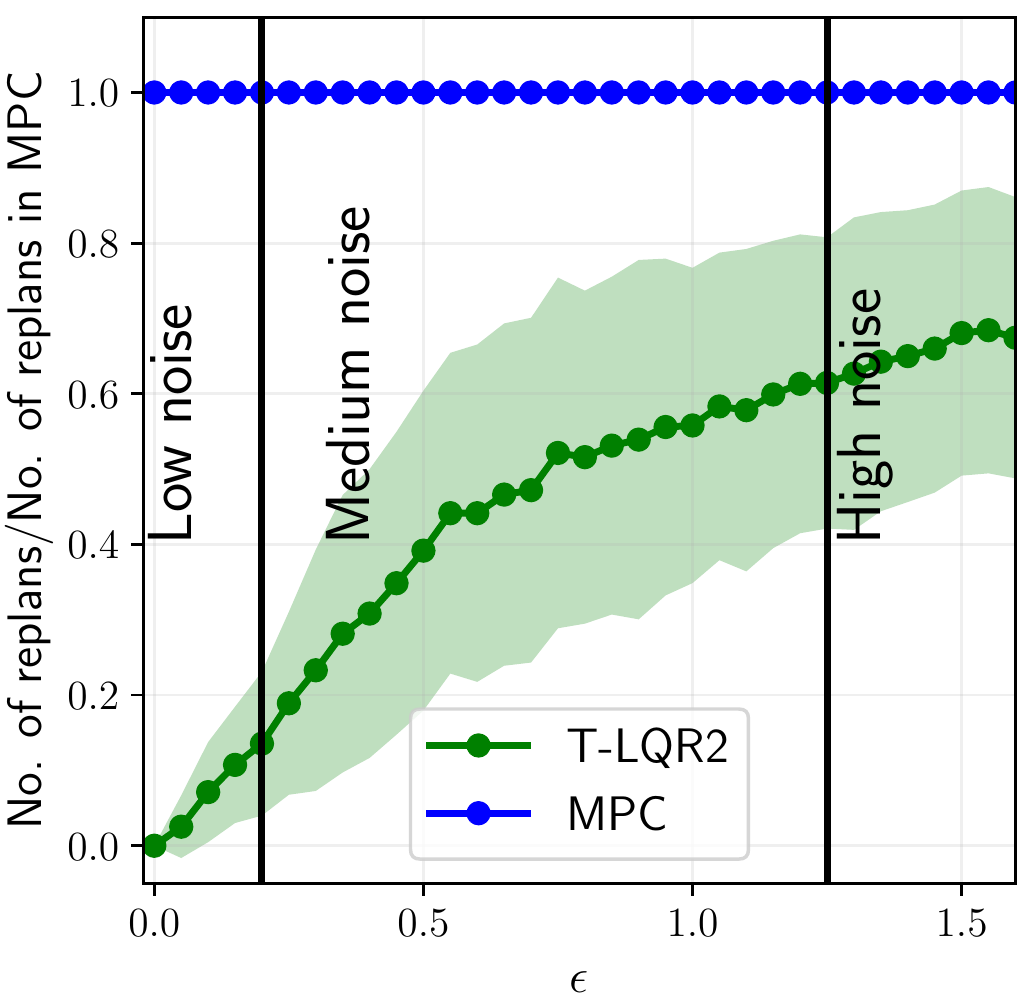}
    \label{fig:car_Obs_replan}}
    \caption{Cost evolution of the different algorithms for a car-like robot in an environment with obstacles.}
    \label{fig:cost_car_obs}
\end{figure}%
\begin{table}[h]
\centering
    \begin{tabular}{|m{1.5cm}|m{3cm}|m{3cm}|m{3cm}|}
        \hline
         & Car-like & Car with trailers & Quadrotor \\
         \hline
         $\V{x}_0$ & $[3,1,0,0]^T$ & $[0,0,\pi/3,0,0,0]^T$ & $[0,0,0,0,0,0,0,0,0,\newline 0,0,0]^T$\\
          \hline
         $\V{x}_f$ & $[3.5,7,m.pi/2,0]^T$ & $[2,2,0,0,0,0]^T$ & $[2,2,2,0,0,0,0,0,0,\newline 0,0,0]^T$\\
          \hline
         T, $\Delta$t & $35, 0.1$ & $40, 0.1$& $30, 0.1$ \\
          \hline
         $\M{W}^x$ & $diag(20,20,0,0)$ & $diag(10,10,1,1, \newline1,1)$ & $diag(10,10,10,1,1,\newline1,1,1,1,1,1,1)$\\
          \hline
         $\M{W}^u$ & $diag(20,200)$ & $diag(5,5)$ & $diag(5,10,10,10)$\\
          \hline
         $\M{W}^{x}_{f}$ & $10^3 diag(7,7,10,1)$ & $10^3diag(1,1,1,0.1,\newline 0.1,0.1)$ & $10^3 diag(1,1,1,1,1,\newline 1,1,1,1,1,1,1)$\\
          \hline
         Control bounds & $v_t = [-4,4],\newline \omega_t = [-\pi/12, \pi/12]$ & $ v_t = [-0.8,0.8], \newline \omega_t = [-\pi/6, \pi/6]$ & $u^{(1)}_t = [0,1.5], \newline u^{(i)}_t = [-0.05, 0.05]\newline i = 2,3,4$\\
          \hline
    \end{tabular}
    \caption{Parameters used in the single agent simulations.}
    \label{table:params1}
\end{table}

\begin{table}[h]
    \centering
    \begin{tabular}{|c|c|}
    \hline
         & Car-like \\
         \hline
         $\V{x}_0$ & Agent 1: $[3,1,\pi/2,0];$ \ Agent 2: $[5,1,0,0]$; \ Agent 3: $[6,8,0,0]$\\
         \hline
         $\V{x}_f$ & Agent 1: $[3.5,7,0,0];$ \ Agent 2: $[2,8,0,0]$; \ Agent 3: $[8,1.5,0,0]$\\
         \hline
         T, $\Delta$t & 35, 0.1 \\
         \hline
          $\M{W}^x$ & $diag(20,20,0,0)$ \\
          \hline
         $\M{W}^u$ & $diag(20,200)$ \\
          \hline
         $\M{W}^{x}_{f}$ & $10^3 diag(7,7,10,1)$ \\
         \hline
         Control bounds & $v_t = [-4,4],\ \omega_t = [-\pi/12, \pi/12]$ \\
         \hline
         
    \end{tabular}
    \caption{Parameters used in the multi-agent simulations.}
    \label{table:params3}
\end{table}
% \bibliographystyleS{abbrv}
% \bibliographyS{bib_files/MAP_refs1,bib_files/MohammadRafi,bib_files/naveed_references}

\end{document}